\documentclass[11pt,reqno]{amsart}
\usepackage{cases}
\usepackage{amscd}
\usepackage{amsfonts}
\usepackage{amssymb}
\usepackage{amsmath}
\usepackage{graphicx}
\usepackage{epstopdf}
\usepackage{subfigure}
\usepackage{amsthm}
\usepackage{stmaryrd}

\theoremstyle{remark}
\newtheorem{example}{\textbf{Example}}[section]
\numberwithin{equation}{section}

\usepackage{color}
\usepackage{datetime}

\makeatletter
\newcommand\figcaption{\def\@captype{figure}\caption}
\newcommand\tabcaption{\def\@captype{table}\caption}
\makeatother

\usepackage{geometry}
\usepackage{algorithm}
\usepackage{multirow}

\usepackage{hyperref}
\usepackage{cleveref}
\hypersetup{
  bookmarksnumbered = true,
  bookmarksopen=false,
  pdfborder=0 0 0,         
  pdffitwindow=true,      
  pdfnewwindow=true, 
  colorlinks=true,           
  linkcolor=blue,            
  citecolor=magenta,    
  filecolor=magenta,     
  urlcolor=cyan              
}

\def\bq{\begin{equation}}
\def\eq{\end{equation}}
\def\bqs{\begin{equation*}}
\def\eqs{\end{equation*}}
\def\bsqs{\begin{subequations}}
\def\esqs{\end{subequations}}
\def\ba{\begin{aligned}}
\def\ea{\end{aligned}}
\def\br{\begin{eqnarray}}
\def\er{\end{eqnarray}}
\def\brr{\bq\begin{array}{rlll}}
\def\err{\end{array}\eq}

\def\pmb#1{\mbox{\boldmath $#1$}}\def\text#1{\hbox{#1}}
\newtheorem{thm}{Theorem}[section]
\newtheorem{lem}{Lemma}[section]

\newtheorem{rem}{Remark}[section]

\newcommand{\bsub}{\begin{subequations}}
\newcommand{\esub}{\end{subequations}$\!$}

\newcommand{\maT}{\mathcal T}

\newcommand{\be}{\begin{eqnarray}}
\newcommand{\ee}{\end{eqnarray}}
\newcommand{\ben}{\begin{eqnarray*}}
\newcommand{\een}{\end{eqnarray*}}

\title[IEQ-FEMs]{Unconditionally energy stable IEQ-FEMs for the Cahn-Hilliard equation and Allen-Cahn equation}

\author[Y.~Chen, H.~Liu, N.~Yi, P. Yin]{Yaoyao Chen$^{\P}$,  Hailiang Liu$^\dag$, Nianyu Yi$^\ddag$, Peimeng Yin$^{\S,*}$}

\address{$\P$ School of Mathematics and Statistics, Anhui Normal University, Wuhu, Anhui 241000, PR China} \email{cyy1012xtu@126.com}
\address{$^\dag$ Department of Mathematics, Iowa State University, Ames, IA 50011, USA.} \email{hliu@iastate.edu}
\address{$\ddag$ Hunan Key Laboratory for Computation and Simulation in Science and Engineering, School of Mathematics and Computational Science, Xiangtan University, Xiangtan 411105, Hunan, P.R.China} \email{yinianyu@xtu.edu.cn}
\address{$^\S$ Department of Mathematical Sciences, The University of Texas at El Paso, El Paso, TX 79968, USA}
\email{pyin@utep.edu}

\begin{document}
\begin{abstract}
In this paper, we present several unconditionally energy-stable invariant energy quadratization (IEQ) finite element methods (FEMs)
with linear, first- and second-order
accuracy for solving both the Cahn-Hilliard equation and the Allen-Cahn equation.  For time discretization, we compare three distinct IEQ-FEM schemes that position the intermediate function introduced by the IEQ approach in different function spaces: finite element space, continuous function space, or a combination of these spaces. Rigorous proofs establishing the existence and uniqueness of the numerical solution, along with analyses of energy dissipation for both equations and mass conservation for the Cahn-Hilliard equation, are provided. The proposed schemes' accuracy, efficiency, and solution properties are demonstrated through numerical experiments.
\end{abstract}

\maketitle

\bigskip


\section{Introduction}
We focus on the application of invariant energy quadratization (IEQ) finite element methods (FEMs) to gradient flows, specifically addressing the Cahn-Hilliard (CH) equation \cite{CH58},
\begin{equation}\label{CH}
\begin{aligned}
 u_t  = \nabla \cdot (M(u) \nabla w) \qquad &{\it in}\ \Omega\ \times (0, T],\\ 
 w = - \epsilon^2 \Delta u + F'(u)\qquad &{\it in}\ \Omega\ \times (0, T],\\ 
\nabla u\cdot \pmb{n}=0,\quad \nabla w\cdot\pmb{n}=0 \qquad &{\it on}\ \partial\Omega\times [0, T],\\
u=u_{0}\qquad &{\it on}\  \Omega\times \{t=0\},
\end{aligned}
\end{equation}
and the Allen-Cahn (AC) equation \cite{AC75},
\begin{equation}\label{AC}
\begin{aligned}
 u_t = \epsilon^2 \Delta u - F'(u) \qquad &{\it in}\ \Omega\ \times (0, T],\\ 
\nabla u\cdot \pmb{n}=0 \qquad &{\it on}\ \partial\Omega\times [0, T],\\
u=u_{0}\qquad &{\it on}\  \Omega\times \{t=0\},
\end{aligned}
\end{equation}
where $\Omega \subseteq \mathbf{R}^d(d=2, 3)$ is a bounded domain, 
$\epsilon$ is a positive parameter, $M(u)\geq 0$ represents the mobility function, $u_0(x)$ is the initial data, and $\pmb{n}$ refers to the unit outward normal to the boundary $\partial \Omega$, $F(u)$ denotes the nonlinear bulk potential with the two well-known potentials being the double well potential \cite{FP2003}
\be\label{dwp}
F(u)=\frac{1}{4}(u^2-1)^2,
\ee
and the logarithmic Flory–Huggins potential \cite{BE92, CH58}
\be\label{fhp}
F(u) = \frac{\theta}{2} \left( u\ln u + (1-u) \ln (1-u) \right) + \frac{\theta_c}{2} u(1-u), 
\ee
where $\theta, \theta_c>0$ are physical parameters.

We consider the total free energy
\bq\label{ch_engdisO}
\mathcal{E}(u) = \int_{\Omega} \left( \frac{\epsilon^2}{2}|\nabla u|^2 + F(u) \right)dx.
\eq
Let $\frac{\delta}{\delta u}\mathcal{E}(u)$ be the variational derivative of the total free energy \eqref{ch_engdisO}.
Then the CH equation (\ref{CH}) is endowed with a $H^{-1}$-gradient flow structure $u_t= \nabla \cdot(M(u)\frac{\delta}{\delta u} \mathcal{E}(u))$, guided by the energy dissipation law
\begin{align}\label{di}
\frac{d}{dt} \mathcal{E}(u) = - \int_\Omega M(u) |\nabla w|^2 \leq 0;
\end{align}
and the AC equation is endowed with the $L^2$-gradient flow $u_t = -\frac{\delta}{\delta u} \mathcal{E}(u)$, guided by  the energy dissipation law
\begin{align}\label{di+}
\frac{d}{dt} \mathcal{E}(u) = - \int_\Omega |\nabla w|^2 \leq 0.
\end{align}

The gradient flows, represented by the CH problem \eqref{CH} and the AC problem \eqref{AC}, are inherently nonlinear, and making analytical solutions a challenging endeavor. Given the significance of steady states in gradient flows, a particular focus is placed on their study. Consequently, the development of precise, efficient, and energy-stable algorithms becomes imperative, especially when aiming for accurate simulations over extended periods. Ensuring adherence to energy dissipation laws, as expressed in (\ref{di}) or \eqref{di+}, emerges as a crucial consideration in formulating various numerical schemes \cite{BE92, CWWW19, Ey98, F01, FKU17, SXY19, SY10, S04, YZ19}.  This adherence is pivotal for guaranteeing the accuracy of long-time simulations.

The CH equation and the AC equation are semilinear parabolic equations with a small parameter $\epsilon$.  Various space discretization methods have been proposed for numerically solving these equations, including finite difference methods \cite{CFWW2019}, spectral methods \cite{DN1991}, finite element methods \cite{FP2003,FP2004}, and discontinuous Galerkin methods \cite{FK07, LY22, Y19, WPY22}.
In this paper, our emphasis is on the finite element method, a high-order space discretization technique rooted in the weak formulation of the gradient flows \eqref{CH} and \eqref{AC}. This method employs a finite element space composed of continuous piecewise polynomials for both trial and test spaces \cite{Ciarlet74}. Another noteworthy high-order numerical approach is the discontinuous Galerkin (DG) method \cite{LY09, LY10, LY18}, wherein the discontinuous Galerkin finite element space is characterized by completely discontinuous piecewise polynomials.

An additional significant challenge in simulating gradient flows \eqref{CH} and \eqref{AC} lies in effectively managing the nonlinear term during time discretization. Various methodologies have been proposed in the literature, such as the convex splitting approach \cite{CWWW19, Ey98, WZZ14}. However, this method involves iterative techniques. Another approach is the stabilization method \cite{SY10, XT06}, where 
selecting a stabilization constant is crucial for ensuring scheme stability. In recent developments, two noteworthy approaches have emerged for tackling gradient flows: the invariant energy quadratization (IEQ) approach \cite{Y16, ZWY17} and the scalar auxiliary variable (SAV) approach \cite{SXY18}. Importantly, both methods are linear and do not use iterative techniques. 

The SAV approach has been explored in conjunction with high-order spatial discretization methods, such as the finite element method \cite{CMS20, LW22} and the discontinuous Galerkin (DG) method \cite{LY23, YX22}. In contrast, research on the IEQ approximation with high-order spatial discretization has predominantly focused on DG methods \cite{LY19, LY21, LY22, WPY22}. Remarkably, there is limited research on the exploration of IEQ-FEM for gradient flows.

Based on the idea in \cite{Y16, ZWY17}, the IEQ approach for \eqref{CH} and \eqref{AC} is to introduce an intermediate function $U=\sqrt{F(u)+B}$, where $B \geq 0$ is a constant ensuring the validity of the new variable $U$, and $F(u)$ represents the nonlinear bulk potential in the gradient flows \eqref{CH} and \eqref{AC}. Consequently, the nonlinear term $F'(u)$ in the model equations is substituted with
\begin{equation*}
F'(u) = H(u)U, \quad U_t = \frac{1}{2}H(u)u_t,
\end{equation*}
where $H(u):=F'(u)/\sqrt{F(u)+B}$.

Motivated by various techniques employed for the intermediate function in DG schemes \cite{LY19, LY21, LY22, Y19, WPY22}, when coupled with the IEQ approach, this study delves into and compares three distinct unconditionally energy-stable IEQ-FEMs. 
These are based on different techniques in approximating the intermediate function $U$. 
Specifically, the nonlinear term $F'(u_h^{n+1})$ in the first order approximation is replaced by $H(u_h^n)\Phi$, where $u_h^n \in V_h$ is the finite element approximation of $u$ at time $t^n$ in the finite element space $V_h$, and $\Phi$ is an approximation of $U$ at time $t^{n+1}$ in certain  subspaces of $C^0(\Omega)$.

{ Three of our main IEQ-FEMs can be stated non-technically as follows:} 
\begin{itemize} 
\item The first method (Method 1) directly incorporates the approximation of the intermediate function $U$ at $t^{n+1}$ into the finite element space $V_h \subset C^0(\Omega)$ and treats it on par with other variables. The resulting fully discrete IEQ-FEM scheme implicitly updates the new variable, leading to an augmented linear system. While this does not amplify the computational cost in the DG setting due to an element-by-element matrix inverse used to eliminate the new variable, it does increase costs in the FEM setting due to the global matrix inverse, specifically, the inversion of the mass matrix.

\item The second method (Method 2) situates the intermediate function in continuous function space $C^0(\Omega)$ rather than the finite element space $V_h$. This approach facilitates independent and pointwise computation of the intermediate function, resulting in a computationally efficient method. Since the intermediate function lies outside the finite element space, the computed energy cannot be directly observed numerically;   
it can only be approximated within the finite element space. Consequently, the approximated energy may not consistently demonstrate energy decay, but conditional energy decay is achievable.  Notably, this condition performs well in the FEM setting, but is less effective in the DG setting.

\item The final method (Method 3), drawing inspiration from IEQ-DG schemes for the Swift-Hohenberg equation \cite{LY19, LY22} and the CH equation \cite{LY21, Y19}, first computes the intermediate function in the continuous function space $C^0$ and subsequently projects it onto the finite element space $V_h$. Despite the additional projection step compared to Method 2, this approach yields numerical solutions that are unconditionally energy-stable. It is noteworthy that this method proves effective in both the FEM and DG settings.
\end{itemize}

The paper is organized as follows. In Section \ref{secCH}, we develop linearly unconditionally energy stable fully discrete IEQ-FEM schemes for numerically solving the CH equation. The corresponding existence and uniqueness of the numerical solutions are rigorously proved as well as the property of energy dissipation. In Section \ref{secAC}, we extend the proposed methods for the AC equation. In Section \ref{secNumer}, several numerical examples demonstrate the accuracy, efficiency, and stability of the proposed schemes for CH and AC equations. We present concluding remarks in Section \ref{secConsi}. Second order IEQ-FEMs for the CH equation and the AC equation are presented in \Cref{CH2nd} and \Cref{AC2nd}, respectively.  

\section{IEQ-FEMs for the CH equation}\label{secCH}

In this section, we present the IEQ-FEMs for the CH equation (\ref{CH}).
Let $\maT$ be a quasi-uniform triangulation of $\Omega=\cup_{i=1}^N T_i$ with $N$ being the total number of triangles $T_i$, and $V_h\subset H^1(\Omega)$ be the $C^0$ Lagrange finite element space associated with $\maT$,
\be\label{eqn.fems}
V_h:=\{v\in C^0(\Omega) \cap H^1(\Omega): \ v|_{T_i}\in P_k, \ \forall T_i \in \maT\},
\ee
where $P_k$ is the space of polynomials of degree no more than $k$. 
Now we are prepared for the semi-discrete FEM scheme for the CH equation.

\subsection{The semi-discrete FEM for the CH equation}

The semi-discrete finite element scheme for the CH problem (\ref{CH}) is to find $(u_h, w_h) \in V_h \times V_h$ such that 
\begin{subequations}\label{CHsemifem}
\begin{align}
(u_{ht}, \phi) = & -A(M(u_h); w_h,\phi) ,\qquad  \phi \in V_h,\\
(w_h, \psi) = & A(\epsilon^2; u_h,\psi)+\left(F'(u_h),\psi \right), \qquad  \psi \in V_h,
\end{align}
\end{subequations}
where the bilinear functional 
\bq\label{bilinearfunc}
\ba
A(a(x);q,v) = & \int_\Omega a(x) \nabla q \cdot \nabla v dx
\ea
\eq
satisfies 
\be\label{Ablidef}
A(a(x);v,v) \geq \inf_{x\in \Omega} a(x) |v|^2_{H^1(\Omega)}, \quad \forall v \in V_h.
\ee
The initial data for $u_h$ is also chosen as 
\be\label{u0init}
u_h(x,0)=\Pi u_0(x), 
\ee
here and in what follows the operator $\Pi$ denotes the $L^2$ projection, i.e.,
\be
\int_\Omega (\Pi u_0(x)-u_0(x))\phi dx =0, \quad \forall \phi \in V_h.
\ee
The discrete  free energy is represented by
\be\label{freeenergy}
E(u_h) = \frac{1}{2}A(\epsilon^2; u_h, u_h) + \int_\Omega F(u_h)dx.
\ee

\begin{lem}\label{CHsemiLem}
The semi-discrete finite element scheme (\ref{CHsemifem}) conserves the total mass
\be\label{mass0}
\frac{d}{dt}\int_{\Omega }u_{h} dx = 0,
\ee
and the solution satisfies the energy dissipation law
\be\label{energy0}
\frac{d}{dt}E(u_h) = -A(M(u_h);w_h,w_h) \leq 0, \text{for all}\; t>0.
\ee
\end{lem}
\begin{proof}
Choosing  $\phi=1$ in (\ref{CHsemifem}), we obtain the total mass conservation (\ref{mass0}). 
By setting $\phi=w_h$ in (\ref{CHsemifem}a) and $\psi=u_{ht}$ in (\ref{CHsemifem}b), the resulting sum leads to (\ref{energy0}).  
\end{proof}

\subsection{The IEQ reformulation}\label{sec:IEQ}
For time discretization, we mimic the IEQ approach introduced in \cite{Y16}. Such time discretization transforms $F(u_h)$ 
in the energy functional into a quadratic form through an intermediate function, 
\bq\label{Uauxi}
U=\sqrt{F(u_h)+B},
\eq
for some constant $B$ such that $F(u_h)+B>0$.
Consequently, the total free energy (\ref{freeenergy}) can be reformulated as:
\bq\label{eng}
E(u_h, U)=\frac{1}{2}A(\epsilon^2; u_h,u_h)+\int_{\Omega}U^2dx -B|\Omega| =E(u_h).
\eq
By utilizing the intermediate function $U$, the nonlinear term $F'(u_h)$ can be expressed as
\be\label{Nonlinear}
F'(u_h) = H(u_h)U,
\ee
where
\begin{align}\label{ch_hw}
 H(w)= \frac{F'(w)}{\sqrt{F(w)+B}}.
\end{align}
This allows us to update $U$ by using  
\be\label{uode}
U_{t} = \frac{1}{2}H(u_h) u_{ht},
\ee
subject to the initial data
\be\label{CHinit+}
U(x, 0)=\sqrt{F(u_0(x))+B}.
\ee

Substituting $F'(u_h)$ in the semi-discrete FEM scheme \eqref{CHsemifem} with \eqref{Nonlinear}, along with \eqref{uode}, will lead to an augmented semi-discrete  formulation. A crucial consideration involves  determining a function space to approximate the intermediate function $U$. This choice will significantly influence the approximation of $F'(u_h)$ in \eqref{CHsemifem}. Our goal is to ensure that the numerical approximation aligns with the properties of the CH equation while maintaining the computational efficiency.

In the upcoming subsection, we will address this consideration by discretizing the augmented scheme  
using different function spaces for the intermediate function $U$. We will analyze the performance of the corresponding fully discrete IEQ-FEM schemes.

\subsection{First order fully discrete BDF1-IEQ-FEMs}\label{BDF1-IEQ}
In this section, we introduce three distinct fully discrete IEQ-FEMs.
In all these schemes, we explore first order backward differentiation (BDF1 for short). 
For $n \geq 0$, let $u_h^n \in V_h$, $w_h^n \in V_h$ represent $u(x,t^n)$, $w(x,t^n)$, where $t^n=n\Delta t$ and $\Delta t>0$ is the time step. 

\subsubsection{Method 1 (BDF1-IEQ-FEM1)}
First, we consider $U_h(x,t) \in V_h$ in the finite element space to approximate $U$ based on the following relation 
\be\label{CHSemiFEM+}
(U_{ht}, \tau) = \frac{1}{2}(H(u_h) u_{ht},\tau) ,\qquad \forall  \tau \in V_h, 
\ee
subject to the initial data
\bq\label{CHinit}
\quad U_h(x, 0)=\Pi\sqrt{F(u_0(x))+B}.
\eq
More specifically, for $n\geq 0$, let $U_h^n \in V_h$ represent $U_h(x,t^n)$. 
Given $u^{n}_h, w_h^{n}, U_h^n \in V_h$, the first order fully discrete BDF1-IEQ-FEM1 scheme is to find  $(u^{n+1}_h, w_h^{n+1}, U_h^{n+1}) \in V_h  \times V_h \times V_h$ such that
\begin{subequations}\label{BDF1-IEQ-FEM1}
\begin{align}
\left(  \frac{u_h^{n+1} - u_h^n}{\Delta t}, \phi \right)= & -A(M(u_h^n); w_h^{n+1}, \phi), \\
(w_h^{n+1},\psi) = & A(\epsilon^2; u_h^{n+1}, \psi)+\left( H(u_h^n)U_h^{n+1}, \psi \right),\\
\left(\frac{U_h^{n+1} - U_h^n}{\Delta t}, \tau\right) = & \frac{1}{2}\left(H(u_h^n) \frac{u_h^{n+1} - u_h^n}{\Delta t}, \tau\right),
\end{align}
\end{subequations}
for any $ \phi, \psi, \tau \in V_h $. The initial data $(u_h^0, U_h^0)$ is given by \eqref{u0init} and \eqref{CHinit}.
It can be verified that the following result holds.
\begin{lem}
The first order fully discrete BDF1-IEQ-FEM1 scheme (\ref{BDF1-IEQ-FEM1}) admits a unique solution $(u^{n+1}_h, w_h^{n+1}, U_h^{n+1}) \in V_h  \times V_h \times V_h$ for any $\Delta t>0$. The solution $u_h^n$ satisfies the total mass conservation
\be\label{BDF1MassCons}
\int_\Omega u_h^ndx = \int_\Omega u_h^0 dx,
\ee
for any $n>0$. Additionally, it satisfies the energy dissipation law
\bqs
E(u_h^{n+1}, U_h^{n+1}) = E(u_h^{n}, U_h^{n}) - \|U_h^{n+1} - U_h^n\|^2 -\Delta t A(M(u_h^n); w_h^{n+1},  w_h^{n+1}) - \frac{1}{2} A(\epsilon^2; u_h^{n+1}-u_h^n, u_h^{n+1} - u_h^n),
\eqs
independent of the time step $\Delta t$.
\end{lem}

The linear system for the BDF1-IEQ-FEM1 scheme \eqref{BDF1-IEQ-FEM1} takes the form 
\begin{equation}\label{BDF1-IEQ-FEM1-mat}
\begin{aligned}
\left(
  \begin{array}{ccc}
    \frac{1}{\Delta t}G &  D_1 &  \mathbf{0} \\
    -\epsilon^{2}D & G & -G_1 \\
    -\frac{G_1}{2\Delta t} & \mathbf{0} & \frac{G}{\Delta t} \\
  \end{array}
\right)\left(
  \begin{array}{ccc}
    \vec u^{n+1} \\
    \vec w^{n+1}\\
    \vec U^{n+1} \\
  \end{array}
\right)=\left(\begin{array}{ccc}
    \vec g_1\\
    \vec g_2\\
    \vec g_3  \\
  \end{array}
\right).
\end{aligned}
\end{equation}
Here $\vec{u}^{n}$, $\vec{w}^{n}$, $\vec{U}^{n}$ represent the coefficient vectors of the finite element approximations $u_h^{n}, w_h^{n}, U_h^{n}$ (we have abused the notation ).  $D$ and $G$ are the stiffness matrix and the mass matrix, respectively. $D_1$ is the stiffness matrix with the weight $M(u_h^n)$, and $G_1$ is the mass matrix with the weight $H(u_h^n)$.
The vectors $\vec g_1$, $\vec g_2$, $\vec g_3$ are given as
$$
\vec g_1=\frac{1}{\Delta t} G \vec u^{n}, \ \ \vec g_2=\vec{0},\ \ \vec g_3=\frac{1}{\Delta t} \left(G  \vec U^{n} -\frac{1}{2}G_1 \vec  u^n \right).   
$$
Solving the linear system \eqref{BDF1-IEQ-FEM1-mat} at each step can be computationally expensive. To mitigate this, an alternative approach is to reduce the dimension of the linear system \eqref{BDF1-IEQ-FEM1-mat}. Specifically, we express the third row of \eqref{BDF1-IEQ-FEM1-mat} as
\be
\vec U^{n+1} = G^{-1}\left( \Delta t \vec g_3 + \frac{1}{2}G_1 \vec u^{n+1} \right).
\ee
Therefore, the matrix equation \eqref{BDF1-IEQ-FEM1-mat} reduces to
\begin{equation}\label{BDF1-IEQ-FEM1-mat1}
\begin{aligned}
\left(
  \begin{array}{cc}
    \frac{1}{\Delta t}G &  D_1 \\
    -\epsilon^{2}D-\frac{1}{2}G_1G^{-1}G_1   & G  \\
  \end{array}
\right)\left(
  \begin{array}{c}
    \vec u^{n+1} \\
    \vec w^{n+1}\\
  \end{array}
\right)=\left(\begin{array}{c}
    \vec g_1\\
    \vec g_2 + G_1 G^{-1}\Delta t \vec g_3 \\
  \end{array}
\right).
\end{aligned}
\end{equation}
To solve the reduced linear system \eqref{BDF1-IEQ-FEM1-mat1}, it is essential to evaluate the inverse of the mass matrix, $G^{-1}$. 
In the DG setting, the matrix $G$ adopts a diagonal block structure, and computing its inverse $G^{-1}$ does not incur additional computational complexity. However, in the finite element framework, the matrix $G$ lacks this local property, leading to a noticeable  
increase in computational cost when computing $G^{-1}$. This computational burden makes  the BDF1-IEQ-FEM1 scheme relatively expensive.

\subsubsection{Method 2 (BDF1-IEQ-FEM2)}
An alternative approach is to preserve the intermediate function $U(x,t)$ within the continuous function space $C^0(\Omega)$,  with its evolution governed by \eqref{uode}.
More specifically, for $n\geq 0$, let $U^n$ represent $U(x,t^n)$.
Given $u^{n}_h, w_h^{n} \in V_h$ and $U^n \in C^0(\Omega)$, the first order fully discrete BDF1-IEQ-FEM2 scheme is to find $(u^{n+1}_h, w_h^{n+1}) \in V_h  \times V_h $ and $U^{n+1} \in C^0(\Omega)$ such that
\begin{subequations}\label{BDF1-IEQ-FEM2}
\begin{align}
\left(  \frac{u_h^{n+1} - u_h^n}{\Delta t}, \phi \right)= & -A(M(u_h^n); w_h^{n+1}, \phi), \\
(w_h^{n+1},\psi) = & A(\epsilon^2; u_h^{n+1}, \psi)+\left( H(u_h^n)U^{n+1}, \psi \right),\\
\frac{U^{n+1} - U^n}{\Delta t} = & \frac{1}{2}H(u_h^n) \frac{u_h^{n+1} - u_h^n}{\Delta t},
\end{align}
\end{subequations}
for any $\phi, \psi \in V_h $. The initial data $(u_h^0, U^0)$ is given by \eqref{u0init} and \eqref{CHinit+}.

Reformulating (\ref{BDF1-IEQ-FEM2}c) gives
\bq\label{CHPwU}
U^{n+1} = U^n+\frac{1}{2}H(u_h^n) (u_h^{n+1} - u_h^n).
\eq
Plugging $U^{n+1}$ into (\ref{BDF1-IEQ-FEM2}b) yields a system expressed solely in terms of  $(u^{n+1}_h, w_h^{n+1})$,
\begin{subequations}\label{CHBDF1Lin}
\begin{align}
& \left(  \frac{u_h^{n+1} - u_h^n}{\Delta t}, \phi \right)=  -A(M(u_h^n); w_h^{n+1}, \phi), \\
& (w_h^{n+1},\psi) = A(\epsilon^2; u_h^{n+1}, \psi)+\left( H(u_h^n)U^{n}, \psi \right) + \frac{1}{2}((H(u_h^n))^2 u_h^{n+1}, \psi) - \frac{1}{2}((H(u_h^n))^2 u_h^{n}, \psi).
\end{align}
\end{subequations}
The first order fully discrete BDF1-IEQ-FEM2 scheme (\ref{BDF1-IEQ-FEM2}) is equivalent to the system formed by (\ref{CHBDF1Lin}) and (\ref{CHPwU}). Here one can first obtain $(u^{n+1}_h, w_h^{n+1})$ from (\ref{CHBDF1Lin}), and then obtain $U^{n+1}$ by (\ref{CHPwU}) or (\ref{BDF1-IEQ-FEM2}c).
Compared with the BDF1-IEQ-FEM1 scheme (\ref{BDF1-IEQ-FEM1}), the scheme (\ref{BDF1-IEQ-FEM2}) can avoid the need to solve a large coupled system or compute inverse matrices. Therefore, the BDF1-IEQ-FEM2 scheme is computationally more econmomical.

Similar to the BDF1-IEQ-FEM1 scheme, the following result holds.

\begin{thm}\label{CHBDF1thm}
The first-order fully discrete BDF1-IEQ-FEM2 scheme (\ref{BDF1-IEQ-FEM2}) produces a unique solution $(u_h^{n+1},w_h^{n+1}, U^{n+1})\in V_h\times V_h \times C^0(\Omega)$ for any $\Delta t>0$. The solution $u_h^n$ also satisfies total mass conservation \eqref{BDF1MassCons}
for any $n>0$. Additionally, it satisfies the energy dissipation law:
\bq\label{CHBDF1Stab}
\ba
E(u_h^{n+1}, U^{n+1}) = &  E(u_h^{n}, U^{n}) - \|U^{n+1} - U^n\|^2 \\
 & -\Delta t A(M(u_h^n); w_h^{n+1},  w_h^{n+1}) - \frac{1}{2} A(\epsilon^2; u_h^{n+1}-u_h^n, u_h^{n+1} - u_h^n) ,
\ea
\eq
independent of the time step $\Delta t$.
\end{thm}
\begin{proof}
We begin by establishing the existence and uniqueness of the solution. Since the scheme (\ref{BDF1-IEQ-FEM2}) is equivalent to the system (\ref{CHBDF1Lin}) along with (\ref{CHPwU}), they share the same existence and uniqueness of the solution. As (\ref{CHBDF1Lin}) represents a finite-dimensional linear system, the existence of the solution is equivalent to its uniqueness.
Now, assume that the linear system (\ref{CHBDF1Lin}) has two solutions, and denote their difference at $t^{n+1}$ by $(\tilde u^{n+1}_h, \tilde w^{n+1}_h)$. This leads to the following system:
\begin{subequations}\label{CHBDF1Diff}
\begin{align}
\left(  \frac{\tilde u_h^{n+1}}{\Delta t}, \phi \right)= & -A(M(u_h^n); \tilde w_h^{n+1}, \phi), \\
(\tilde w_h^{n+1},\psi) = & A(\epsilon^2; \tilde u_h^{n+1}, \psi) + \frac{1}{2}((H(u_h^n))^2 \tilde u_n^{n+1}, \psi).
\end{align}
\end{subequations}
Taking $\phi=\Delta t \tilde w_h^{n+1}$ and $\psi=\tilde u_h^{n+1}$ in (\ref{CHBDF1Diff}) yields
$$
A(\epsilon^2; \tilde u_h^{n+1}, \tilde u_h^{n+1}) + \frac{1}{2}((H(u_h^n))^2 \tilde u_n^{n+1}, \tilde u_n^{n+1})+\Delta t A(M(u_h^n); \tilde w_h^{n+1}, \tilde w_h^{n+1}) =0,
$$
which implies
$$
\epsilon^2|\tilde u_h^{n+1}|^2_{H^1(\Omega)}+ \frac{1}{2}\|H(u_h^n) \tilde u_n^{n+1}\|^2 + \Delta t \min_{x\in \Omega} M(u_h^n) |\tilde w_n^{n+1}|^2_{H^1(\Omega)} \leq 0.
$$
Consequently, we find that $\tilde u_h^{n+1}=\text{const}$ and $\tilde w_h^{n+1}=\text{const}$. From (\ref{CHBDF1Diff}a), we have
$$
\left(  \tilde u_h^{n+1}, \phi \right)=0, \quad \phi \in V_h,
$$
which implies $\tilde u_h^{n+1}=0$. By (\ref{CHBDF1Diff}b), this also leads to $\tilde w_h^{n+1}=0$.
These establish the existence and uniqueness of the solution for the system (\ref{CHBDF1Lin}). Since $U^{n+1}$ is uniquely determined by $u_h^{n+1}$ from (\ref{CHPwU}) or (\ref{BDF1-IEQ-FEM2}c), the scheme (\ref{BDF1-IEQ-FEM2}) admits a unique solution.

By setting $\phi=1$ in (\ref{BDF1-IEQ-FEM2}a), we obtain the mass conservation (\ref{BDF1MassCons}).

Finally, to demonstrate energy stability (\ref{CHBDF1Stab}), we select $\phi=\Delta t w_h^{n+1}$, $\psi=-(u_h^{n+1} - u_h^n)$ in (\ref{BDF1-IEQ-FEM2}ab). Taking the $L^2$ inner product of (\ref{BDF1-IEQ-FEM2}c) with $2\Delta t U^{n+1}$, the summation of (\ref{BDF1-IEQ-FEM2}) upon regrouping gives:
\bqs
\ba
 -\Delta t A(M(u_h^n); w_h^{n+1},  w_h^{n+1})= & A(\epsilon^2; u_h^{n+1}, u_h^{n+1} - u_h^n) + 2(U^{n+1} - U^n, U^{n+1})\\
= & \frac{1}{2}\left(A(\epsilon^2; u_h^{n+1}, u_h^{n+1})-A(\epsilon^2; u_h^{n}, u_h^{n}) + A(\epsilon^2; u_h^{n+1}-u_h^n, u_h^{n+1} - u_h^n)\right)\\
& + \left( \|U^{n+1}\|^2 - \|U^{n}\|^2 + \|U^{n+1} - U^n\|^2 \right), 
\ea
\eqs
which establishes the energy stability (\ref{BDF1-IEQ-FEM2}).
\end{proof}

\begin{rem}
The BDF1-IEQ-FEM2 scheme provides a computationally efficient solution  while preserving the energy dissipation law \eqref{CHBDF1Stab}. However, the exact computation of the energy $E(u_h^{n}, U^{n})$ is unattainable, given that $U^{n} \not \in V_h$. Instead, it can only be approximated by $E(u_h^{n}, \Pi U^{n})$, where $\Pi$ is an $L^2$ projection onto the finite element space $V_h$. Here, the approximated energy may not 
exhibit dissipation, as demonstrated in \Cref{cexm2-n} for the CH equation.
\end{rem}

\begin{rem}\label{Enerdis}
While the energy $E(u_h^{n}, U^{n})$ is unconditional  stable, the approximated energy $E(u_h^{n}, \Pi U^{n})$ may exhibit only conditional stability.  To ensure the energy dissipation laws to hold also for the approximated energy, a natural strategy is to minimize the error $|E(u_h^{n}, U^{n})-E(u_h^{n}, \Pi U^{n})|$. Let $\Delta E^n = |E(u_h^{n}, U^{n})-E(u_h^{n-1}, U^{n-1})|$.
For any $n\geq 1$, if the condition
\be\label{projenergy}
|E(u_h^{n}, U^{n})-E(u_h^{n}, \Pi U^{n})| \leq \frac{1}{2}\max\{ \Delta E^{n+1}, \Delta E^n \},
\ee
is met, then the approximated energy $E(u_h^{n}, \Pi U^{n})$ will satisfy an energy dissipation law.
Minimizing the approximated energy error is equivalent to minimizing the projection error $\|U^n-\Pi U^n\|$.
Various numerical techniques, such as using refined meshes and adjusting the constant $B$ in the definition \eqref{Uauxi} with a larger value (as recommended in Remark 4.1 in \cite{LY22}), can effectively reduce this projection error.
\end{rem}

\begin{rem}
In the FEM framework, $U^n \in H^1(\Omega)$, providing sufficient regularity to minimize the projection error $\|U^n-\Pi U^n\|$ using the techniques outlined in \Cref{Enerdis}. This ensures the
fulfillment of the condition \eqref{projenergy}.
In contrast, a similar IEQ-DG method was evaluated in the DG setting in \cite{LY22, Y19} for the CH equation. However,  numerical tests revealed that the DG solution does not adhere to the energy dissipation law. 
\end{rem}

\begin{rem}
In reference to the techniques outlined in \Cref{Enerdis}, it is crucial to assess their numerical performance. Refining meshes has proven effective in specific scenarios, as demonstrated in \Cref{cexmAC2-n}. However, adopting an appropriately large constant $B$ appears to be effective across a broad spectrum of cases, although it is important to note that the optimal constant may be problem-specific. 
\end{rem} 

In addition to the techniques discussed in \Cref{Enerdis}, we propose a method that is independent of the numerical setting to ensure  unconditional energy dissipation.

\subsubsection{Method 3 (BDF1-IEQ-FEM3)}
It is natural to explore whether we can adapt the BDF1-IEQ-FEM2 scheme to preserve its accuracy, achieve computational efficiency, and ensure the numerical satisfaction of the energy dissipation law with the computed energy $E(u_h^{n}, \Pi U^{n})$. Essentially, this entails transforming the intermediate function $U^n$ into a representation in $V_h \subset C^0(\Omega)$. Motivated by the techniques for the IEQ-DG method \cite{LY19, LY21, LY22, Y19}, we introduce the following modified IEQ-FEM scheme.

Consider $u^{n}_h, w_h^{n} \in V_h$ and $U^n \in C^0(\Omega)$, the modified first-order fully discrete scheme, referred to as the BDF1-IEQ-FEM3 scheme, seeks to find $(u^{n+1}_h, w_h^{n+1}) \in V_h \times V_h$ and $U^{n+1} \in C^0(\Omega)$ such that for $\forall \phi, \psi \in V_h $,
\begin{subequations}\label{BDF1-IEQ-FEM3}
\begin{align}
\left(  \frac{u_h^{n+1} - u_h^n}{\Delta t}, \phi \right)= & -A(M(u_h^n); w_h^{n+1}, \phi), \\
(w_h^{n+1},\psi) = & A(\epsilon^2; u_h^{n+1}, \psi)+\left( H(u_h^n)U^{n+1}, \psi \right),\\
U^n_h= & \Pi U^n, \label{BDF1-IEQ-FEM33}\\
\frac{U^{n+1} - U_h^n}{\Delta t} = & \frac{1}{2}H(u_h^n) \frac{u_h^{n+1} - u_h^n}{\Delta t},
\end{align}
\end{subequations}
subject to the initial data $(u_h^0, U^0)$ given by \eqref{u0init} and \eqref{CHinit+}.

Note that (\ref{BDF1-IEQ-FEM3}d) can be reformulated as 
\bq\label{CHPwU3}
U^{n+1} = U_h^n+\frac{1}{2}H(u_h^n) (u_h^{n+1} - u_h^n).
\eq
Plugging \eqref{CHPwU3} into (\ref{BDF1-IEQ-FEM3}b) gives the linear system 
\begin{subequations}\label{CHBDF1Lin3}
\begin{align}
&\left(  \frac{u_h^{n+1} - u_h^n}{\Delta t}, \phi \right)=  -A(M(u_h^n); w_h^{n+1}, \phi), \\
&(w_h^{n+1},\psi) = A(\epsilon^2; u_h^{n+1}, \psi)+\left( H(u_h^n)U^{n}_h, \psi \right) + \frac{1}{2}((H(u_h^n))^2 u_h^{n+1}, \psi) - \frac{1}{2}((H(u_h^n))^2 u_h^{n}, \psi).
\end{align}
\end{subequations}
Now, the BDF1-IEQ-FEM3 scheme (\ref{BDF1-IEQ-FEM3}) can be  equivalently expressed as the linear system \eqref{CHBDF1Lin3} along with  (\ref{BDF1-IEQ-FEM3}cd). Solving  scheme \eqref{BDF1-IEQ-FEM3} involves first solving $(u^{n+1}_h, w_h^{n+1})$ from \eqref{CHBDF1Lin3} and then solving $U^{n+1}$ or $U_h^{n+1}$ from (\ref{BDF1-IEQ-FEM3}cd). 
In addition, the following result holds.
\begin{thm}
The BDF1-IEQ-FEM3 scheme (\ref{BDF1-IEQ-FEM3}) admits a unique solution $(u_h^{n+1},w_h^{n+1})\in V_h$ and $U^{n+1} \in C^0(\Omega)$ for any $\Delta t>0$. The solution $u_h^n$ also satisfies total mass conservation \eqref{BDF1MassCons} 
for any $n>0$. Additionally, it satisfies the energy dissipation law:
\bq\label{BDF1-IEQ-FEM3-energy}
\ba
 E(u_h^{n+1}, U^{n+1}_h) \leq &  E(u_h^{n+1}, U^{n+1})  =  E(u_h^{n}, U_h^{n}) - \|U^{n+1} - U_h^n\|^2 \\
& -\Delta t A(M(u_h^n); w_h^{n+1},  w_h^{n+1}) - \frac{1}{2} A(\epsilon^2; u_h^{n+1}-u_h^n, u_h^{n+1} - u_h^n),
\ea
\eq
independent of the time step $\Delta t$.
\end{thm}
In equation \eqref{BDF1-IEQ-FEM3-energy}, the inequality arises due to the property that the $L^2$ projection is a contraction  mapping. The subsequent steps of the proof closely resemble those in  \Cref{CHBDF1thm}.

\begin{rem}
Unlike the BDF1-IEQ-FEM1 scheme \eqref{BDF1-IEQ-FEM1}, scheme \eqref{BDF1-IEQ-FEM3} also avoids the need to solve a large coupled system or computing inverse matrices. Furthermore, in contrast to the BDF1-IEQ-FEM2 scheme \eqref{BDF1-IEQ-FEM2},  scheme \eqref{BDF1-IEQ-FEM3} not only achieves computational efficiency but also ensures a numerical energy dissipation law, which can be directly computed.  
\end{rem}

\section{IEQ-FEMs for the AC equation}\label{secAC}

In this section, we extend the IEQ-FEMs to the Allen-Cahn (AC) equation \eqref{AC}, while continuing to use the finite element space $V_h$ defined in \eqref{eqn.fems}.

\subsection{The semi-discrete finite element scheme}
The semi-discrete finite element scheme for the AC equation (\ref{AC}) is to find $u_h \in V_h$ such that 
\bq\label{ACsemifem}
\ba
(u_{ht}, \phi) = -A(\epsilon^2; u_h,\phi)-\left(F'(u_h),\phi \right),\qquad \forall \phi \in V_h,
\ea
\eq
where the bilinear functional
is given by (\ref{bilinearfunc}). The initial data for $u_h$ is taken as $u_h(x,0)=\Pi u_0(x)$. Recall the free energy $E(\cdot)$ defined in \eqref{freeenergy}, we have the following result.
\begin{lem}
The solution of the semi-discrete finite element scheme (\ref{ACsemifem}) satisfies the energy dissipation law
\be
\frac{d}{dt}E(u_h) = -\|u_{ht}\|^2 \leq 0.
\ee
\end{lem}
\begin{proof}
Setting $\phi=-u_{ht}$ in (\ref{ACsemifem}), then the conclusion holds.
\end{proof}

Building upon the IEQ formulation presented in \Cref{sec:IEQ} and the IEQ-FEMs outlined in \Cref{BDF1-IEQ} for the CH equation, depending on how to approximate the intermediate function $U(x,t)$, the semi-discrete scheme \eqref{ACsemifem} can be equivalently expressed in two ways, depending on how the intermediate function $U(x,t)$ is approximated:

(i) Use $U_h\in V_h$ to approximate $U$. The semi-discrete IEQ-FEM scheme is to find $(u_h, U_h) \in V_h \times V_h$ such that 
\begin{subequations}\label{ACSemiFEM}
\begin{align}
(u_{ht}, \phi) = & -  A(\epsilon^2;u_h,\phi)-\left(H(u_h)U_h,\phi \right),\qquad \forall \phi\in V_h,\\
(U_{ht}, \psi) =&  \left(\frac{1}{2}H(u_h) u_{ht}, \psi \right),\qquad \forall \psi\in V_h,
\end{align}
\end{subequations}
subject to the initial data
\be\label{ACinit}
u_h(x, 0)=\Pi u_0(x), \quad U_h(x, 0)=\Pi\sqrt{F(u_0(x))+B}.
\ee
 (ii) Preserve $U$ in continuous function space $C^0(\Omega)$. The IEQ-FEM scheme is to find $u_h \in V_h$ and $U(x, t)\in C^0(\Omega)$ such that 
\begin{subequations}\label{ACSemiFEM+}
\begin{align}
(u_{ht}, \phi) = & -  A(\epsilon^2;u_h,\phi)-\left(H(u_h)U,\phi \right),\qquad \forall \phi\in V_h,\\
U_{t} =&  \frac{1}{2}H(u_h) u_{ht},
\end{align}
\end{subequations}
subject to the initial data
\be\label{ACinit+}
u_h(x, 0)=\Pi u_0(x), \quad U(x, 0)=\sqrt{F(u_0(x))+B}.
\ee
By taking $\phi=u_{ht}$ and $\psi=2U_h$ in (\ref{ACSemiFEM}) and summing (\ref{ACSemiFEM}a) and (\ref{ACSemiFEM}b), it follows that the semi-discrete IEQ-FEM (\ref{ACSemiFEM}) satisfies the energy dissipation law,
\bqs
\frac{d}{dt}E(u_h, U_h) = -\|u_{ht}\|^2 \leq 0.
\eqs
Similarly, by setting $\phi=u_{ht}$ in (\ref{ACSemiFEM+}a), (\ref{ACSemiFEM+}b),  and taking the $L^2$ inner product with $2U$, and by summing (\ref{ACSemiFEM+}a) and (\ref{ACSemiFEM+}b), we can establish that the semi-discrete IEQ-FEM (\ref{ACSemiFEM+}) adheres to the energy dissipation law:
\bqs
\frac{d}{dt}E(u_h, U) = -\|u_{ht}\|^2 \leq 0.
\eqs
Here, it's important to note that $U(x,t) \not \in V_h$ in scheme (\ref{ACSemiFEM+}) unless $H(u_{h})$ is a constant.

\subsection{BDF1-IEQ-FEMs}

The selection of the space for  approximating the intermediate function $U(x, t)$ in the semi-discrete IEQ-FEM schemes \eqref{ACSemiFEM} and \eqref{ACSemiFEM+} gives rise to different fully discrete IEQ-FEM schemes. For $n \geq 0$, let $u_h^n \in V_h$ approximate $u(x, t^n)$, where $t^n = n\Delta t$ and $\Delta t > 0$ is the time step. We denote the approximation of $U(x, t^n)$ by $U_h^n \in V_h$ or $U^n \in C^0(\Omega)$. The choice of this approximation space significantly influences the resulting numerical methods and their inherent properties.

\subsubsection{Method 1 (BDF1-IEQ-FEM1)}

Given $(u_h^n, U_h^n) \in V_h  \times V_h$, the first order BDF1-IEQ-FEM1 scheme based on the semi-discrete FEM scheme (\ref{ACSemiFEM}) is to find $(u^{n+1}_h, U_h^{n+1}) \in V_h  \times V_h $ such that for $\phi, \psi \in V_h$,
\begin{subequations}\label{BDF1-FEM1}
\begin{align}
\left( \frac{u_h^{n+1} - u_h^n}{\Delta t}, \phi \right)= & - A(\epsilon^2; u_h^{n+1}, \phi)-\left( H(u_h^n)U_h^{n+1}, \phi \right),\\
\left(\frac{U_h^{n+1} - U_h^n}{\Delta t},\psi \right) = & \left(\frac{1}{2}H(u_h^n) \frac{u_h^{n+1} - u_h^n}{\Delta t}, \psi\right).
\end{align}
\end{subequations}
Here, $u^{n+1}_h, U_h^{n+1}$ constitute a coupled system, and the initial data is given by (\ref{ACinit}). Similar to the CH equation, the linear system for the scheme \eqref{BDF1-FEM1} takes the form
\begin{equation}\label{AC-BDF1-IEQ-FEM1-mat}
\begin{aligned}
\left(
  \begin{array}{ccc}
    \frac{G}{\Delta t}+\epsilon^2D &  G_1 \\
    -\frac{G_1}{2\Delta t} & \frac{G}{\Delta t} \\
  \end{array}
\right)\left(
  \begin{array}{ccc}
    \vec u^{n+1} \\
    \vec U^{n+1} \\
  \end{array}
\right)=\left(\begin{array}{ccc}
    \vec g_1\\
    \vec g_3
  \end{array}
\right),
\end{aligned}
\end{equation}
with the notations $\vec{u}^{n}$, $\vec{U}^{n}$, the matrices $D$, $G$, $G_1$ and the vectors $\vec g_1$, $\vec g_3$ defined as those in \eqref{BDF1-IEQ-FEM1-mat}.
One approach to solving 
the linear system for the scheme \eqref{AC-BDF1-IEQ-FEM1-mat} is direct,  and the other approach is to reduce the matrix equation \eqref{AC-BDF1-IEQ-FEM1-mat} to
\begin{equation}\label{AC-BDF1-IEQ-FEM1-mat1}
\begin{aligned}
\left(
    \frac{G}{\Delta t}+\epsilon^2D-\frac{1}{2}G_1G^{-1}G_1 
\right)
    \cdot\vec u^{n+1} 
=
    \vec g_1+ G_1 G^{-1}\Delta t \vec g_3 .
\end{aligned}
\end{equation}

\begin{lem}
The first order fully discrete BDF1-IEQ-FEM1 scheme (\ref{BDF1-FEM1}) admits a unique solution $(u_h^{n+1}, U_h^{n+1}) \in V_h\times V_h$ satisfying the energy dissipation law,
\bq\label{ACCoStab}
E(u_h^{n+1}, U_h^{n+1}) = E(u_h^{n}, U_h^{n}) - \frac{1}{\Delta t}\|u_h^{n+1} - u_h^n\|^2 - \frac{1}{2} A(\epsilon^2; u_h^{n+1}-u_h^n, u_h^{n+1} - u_h^n) - \|U_h^{n+1} - U_h^n\|^2.
\eq
\end{lem}

In the scheme \eqref{BDF1-FEM1}, a similar challenge arises, reminiscent of the BDF1-IEQ-FEM scheme for the CH equation. Specifically, the solution of the BDF1-IEQ-FEM scheme \eqref{BDF1-FEM1} can be obtained either by solving the augmented  linear system with two unknowns or by utilizing a simplified linear system that involves computing the inverse of the mass matrix, $G^{-1}$.

\subsubsection{Method 2 (BDF1-IEQ-FEM2)}
Note that the intermediate function $U(x,t)$ in the semi-discrete IEQ-FEM scheme (\ref{ACSemiFEM+}) is defined pointwise. For a given  $u^{n}_h \in V_h$  and $U^n(x)\in C^0(\Omega)$,  the first order fully discrete BDF1-IEQ-FEM2 scheme seeks to find $u^{n+1}_h \in V_h$ and $U^{n+1} \in C^0(\Omega)$ such that,  for $  \phi \in V_h $,
\begin{subequations}\label{BDF1-FEM2}
\begin{align}
\left( \frac{u_h^{n+1} - u_h^n}{\Delta t}, \phi \right)= & - A(\epsilon^2; u_h^{n+1}, \phi)-\left( H(u_h^n)U^{n+1}, \phi \right),\\
\frac{U^{n+1} - U^n}{\Delta t} = & \frac{1}{2}H(u_h^n) \frac{u_h^{n+1} - u_h^n}{\Delta t}.
\end{align}
\end{subequations}
The initial data is given by (\ref{ACinit+}).

By (\ref{BDF1-FEM2}b), we can derive 
\bq\label{ACPwU}
U^{n+1} = U^n+\frac{1}{2}H(u_h^n) (u_h^{n+1} - u_h^n).
\eq
Substituting $U^{n+1}$ from (\ref{ACPwU}) into (\ref{BDF1-FEM2}a) gives the following equation in terms of $u_h^{n+1}$,
\bq\label{ACPWuh}
\left( \frac{u_h^{n+1} - u_h^n}{\Delta t}, \phi \right)=  - A(\epsilon^2; u_h^{n+1}, \phi)-\left( H(u_h^n)U^{n}, \phi \right) - \frac{1}{2}(H(u_h^n)^2 u_n^{n+1}, \phi) + \frac{1}{2}(H(u_h^n)^2 u_n^{n}, \phi). 
\eq
In this context, the system (\ref{BDF1-FEM2}) can be considered equivalent to the system formed by (\ref{ACPwU}) and (\ref{ACPWuh}), where the unknowns in equation (\ref{ACPWuh}) involve only $u_h^{n+1}$. Consequently, one can initially obtain $u_h^{n+1}$ by solving (\ref{ACPWuh}). Subsequently, the evaluation of $U^{n+1}$ from (\ref{ACPwU}) becomes straightforward. Compared to the IEQ-FEM scheme (\ref{BDF1-FEM1}) based on the semi-discrete FEM scheme (\ref{ACSemiFEM}), the BDF1-IEQ-FEM2 scheme (\ref{BDF1-FEM2}) eliminates the need to solve a coupled linear system with double unknowns.

\begin{thm}\label{ACPWLem}
Given $u_h^n \in V_h$ and $U^{n}\in C^0(\Omega)$, the first order fully discrete IEQ-FEM scheme (\ref{BDF1-FEM2}) admits a unique solution $u_h^{n+1}\in V_h$ and $U^{n+1}\in C^0(\Omega)$ satisfying the energy dissipation law,
\bq\label{ACPWStab}
E(u_h^{n+1}, U^{n+1}) = E(u_h^{n}, U^{n}) - \frac{1}{\Delta t}\|u_h^{n+1} - u_h^n\|^2 - \frac{1}{2} A(\epsilon^2; u_h^{n+1}-u_h^n, u_h^{n+1} - u_h^n) - \|U^{n+1} - U^n\|^2,
\eq
for any $\Delta t>0$.
\end{thm}
\begin{proof}
We first prove the energy stability. By setting $\phi=u_h^{n+1} - u_h^n$ in (\ref{BDF1-FEM2}a), (\ref{BDF1-FEM2}b),  taking $L^2$ inner product with $2\Delta t U^{n+1}$, and taking summation, we have
\bqs
\ba
-\frac{1}{\Delta t}\|u_h^{n+1} - u_h^n\|^2 = & A(\epsilon^2; u_h^{n+1}, u_h^{n+1} - u_h^n) + 2(U^{n+1} - U^n, U^{n+1}),\\
= & \frac{1}{2}\left(A(\epsilon^2; u_h^{n+1}, u_h^{n+1})-A(\epsilon^2; u_h^{n}, u_h^{n}) + A(\epsilon^2; u_h^{n+1}-u_h^n, u_h^{n+1} - u_h^n)\right)\\
& + \left( \|U^{n+1}\|^2 - \|U^{n}\|^2 + \|U^{n+1} - U^n\|^2 \right), 
\ea
\eqs
which gives the energy dissipation law (\ref{ACPWStab}).

It is important to note that the scheme (\ref{BDF1-FEM2}) shares the same solution as the system formed by (\ref{ACPwU}) and (\ref{ACPWuh}). As (\ref{ACPWuh}) constitutes a linear system in a finite-dimensional space, the existence of the solution $u_h^{n+1}$ is equivalent to its existence. This, coupled with (\ref{ACPwU}), further implies the existence and uniqueness of the solution for the scheme (\ref{BDF1-FEM2}).
Assuming that (\ref{ACPWuh}) has two solutions, we denote their differences as $\tilde u_h^{n+1}$. This leads to the expression:
\bq\label{ACPointwise1}
\ba
\left( \frac{\tilde u_h^{n+1}}{\Delta t}, \phi \right)=  - A(\epsilon^2;\tilde u_h^{n+1}, \phi)-\frac{1}{2}(H(u_h^n)^2 \tilde u_n^{n+1}, \phi).
\ea
\eq
By taking $\phi=\tilde u_h^{n+1}$ in (\ref{ACPointwise1}), it follows
$$
\frac{1}{\Delta t}\|\tilde u_h^{n+1}\|^2+\epsilon^2|\tilde u_h^{n+1}|^2_{H^1(\Omega)} + \frac{1}{2}\|H(u_h^n)\tilde u_h^{n+1}\|^2=0,
$$
which gives $\tilde u_n^{n+1} = 0$. Thus, it follows the existence and uniqueness of the scheme (\ref{BDF1-FEM2}).
\end{proof}

\begin{rem}
The scheme \eqref{BDF1-FEM2} is computationally cheap, but the computed energy $E(u_h^n, U^n)$ can only be approximated by $E(u_h^n, \Pi U^n)$, which may not exhibit an energy dissipation, as demonstrated by the counter-example given in \Cref{cexmAC2-n} for the AC equation. Several numerical strategies, outlined in \Cref{Enerdis}, can be employed to effectively manage the energy dissipation of $E(u_h^n, \Pi U^n)$. 
\end{rem}
Next, we introduce a methods that is independent of its numerical setting, while  ensuring energy dissipation.

\subsubsection{Method 3 (BDF1-IEQ-FEM3)}

Following the IEQ-DG schemes introduced in \cite{LY19, LY21, LY22, Y19} for the Swift-Hohenberg equation and the CH equation, we introduce the first order fully discrete BDF1-IEQ-FEM3 scheme, which aims to find  $u^{n+1}_h \in V_h $ and $U^{n+1} \in C^0(\Omega)$ such that for $\phi \in V_h $,
\begin{subequations}\label{BDF1-FEM3}
\begin{align}
\left(  \frac{u_h^{n+1} - u_h^n}{\Delta t}, \phi \right)= & -\epsilon^2 A( u_h^{n+1}, \psi)-\left( H(u_h^n)U^{n+1}, \psi \right),\\
U^n_h= & \Pi U^n,\\
\frac{U^{n+1} - U_h^n}{\Delta t} = & \frac{1}{2}H(u_h^n) \frac{u_h^{n+1} - u_h^n}{\Delta t}.
\end{align}
\end{subequations}
The initial data is given by (\ref{ACinit+}). 

By (\ref{BDF1-FEM3}c), we have
\bq\label{ACPJU}
U^{n+1} = U_h^n+\frac{1}{2}H(u_h^n) (u_h^{n+1} - u_h^n).
\eq
Substituting $U^{n+1}$ from (\ref{ACPJU}) into (\ref{BDF1-FEM3}a) gives
\bq\label{ACPJuh}
\left( \frac{u_h^{n+1} - u_h^n}{\Delta t}, \phi \right)=  - A(\epsilon^2; u_h^{n+1}, \phi)-\left( H(u_h^n)U_h^{n}, \phi \right) - \frac{1}{2}(H(u_h^n)^2 u_n^{n+1}, \phi) + \frac{1}{2}(H(u_h^n)^2 u_n^{n}, \phi). 
\eq
Note that the scheme (\ref{BDF1-FEM3}) is equivalent to the system formed by (\ref{ACPJuh}) and (\ref{ACPJU}). Here, (\ref{ACPJuh}) is an equation concerning the solution $u_h^{n+1}$ only.
Similar to the BDF1-IEQ-FEM2 scheme (\ref{BDF1-FEM2}), the scheme (\ref{BDF1-FEM3}) can also avoid the need to solve a coupled linear system with double unknowns. 
Moreover, the solution of \eqref{BDF1-FEM3} satisfies the following result.

\begin{thm}\label{ACPJLem}
Given $u_h^n \in V_h$ and $U^{n} \in C^0(\Omega)$, the first order fully discrete BDF1-IEQ-FEM3 scheme (\ref{BDF1-FEM3}) admits a unique solution $u_h^{n+1}\in V_h$ and $U^{n+1}\in C^0(\Omega)$ satisfying the energy dissipation law,
\bq\label{ACPJStab}
\ba
E(u_h^{n+1}, U_h^{n+1}) \leq & E(u_h^{n+1}, U^{n+1}) = E(u_h^{n}, U_h^{n}) - \frac{1}{\Delta t}\|u_h^{n+1} - u_h^n\|^2 \\
&- \frac{1}{2} A(\epsilon^2; u_h^{n+1}-u_h^n, u_h^{n+1} - u_h^n) - \|U^{n+1} - U_h^n\|^2.
\ea
\eq
\end{thm}
The proof of \Cref{ACPJLem} is similar to that of \Cref{ACPWLem}. The first inequality in (\ref{ACPJStab}) is derived from the fact that $L^2$ projection is a contraction operation. 

We present the second order fully discrete IEQ-FEMs for the AC equation \eqref{AC} in \Cref{AC2nd}.




\section{Numerical examples}
\label{secNumer}
Upon rescaling with $t' = \epsilon t$ and $t' = \epsilon^2 t$ for the CH equation and AC equation, respectively, the  CH problem (\ref{CH}) can be equivalently expressed as the following rescaled CH problem \cite{FW2008}:
\bq\label{CH1}
\begin{aligned}
u_{t}-\nabla \cdot (M(u) \nabla w)&=0,\qquad  (x,t)\in \Omega\times(0,T],\\
w+\epsilon\Delta u-\frac{1}{\epsilon}f(u)&=0,\qquad  (x,t)\in  \Omega\times(0,T],\\
u &=u_{0},\quad\,\ (x,t)\in \Omega\times\{t=0\},
\end{aligned}
\eq
and AC problem (\ref{AC}) is equivalent to the rescaled AC problem \cite{FP2003, ZD2009}:
\begin{equation}\label{AC1}
\begin{aligned}
 u_t = \Delta u - \frac{1}{\epsilon^2}F'(u) \qquad &{\it in}\ \Omega\ \times (0, T],\\ 
\nabla u\cdot \pmb{n}=0 \qquad &{\it on}\ \partial\Omega\times [0, T],\\
u=u_{0}\qquad &{\it on}\  \Omega\times \{t=0\}.
\end{aligned}
\end{equation}

In this section, we employ several numerical examples to showcase the effectiveness of the proposed first order and second order (in Appendix) fully discrete IEQ-FEM schemes in solving both the CH problem (\ref{CH}) or (\ref{CH1}) and the AC problem (\ref{AC}) or (\ref{AC1}). For both temporal and spatial accuracy tests, we introduce a source term into the AC equation in (\ref{AC1}) as follows:
\begin{equation}\label{2e1}
u_{t}-\Delta u+\frac{1}{\epsilon^{2}}F'(u)=s(x,t),
\end{equation}
and into the CH equation in (\ref{CH}) as:
\begin{equation}\label{2e1ac}
u_{t}-\nabla\cdot(M(u)\nabla w)=s(x,t).
\end{equation}
Additionally, unless otherwise specified, we will use the default value $M(u)=1$ for mobility and the double well potential \eqref{dwp}.

\subsection{The CH equation}\label{secNum}
In this section, we initially present an example to validate the temporal and spatial convergence rates of the proposed methods for the CH equation (\ref{2e1ac}). Subsequently, we compare the CPU time of first-order schemes for the CH equation (\ref{CH}). Furthermore, we provide more examples to evaluate the performance of proposed methods for solving the CH equations, as illustrated in $\mathbf{Example\ \ref{cexm2-n}}$-$\mathbf{Example\ \ref{cexm5}}$.
\subsubsection{Convergent rates}
\begin{example}\label{2e1a} \cite{LY21} (Temporal and spatial accuracy test with constant mobility and double-well potential)
Consider the CH equation \eqref{2e1ac} in  $\Omega=[-\pi,3\pi]^{2}$ with the parameters $\epsilon=1$, $B=1$. 
The exact solution $u$ satisfies
\bq
u(x,t)=\mu(x,t),
\eq
and the corresponding right term $s(x,t)$ in \eqref{2e1ac} is taken as:
\bq
s(x,t)=-\frac{\mu(x,t)}{4}+\frac{\epsilon^{2}\mu(x,t)}{4}-\frac{3\mu(x,t)\nu(x,t)}{2}+\frac{3(\mu(x,t))^{3}}{2}-\frac{\mu(x,t)}{2},
\eq
where
\bq
\begin{aligned}
\mu(x,t)&=0.1e^{-t/4}\sin(x/2)\sin(y/2),\\
\nu(x,t)&=(0.1e^{-t/4}\cos(x/2)\sin(y/2))^{2}+(0.1e^{-t/4}\sin(x/2)\cos(y/2))^{2}.
\end{aligned}
\eq

 \begin{figure}[!htbp]
 $\begin{array}{cc}
 \includegraphics[width=5.5cm,height=4.5cm]{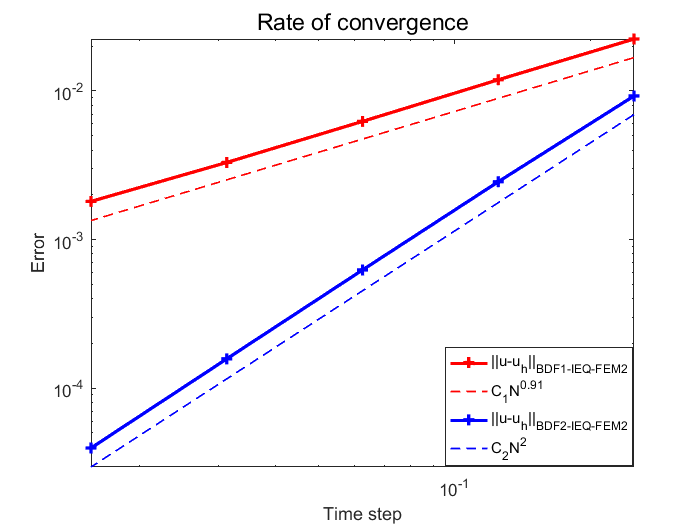}
 \includegraphics[width=5.5cm,height=4.5cm]{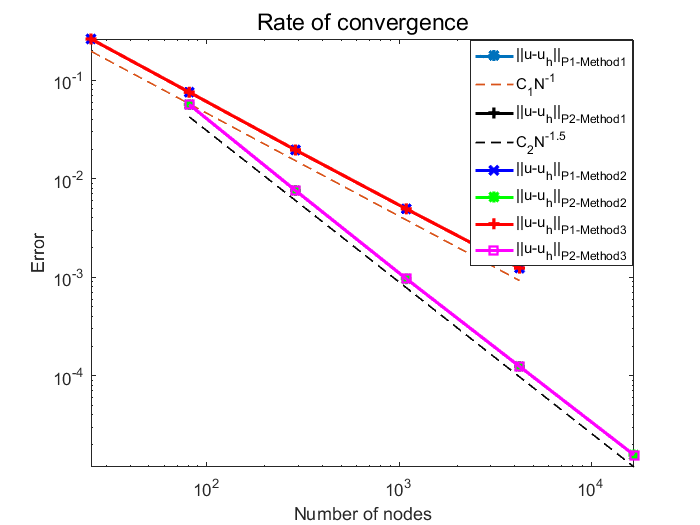}
 \end{array}$
\caption{$\mathbf{Example\ \ref{2e1a}}$, $L_2$ errors and convergent rates for the CH equation, Left: Time accuracy test ($P_{2}$ element, $T=1$, uniform mesh with $(h_x, h_y)=(\frac{4 \pi}{200},\frac{4 \pi}{200})$); Right: Spatial accuracy test ($\Delta t=1e^{-6},\, T=1e^{-3}$).}\label{convergentrateCH}
 \end{figure}
In Figure \ref{convergentrateCH}, the $L_2$ errors and convergent rates for the CH equation are shown.  
From the left figure, we find that the time accuracy is first order for the BDF1-IEQ-FEM2 scheme (\ref{BDF1-IEQ-FEM2}), and second order for the BDF2-IEQ-FEM2 scheme (\ref{BDF2-IEQ-FEM2}). 
For spatial accuracy, it is second order and third order for all proposed schemes based on $P_1$ elements and $P_2$ elements, respectively.
\end{example}

\subsubsection{Comparison between different schemes}
\begin{example}\label{exmS3} \cite{LY21}  Consider the CH equation $(\ref{CH})$ in $\Omega=[-0.5,0.5]^{2}$ with constant mobility $M(u)=1$, the logarithmic Flory-Huggins potential
$$F(u)=600(u\ln u+(1-u)\ln(1-u))+1800u(1-u),$$
and the parameters $\epsilon=1$. The equation is subject to the following initial data
\bq\label{exmS3e5a}
u_{0}=\left\{\begin{aligned}
&0.71, \qquad |x|\leq 0.2, \ |y|\leq 0.2,\\
&0.69, \qquad \text{otherwise}.
\end{aligned}\right.
\eq
The time step is chosen as $\Delta t=10^{-7}$, and the constant $B=100$ for all the involved IEQ-FEMs.

\begin{figure}[!htbp]
$\begin{array}{cccc}
 \includegraphics[width=5.5cm,height=4.5cm]{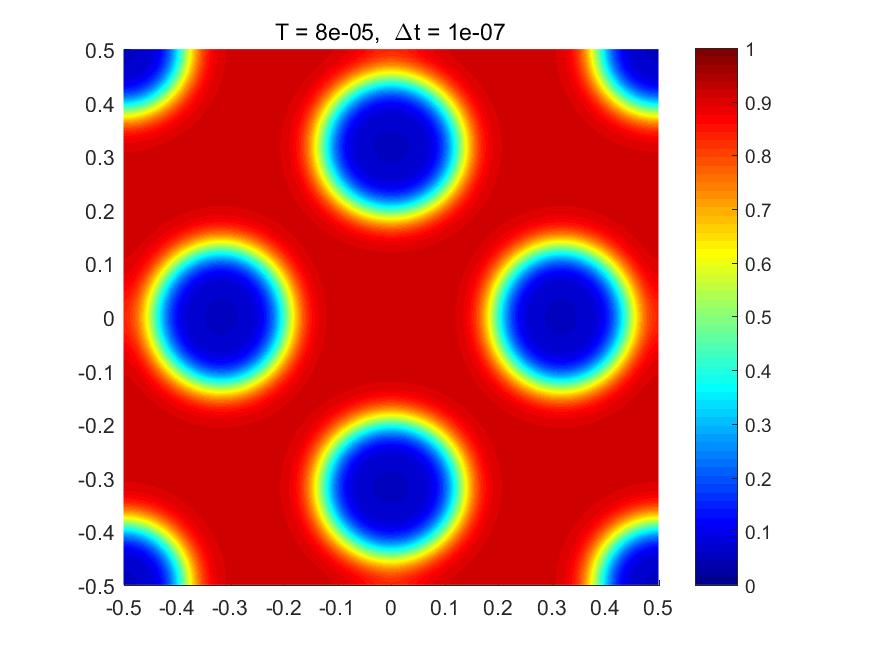}&\vspace{-0.1cm}
 \includegraphics[width=5.5cm,height=4.5cm]{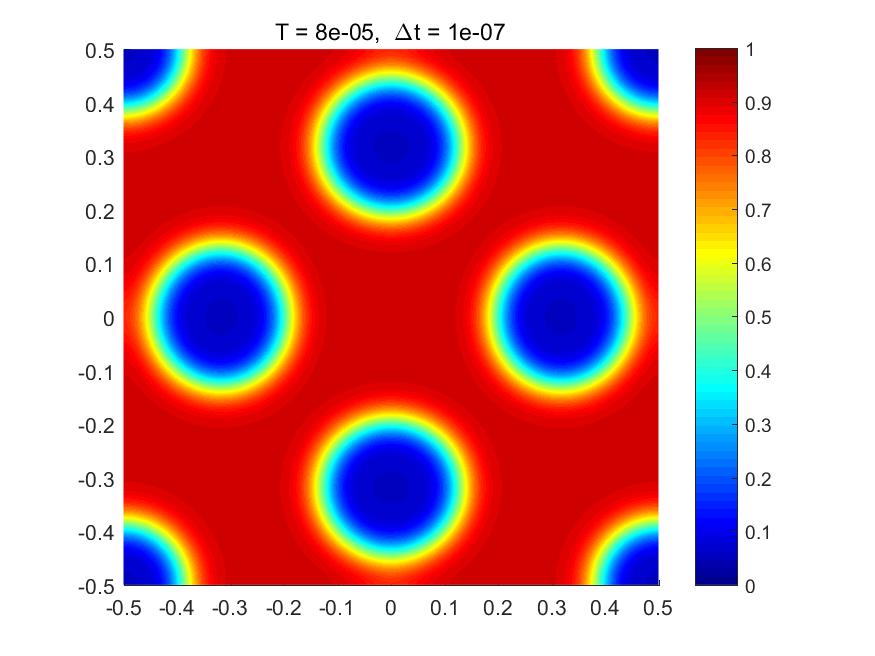}&\vspace{-0.1cm}
 \includegraphics[width=5.5cm,height=4.5cm]{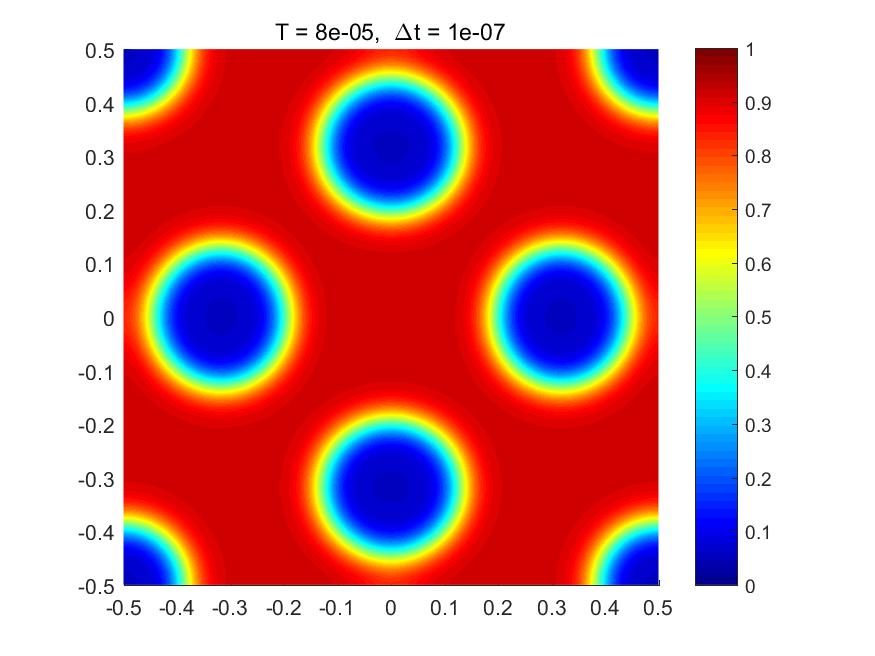}\\
\includegraphics[width=5.5cm,height=4.5cm]{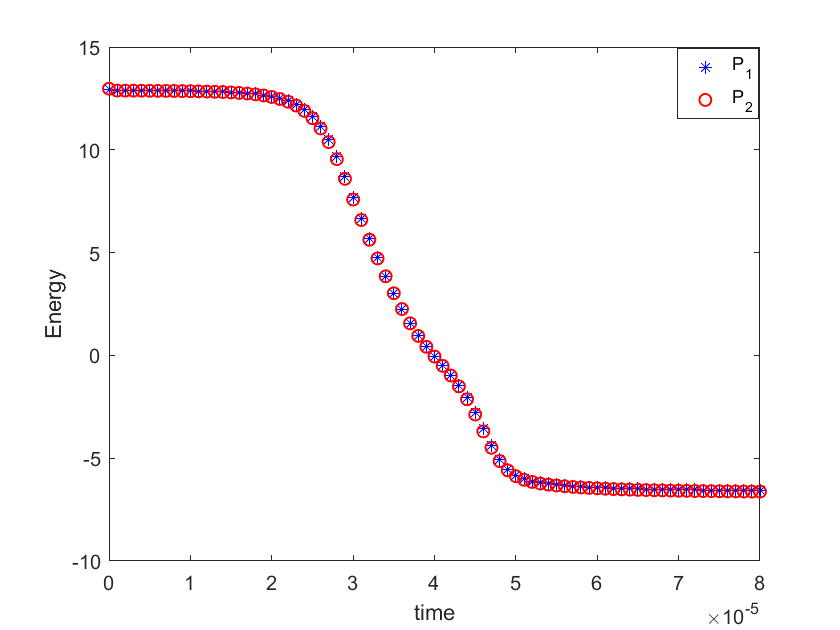}&\vspace{-0.1cm}
\includegraphics[width=5.5cm,height=4.5cm]{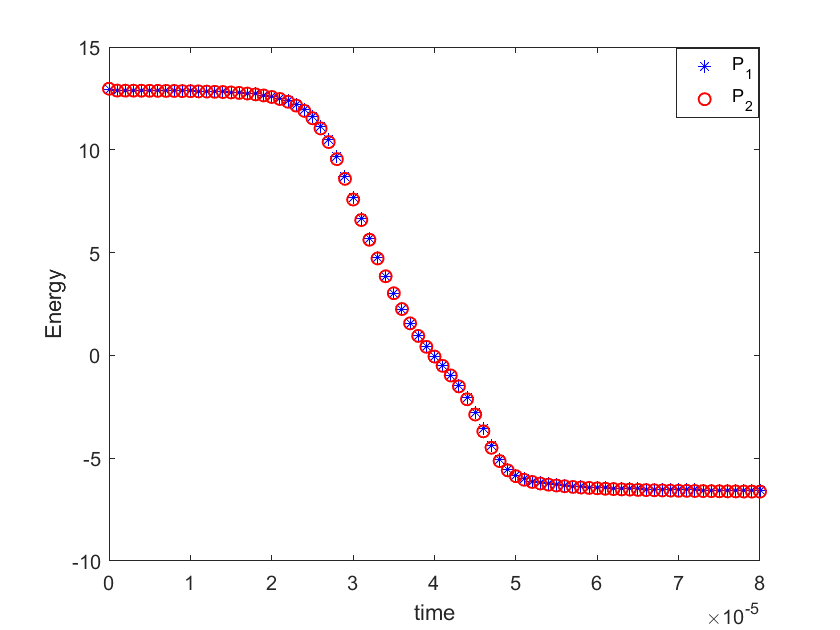}&\vspace{-0.1cm}
\includegraphics[width=5.5cm,height=4.5cm]{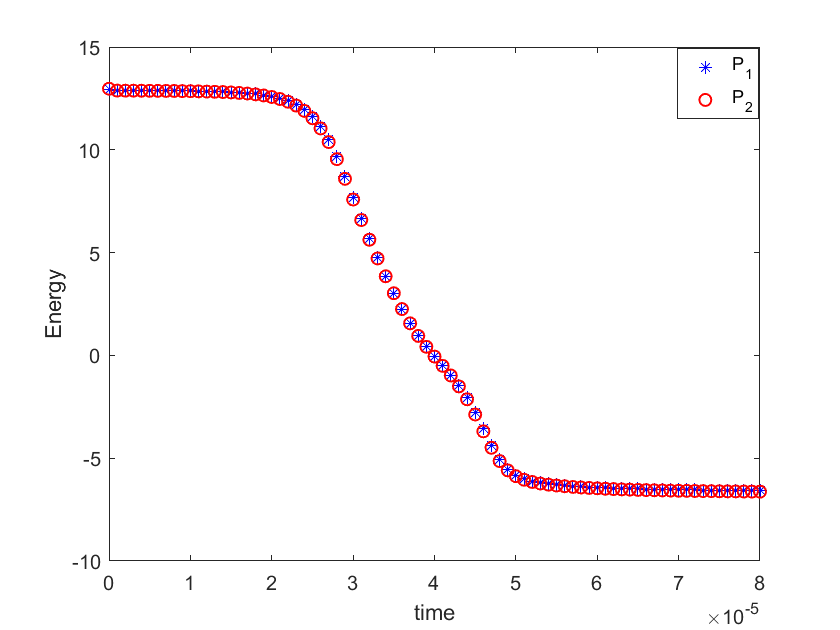}&\vspace{-0.1cm}
\end{array}$\vspace{-0.3cm}
\caption{$\mathbf{Example\ \ref{exmS3}\ (CH, T=8e-5)}$, First order scheme; Left: BDF1-IEQ-FEM1 scheme (\ref{BDF1-IEQ-FEM1}), Middle: BDF1-IEQ-FEM2 scheme (\ref{BDF1-IEQ-FEM2}), Right: BDF1-IEQ-FEM3 scheme (\ref{BDF1-IEQ-FEM3}); First row: Numerical solution; Second row: Time evolution of the discrete energy.}\label{expsu111}
\end{figure}


\begin{table}[!htbp]
\centering
\caption{$\mathbf{Example\ \ref{exmS3}\ (CH, T=8e-5)}$, CPU time of the BDF1-IEQ-FEMs (11th Gen Intel(R) Core(TM) i5-1135G7 @ 2.40GHz 2.42GHz). (`$--$' represents running out of memory)}\label{Table1:S3}
\begin{tabular}{| c | c | c | c |}
\hline
Methods    & scheme (\ref{BDF1-IEQ-FEM1})  &  scheme (\ref{BDF1-IEQ-FEM2}) & scheme (\ref{BDF1-IEQ-FEM3})                      \\  \hline
CPU time ($P_{1}$)    & 822s   &  195s     & 217s                             \\  \hline
CPU time ($P_{2}$)   & 6815s   &     1078s   & 1209s                             \\  \hline
\end{tabular}
\end{table}

The snapshots depicting numerical solutions and discrete energy generated by the three BDF1-IEQ-FEMs (\ref{BDF1-IEQ-FEM1}), (\ref{BDF1-IEQ-FEM2}), and (\ref{BDF1-IEQ-FEM3}) are illustrated in Figure \ref{expsu111}. It is evident that the discrete energy decreases over time, and there are no discernible differences in the solution skeletons.
The corresponding CPU time is presented in Table \ref{Table1:S3}, indicating that the CPU time for the BDF1-IEQ-FEM2 scheme (\ref{BDF1-IEQ-FEM2}) and the BDF1-IEQ-FEM3 scheme (\ref{BDF1-IEQ-FEM3}) is significantly smaller than that of the BDF1-IEQ-FEM1 scheme (\ref{BDF1-IEQ-FEM1}). 
In addition, solving BDF1-IEQ-FEM1 scheme (\ref{BDF1-IEQ-FEM1}) by evaluating the inverse of the mass matrix \eqref{BDF1-IEQ-FEM1-mat1} is more expensive than solving the augmented system \eqref{BDF1-IEQ-FEM1-mat}.
Based on our observation, the focus in the subsequent tests will be primarily on Method 2 and Method 3.
\end{example}

\subsubsection{Numerical solutions for the CH equation}

In the following three examples, we numerically investigate the performance of the proposed IEQ-FEMs for solving the CH equation.

\begin{example}\label{cexm2-n} \cite{FW2008} In this example, we consider the CH equation (\ref{CH1}) in $\Omega=[-1,1]^{2}$
with the parameter $\epsilon=0.01,\,\Delta t=10^{-7}$, and the initial condition
\[\begin{aligned}
u_{0}=&\tanh\Big(\big((x-0.3)^{2}+y^{2}-0.2^{2}\big)/\epsilon\Big)\tanh\Big(\big((x+0.3)^{2}+y^{2}-0.2^{2}\big)/\epsilon\Big)\times\\
&\tanh\Big(\big(x^{2}+(y-0.3)^{2}-0.2^{2}\big)/\epsilon\Big)\tanh\Big(\big(x^{2}+(y+0.3)^{2}-0.2^{2}\big)/\epsilon\Big).
\end{aligned}\]

\noindent\textbf{Test Case 1.} 
We solve this problem by the BDF1-IEQ-FEM2 scheme \eqref{BDF1-IEQ-FEM2} and the BDF1-IEQ-FEM3 scheme \eqref{BDF1-IEQ-FEM3} with mesh size $h=0.01$ but different constants $B$. The evolution of the discrete energy $E(u_h^n, \Pi U^n)$ for BDF1-IEQ-FEM2 and $E(u_h^n, U_h^n)$ for BDF1-IEQ-FEM3 is reported in \Cref{Cexp3u-B}, from which we observe that solutions of the BDF1-IEQ-FEM2 satisfy the energy dissipation law only when $B$ is appropriately large, as explained in \Cref{Enerdis}, while solutions of BDF1-IEQ-FEM3 satisfy the energy dissipation law independent of the choice of $B$.

\begin{figure}[!htbp]
$\begin{array}{cc}
\includegraphics[width=5.5cm,height=4.5cm]{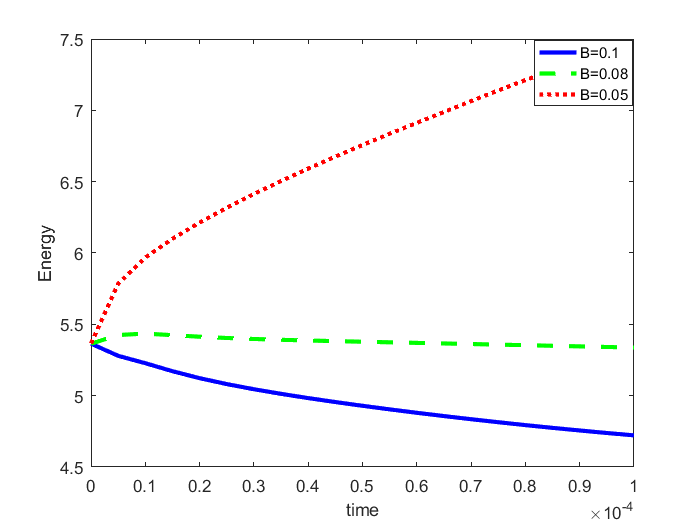}&
\includegraphics[width=5.5cm,height=4.5cm]{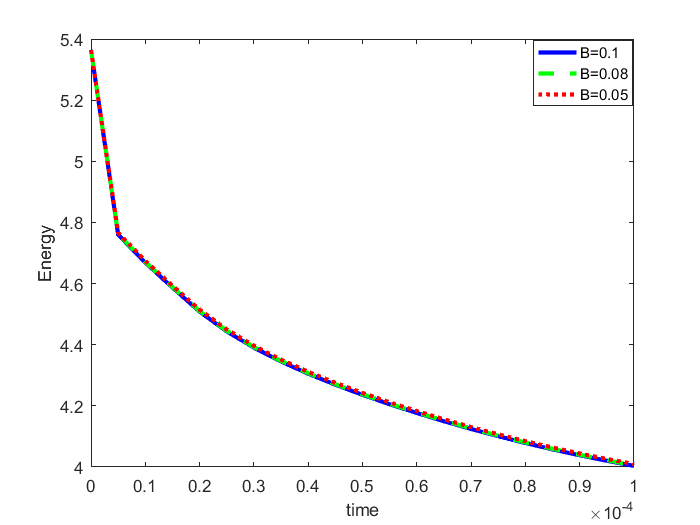}
\end{array}$\vspace{-0.2cm}
\caption{$\mathbf{Example\ \ref{cexm2-n}\ (CH)}$, Effect of constant B on the discrete energy, Left: BDF1-IEQ-FEM2 scheme (\ref{BDF1-IEQ-FEM2}), Right: BDF1-IEQ-FEM3 scheme (\ref{BDF1-IEQ-FEM3}).}\label{Cexp3u-B}
\end{figure}

\noindent\textbf{Test Case 2.} 
We continue with this problem by testing BDF2-IEQ-FEMs with $B=1$ and mesh size $h=0.01$. The solutions produced by the BDF2-IEQ-FEM2 scheme (\ref{BDF2-IEQ-FEM2}) and the BDF2-IEQ-FEM3 scheme (\ref{BDF2-IEQ-FEM3}) are shown in \Cref{Cexp3u}, from which we see that the BDF2-IEQ-FEMs produce high-quality numerical solutions for the CH equation.
The evolution of the discrete energy and  total mass are shown in \Cref{Cexp3u-B-mass}. The results indicate that both methods preserve  the total mass and satisfy the energy dissipation law.
\begin{figure}[!htbp]
$\begin{array}{ccc}
\includegraphics[width=5.5cm,height=4.5cm]{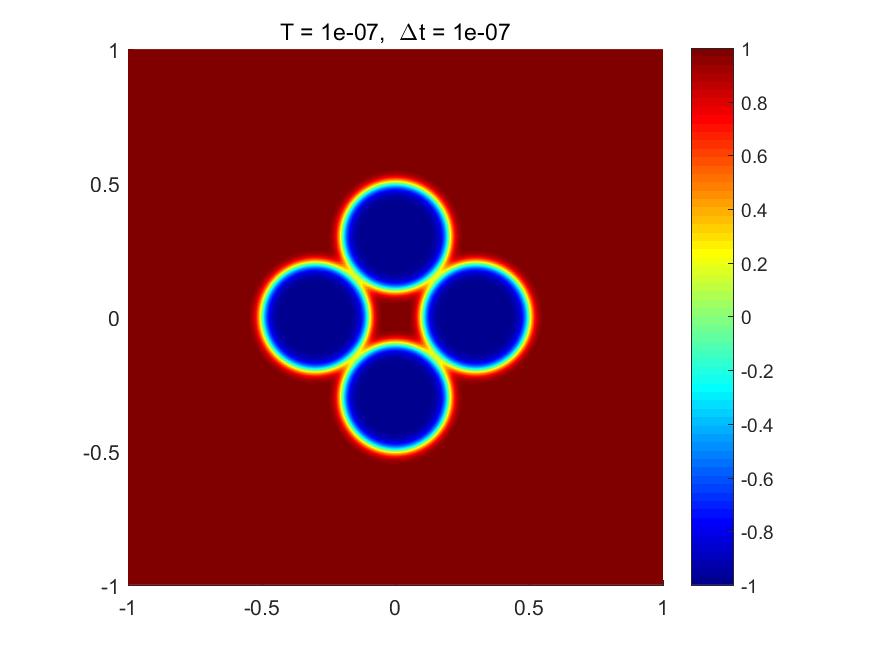}&
\includegraphics[width=5.5cm,height=4.5cm]{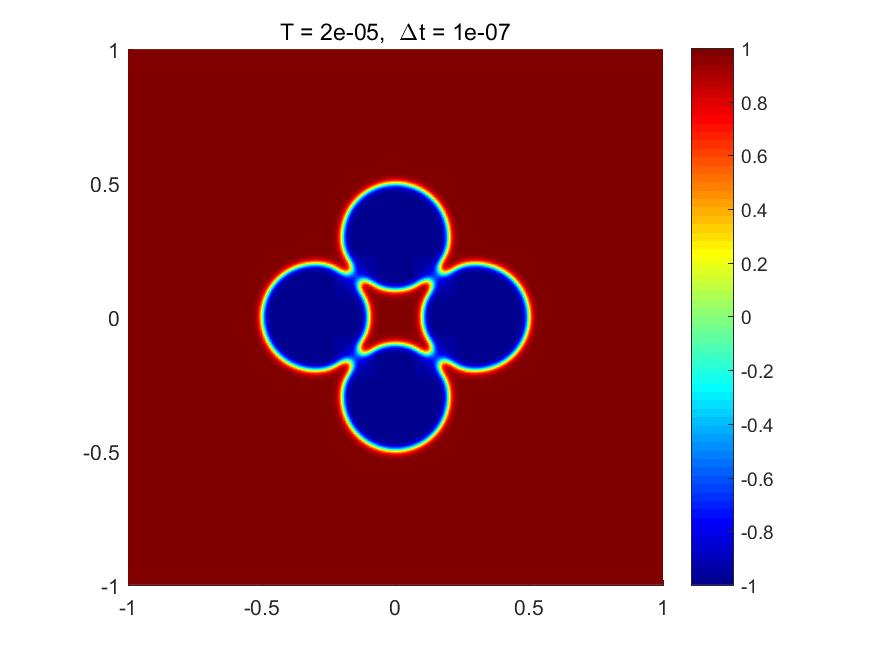}&
\includegraphics[width=5.5cm,height=4.5cm]{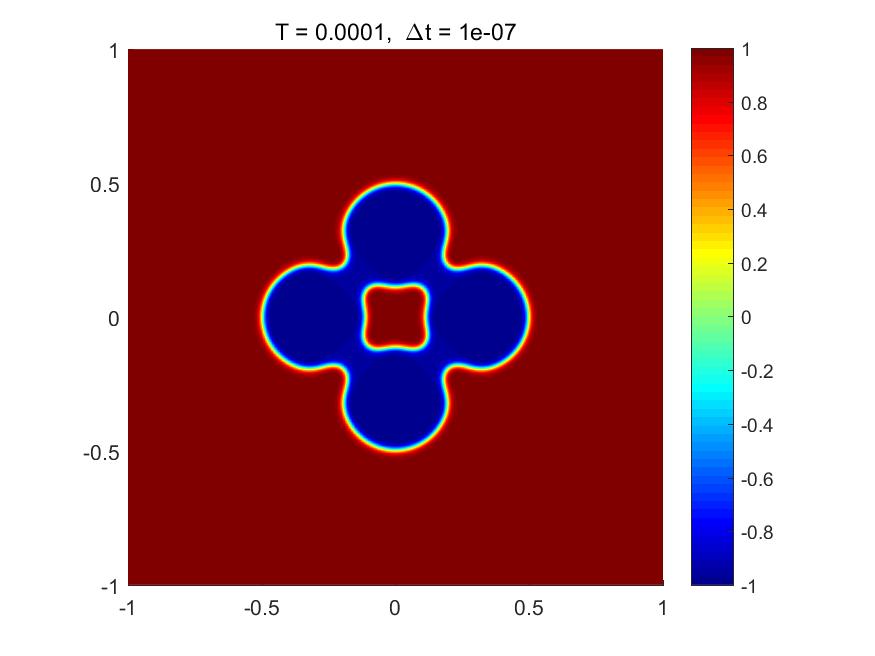}\\
\includegraphics[width=5.5cm,height=4.5cm]{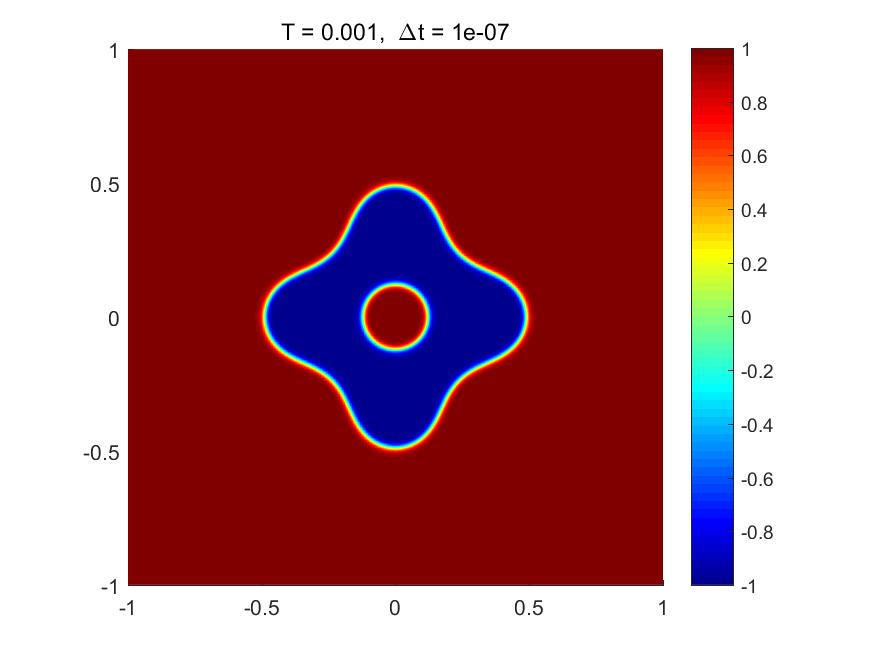}&
\includegraphics[width=5.5cm,height=4.5cm]{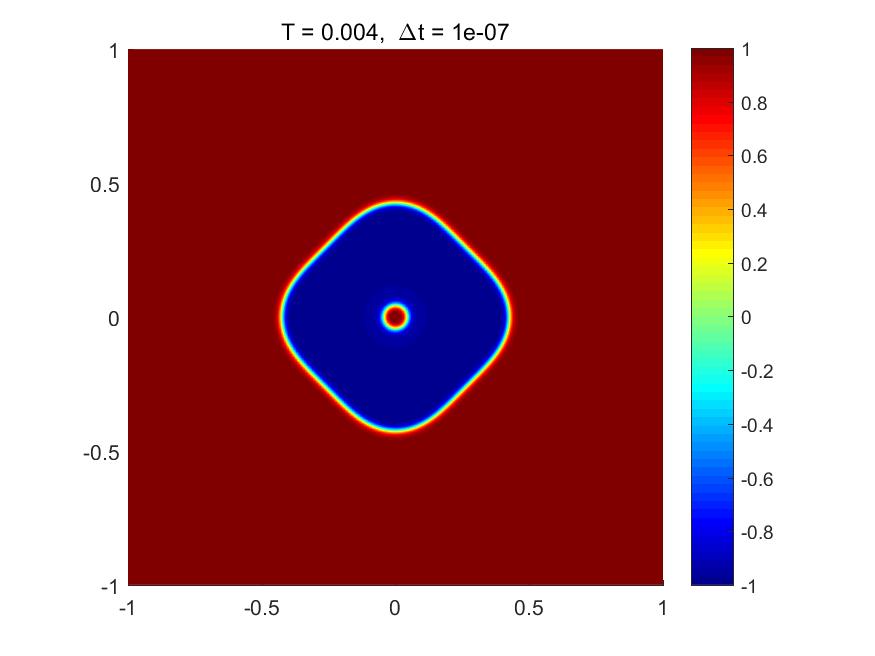}&
\includegraphics[width=5.5cm,height=4.5cm]{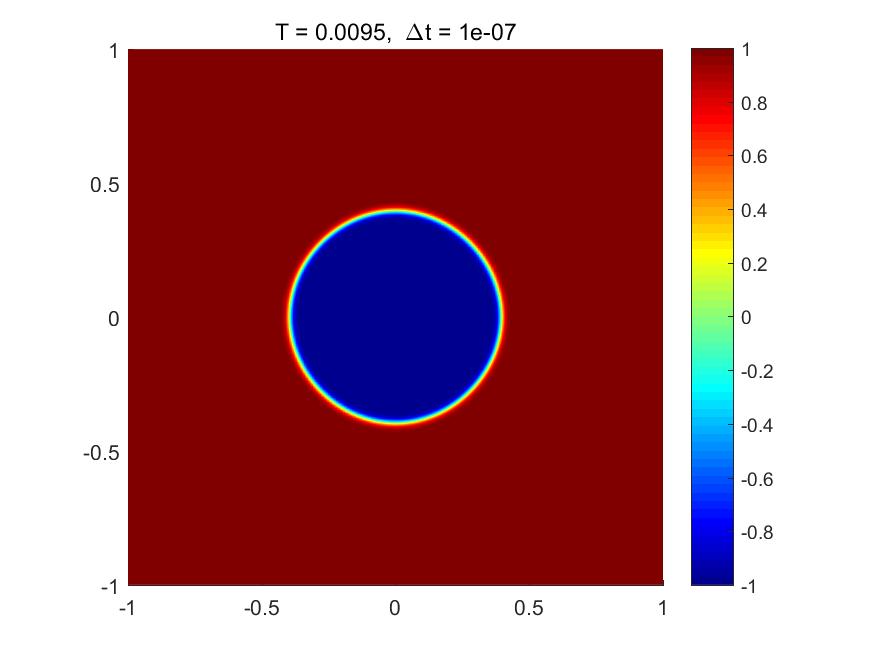}\\
\includegraphics[width=5.5cm,height=4.5cm]{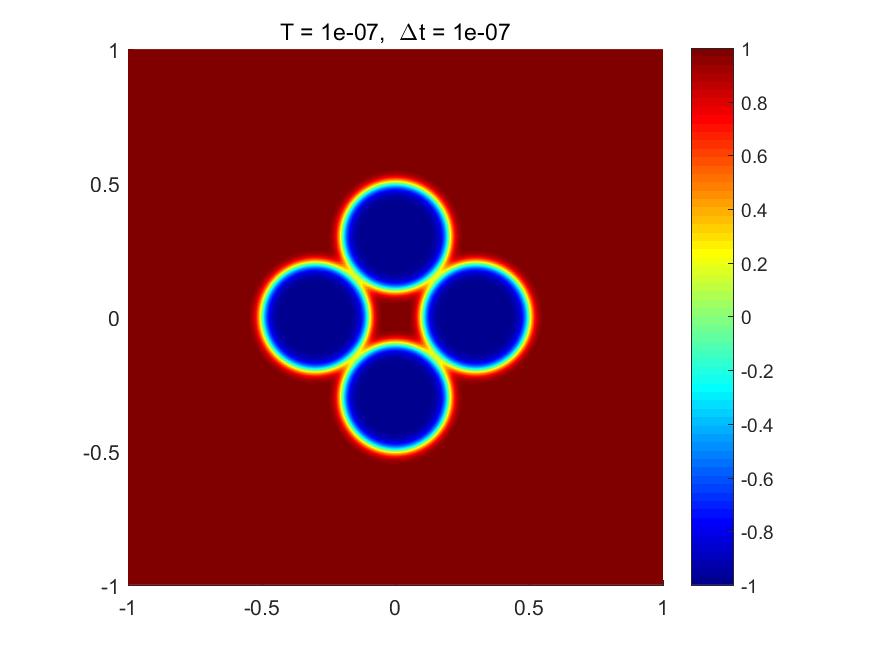}&
\includegraphics[width=5.5cm,height=4.5cm]{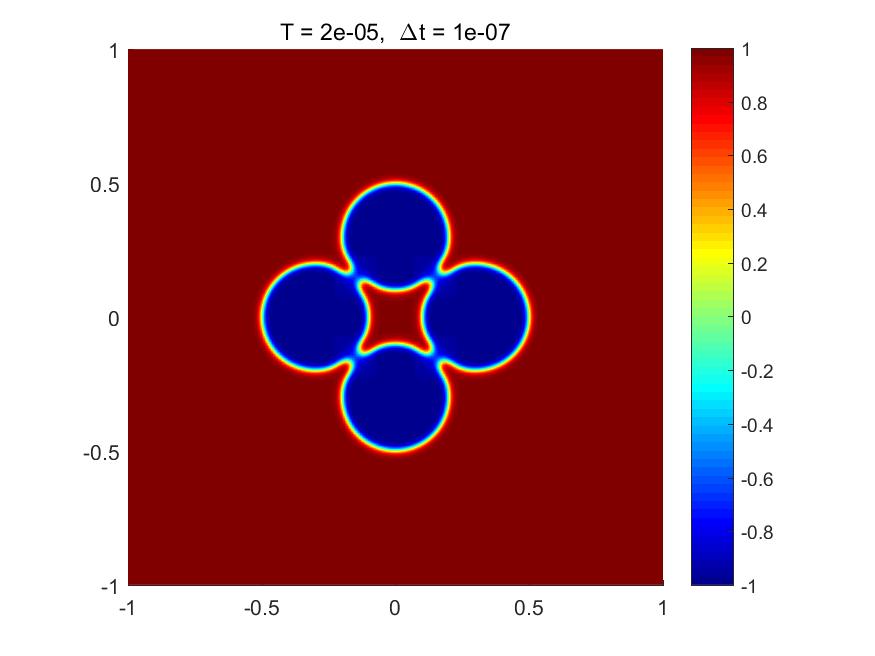}&
\includegraphics[width=5.5cm,height=4.5cm]{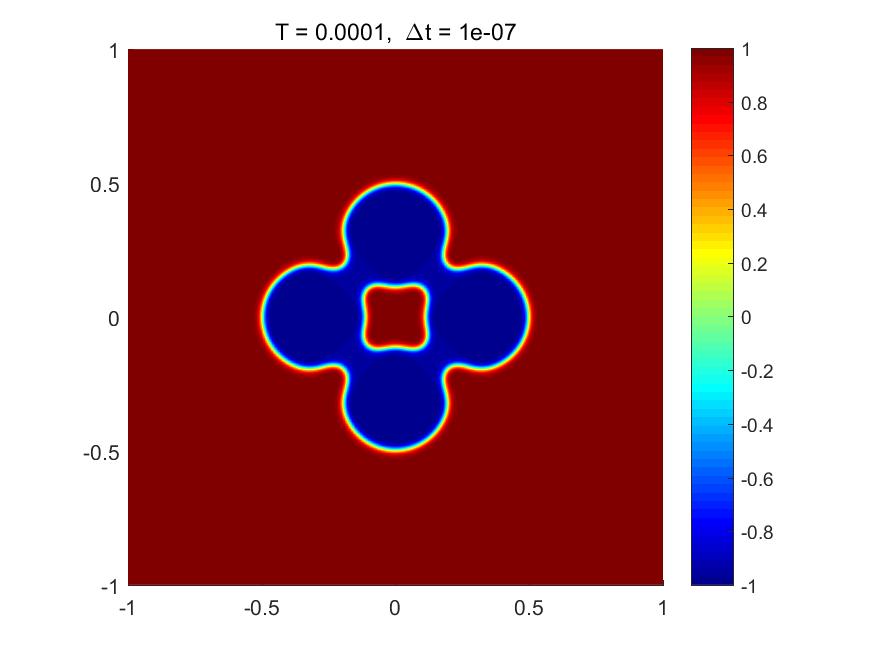}\\
\includegraphics[width=5.5cm,height=4.5cm]{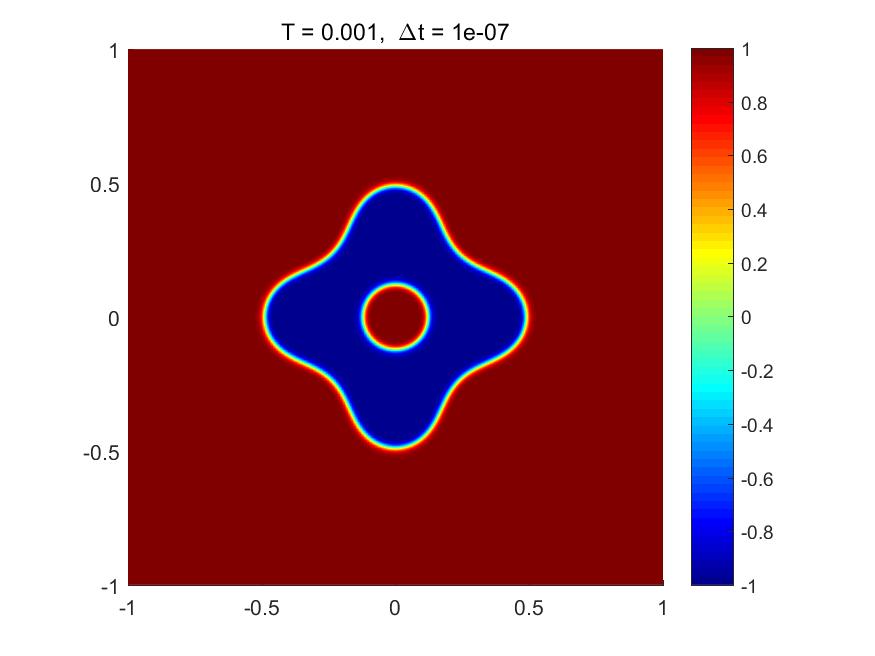}&
\includegraphics[width=5.5cm,height=4.5cm]{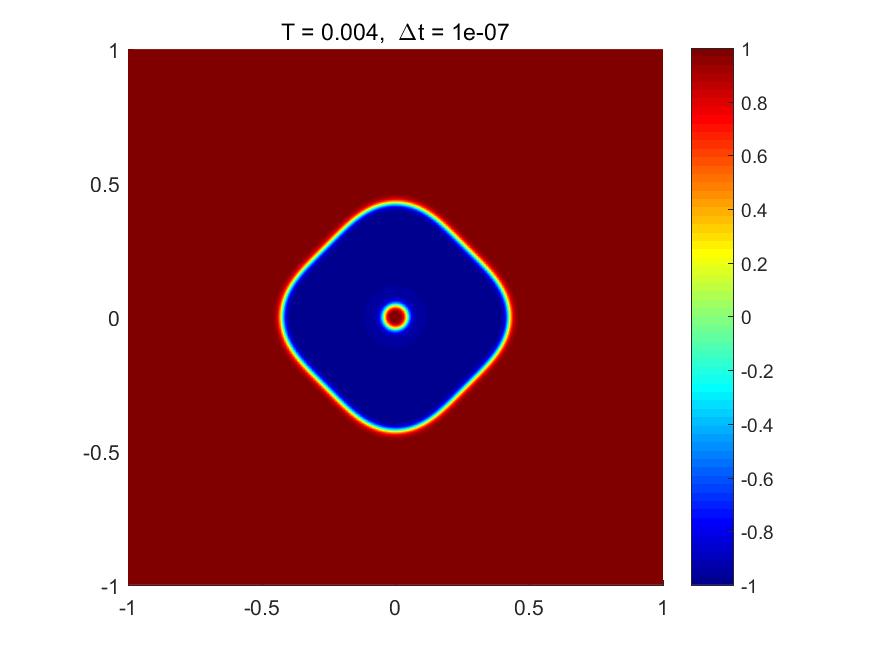}&
\includegraphics[width=5.5cm,height=4.5cm]{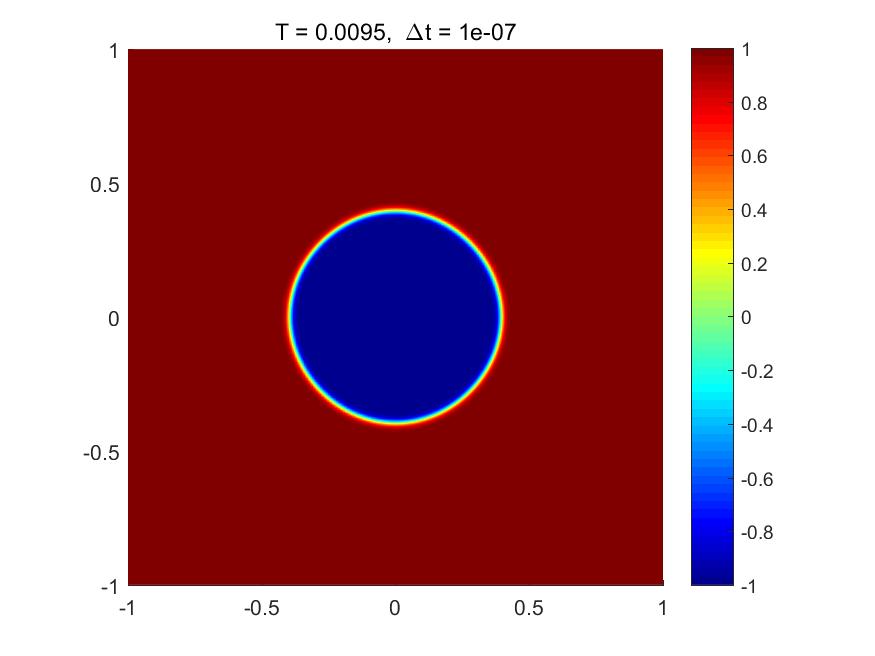}
\end{array}$\vspace{-0.2cm}
\caption{$\mathbf{Example\ \ref{cexm2-n}\ (CH)}$, numerical solutions, First and second line: BDF2-IEQ-FEM2 scheme (\ref{BDF2-IEQ-FEM2}), Third and forth line: BDF2-IEQ-FEM3 scheme (\ref{BDF2-IEQ-FEM3}).}\label{Cexp3u}
\end{figure}

In the following, we further present more examples based on Method 2 with a sufficiently large constant $B$ to investigate the method performance.
\begin{figure}[!htbp]
$\begin{array}{cc}
\includegraphics[width=5.5cm,height=4.5cm]{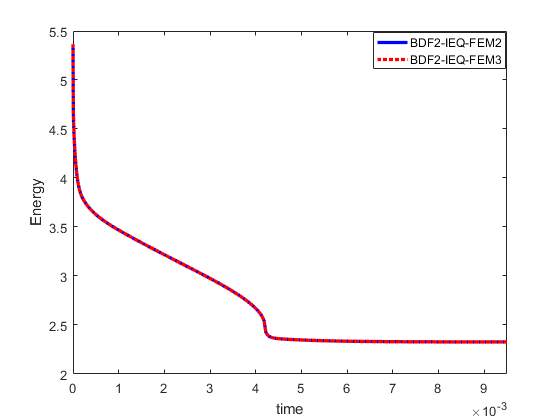}&
\includegraphics[width=5.5cm,height=4.5cm]{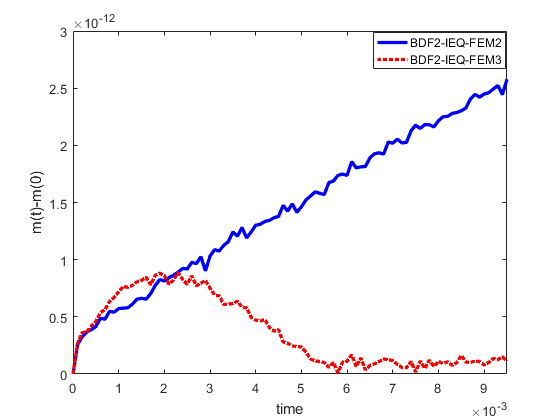}
\end{array}$\vspace{-0.2cm}
\caption{$\mathbf{Example\ \ref{cexm2-n}\ (CH)}$, Left: The discrete energy, Right: the total mass.}\label{Cexp3u-B-mass}
\end{figure}
 \end{example}
 \begin{example} \label{cexm5} \cite{LCK2017}  Consider the 3D CH equation (\ref{CH}) in $\Omega=[0,256]^{3}$ with parameter $\epsilon=\frac{4}{\sqrt{8}\cdot \tanh(0.9)}$, and  
 initial data  
 $$
 u_{0}=0.1\text{rand}(x,y,z),
 $$
 where $\text{rand}(x,y,z)$ generates random values between $-1$ and $1$.
 \begin{figure}[!htbp]
 $\begin{array}{ccc}
 \includegraphics[width=5.5cm,height=4.5cm]{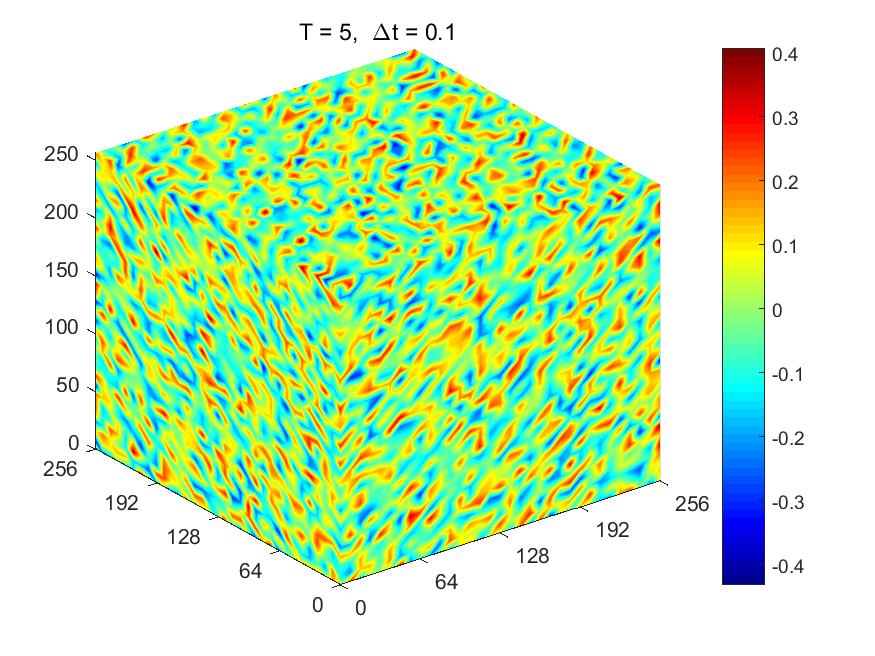}&
 \includegraphics[width=5.5cm,height=4.5cm]{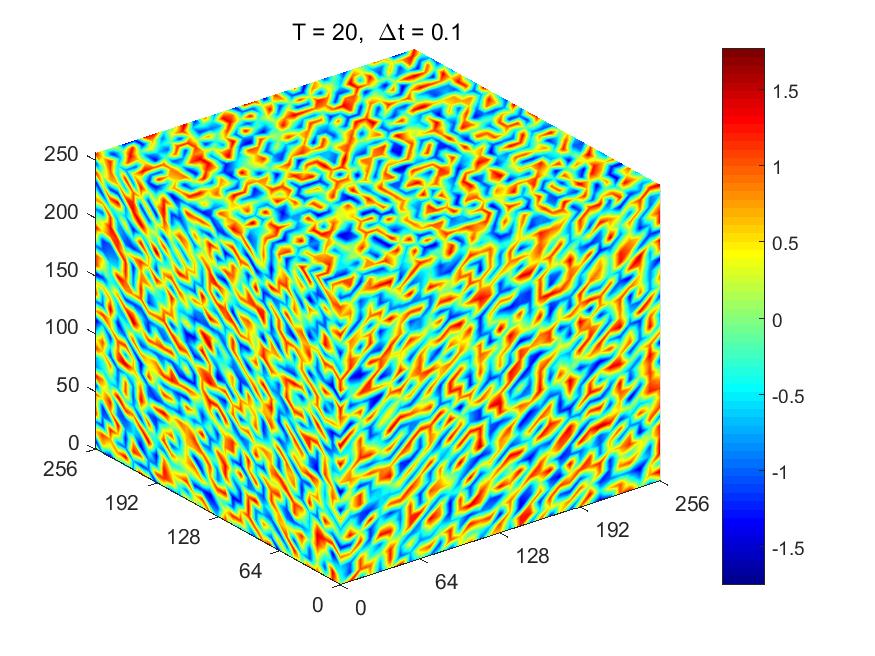}&
 \includegraphics[width=5.5cm,height=4.5cm]{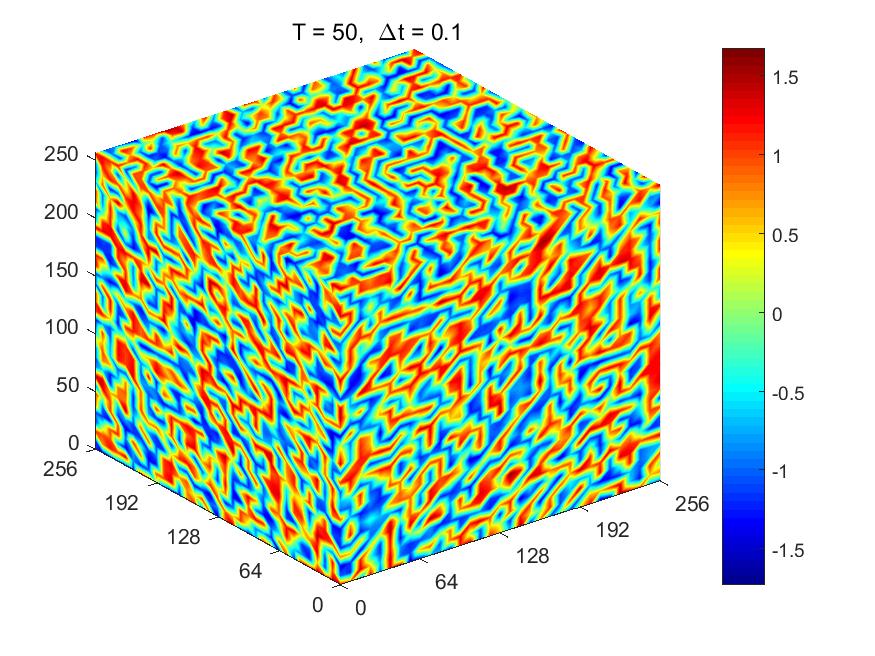}\\
 \includegraphics[width=5.5cm,height=4.5cm]{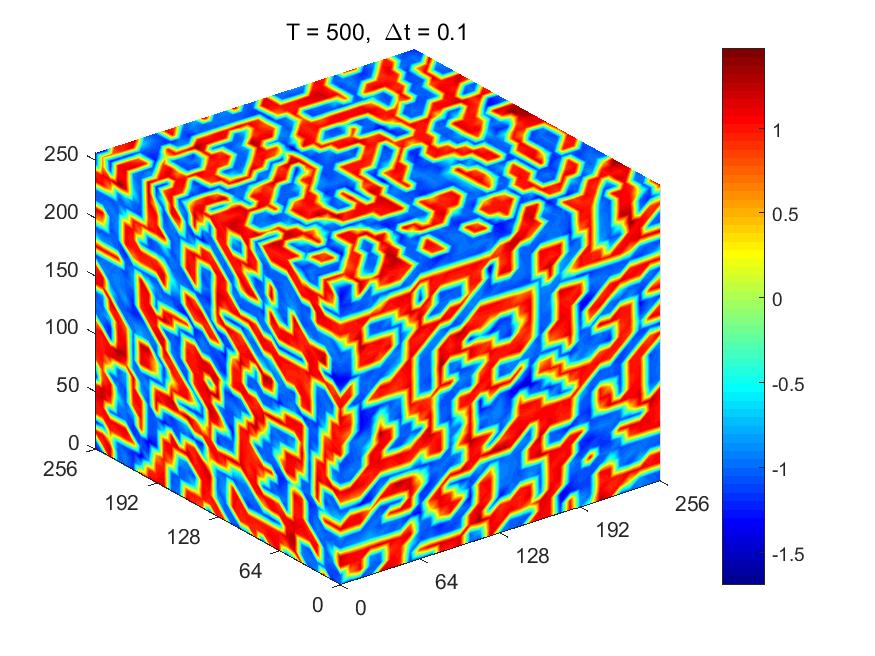}&
 \includegraphics[width=5.5cm,height=4.5cm]{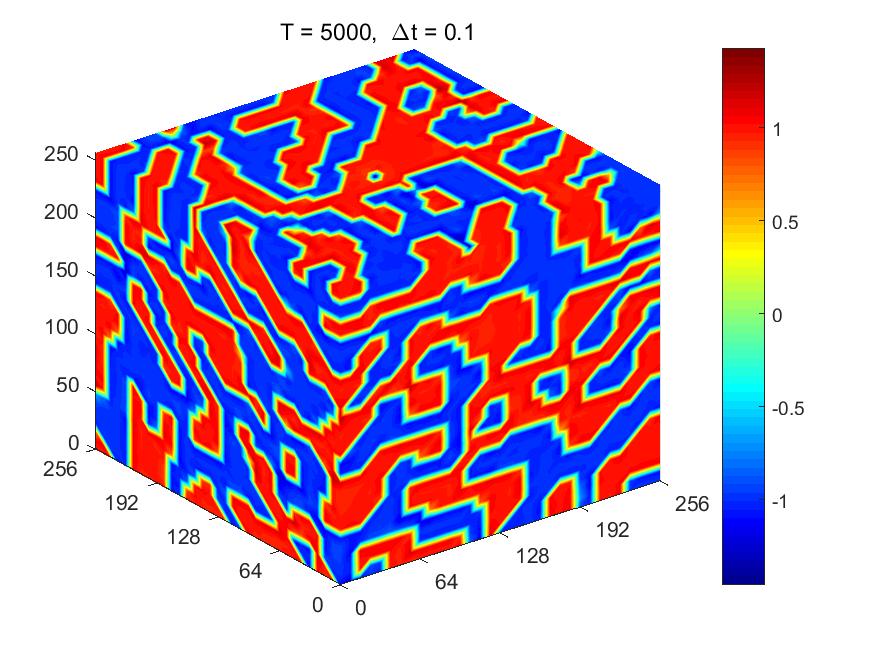}&
 \includegraphics[width=5.5cm,height=4.5cm]{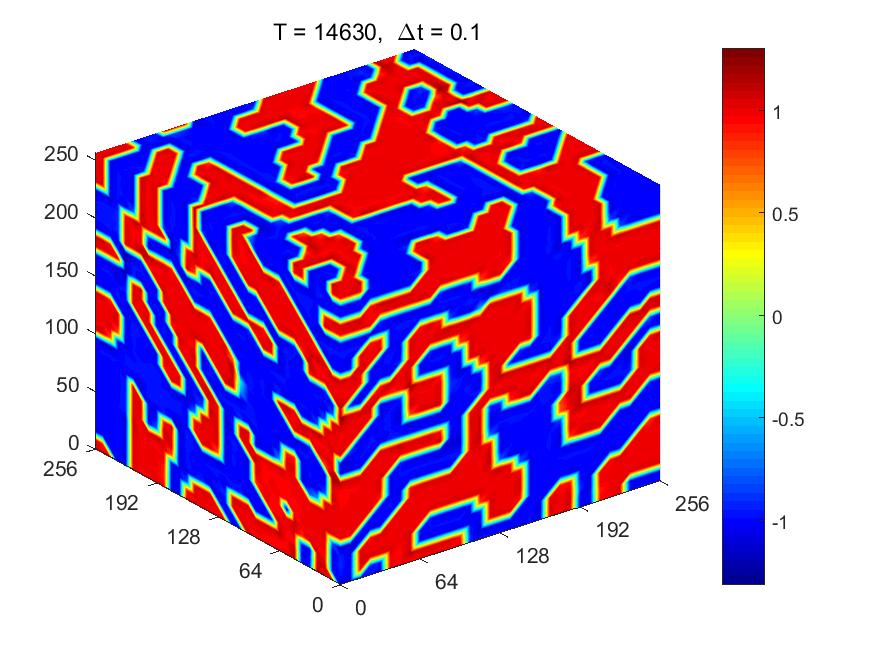}
 \end{array}$\vspace{-0.2cm}
 \caption{$\mathbf{Example\ \ref{cexm5}\ (CH3D)}$,  BDF2-IEQ-FEM2 scheme \eqref{BDF2-IEQ-FEM2}, numerical solutions.}\label{Cexp2u}
 \end{figure}


In Figure \ref{Cexp2u}, the contour plots of  solutions based on the BDF2-IEQ-FEM2 scheme \eqref{BDF2-IEQ-FEM2} with $B=20000,\, \Delta t=0.1$ are presented. The observed coarsening phenomenon indicates that the proposed method is effective for solving the 3D CH equation.

Finally, the graphs depicting the evolution of discrete energy and mass for $\mathbf{Example\ \ref{cexm5}}$ are presented in Figure \ref{energyCH}. As illustrated, the discrete energy decreases over time, and the mass is conserved.

\begin{figure}[!htbp]
$\begin{array}{ccc}
\includegraphics[width=5.5cm,height=4.5cm]{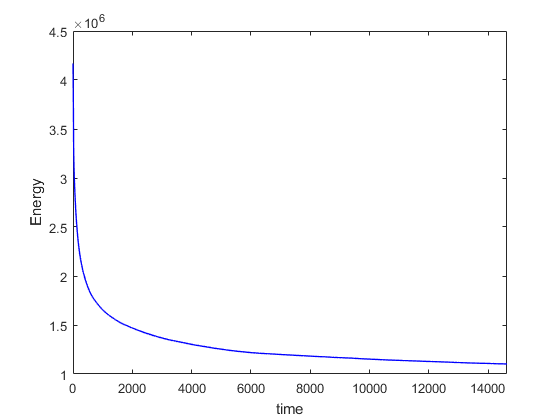} & 
\includegraphics[width=5.5cm,height=4.5cm]{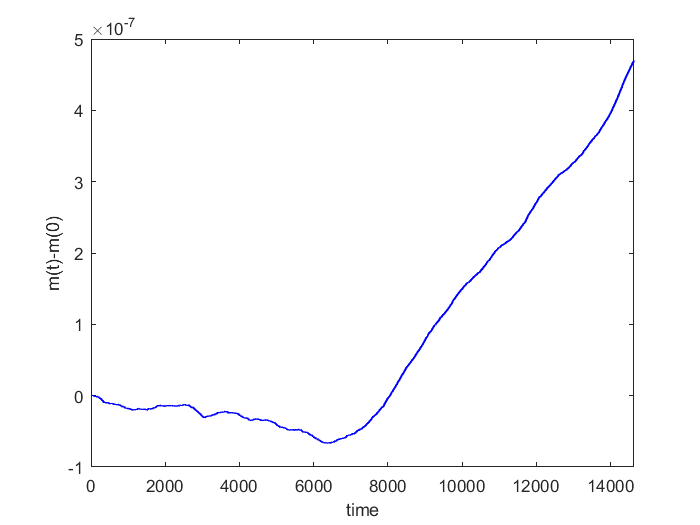}
\end{array}$\vspace{-0.2cm}
\caption{$\mathbf{Example\ \ref{cexm5}}$. Left: The discrete energy; Right: the total mass.}\label{energyCH}
\end{figure}

\end{example}

\subsection{AC equation}

In this section, we will provide several examples to validate the proposed methods for the AC equation \eqref{AC} and the scaled AC equation \eqref{AC1}, including the tests for temporal and spatial accuracy, a comparison of the first-order schemes, and the performance of the proposed IEQ-FEMs.

\subsubsection{Convergent rates}
 In this part, convergent rates of the proposed methods for the AC equation (\ref{2e1}) are reported in $\mathbf{Example\ \ref{2e2a}}$.
\begin{example}\label{2e2a}
In this example, we examine the AC equation (\ref{2e1}) in $\Omega=[-1,1]^2$ with parameters $\epsilon=1$. The exact solution is given by $u=e^{t}\cos(\pi x)\cos(\pi y)$, and the corresponding right term $s(x,t)$ can be determined by substituting $u$ into equation \eqref{2e1}.

 \begin{figure}[!htbp]
 $\begin{array}{cccc}
 \includegraphics[width=5.5cm,height=4.5cm]{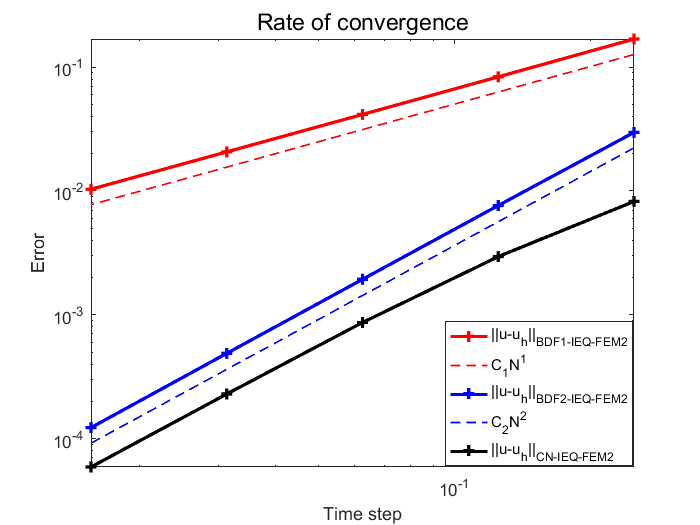}
 \includegraphics[width=5.5cm,height=4.5cm]{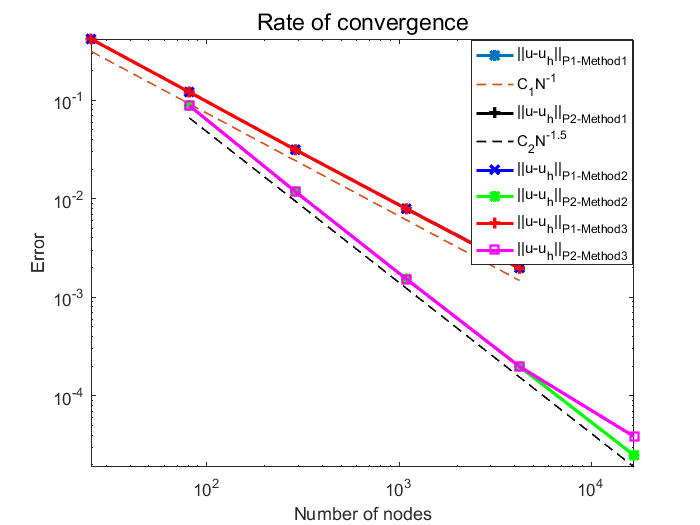}
 \end{array}$
\caption{$\mathbf{Example\ \ref{2e2a}}$, $L_2$ errors and rates of convergence for AC equation, Left: Time accuracy test ($P_{2}$ element, $T=1$, uniform mesh with $(h_x, h_y)=(\frac{4 \pi}{200},\frac{4 \pi}{200})$); Right: Spatial accuracy test ($\Delta t=10^{-6},\, T=10^{-3}$).}\label{convergentrateAC}
 \end{figure}

For the involved IEQ-FEMs, we take the constant $B=1$.
The time and spatial convergent rates of the proposed methods for the AC equation are shown in Figure \ref{convergentrateAC}. From the pictures, we can see that the time accuracy is first order for the BDF1-FEM-FEM2 scheme (\ref{BDF1-FEM2}), and second order for both the CN-FEM-FEM2 scheme (\ref{CN-FEM2}) and the BDF2-FEM-FEM2 scheme (\ref{BDF2-FEM2}). 
For spatial accuracy, the IEQ-FEMs exhibit second order accuracy for $P_1$ approximations and third order accuracy for $P_2$ approximations.

\end{example}
\subsubsection{Comparison between different schemes}
In this part, we compare the first order fully discrete BDF1-IEQ-FEMs (\ref{BDF1-FEM1}), \eqref{BDF1-FEM2}, and \eqref{BDF1-FEM3}.
\begin{example}\label{exmS1} \cite{CHY2019} We consider the AC equation $(\ref{AC1})$ in $\Omega=[-2,2]^{2}$ with $\epsilon=\frac{1}{16}$. Define $m_{1}=(0,2)$, $m_{2}=(0,0)$, $m_{3}=(0,-2)$. Let $r_{1}=r_{3}=2-3\epsilon/2 $, $r_{2}=1$ and set $d_{j}(x)=|x-m_{j}|-r_{j}$ for $x\in \Omega$ and $j=1,2,3$. Then we take the initial condition
\begin{center}
\bq\label{4e1a}
u_{0}(x,y)=-\tanh\left( \frac{d(x)}{\sqrt{2}\epsilon}\right),\ d(x)=\max\{-d_{1}(x),d_{2}(x),-d_{3}(x)\}.
\eq
\end{center}

\begin{figure}[!htbp]
$\begin{array}{cccc}
\includegraphics[width=5.5cm,height=4.5cm]{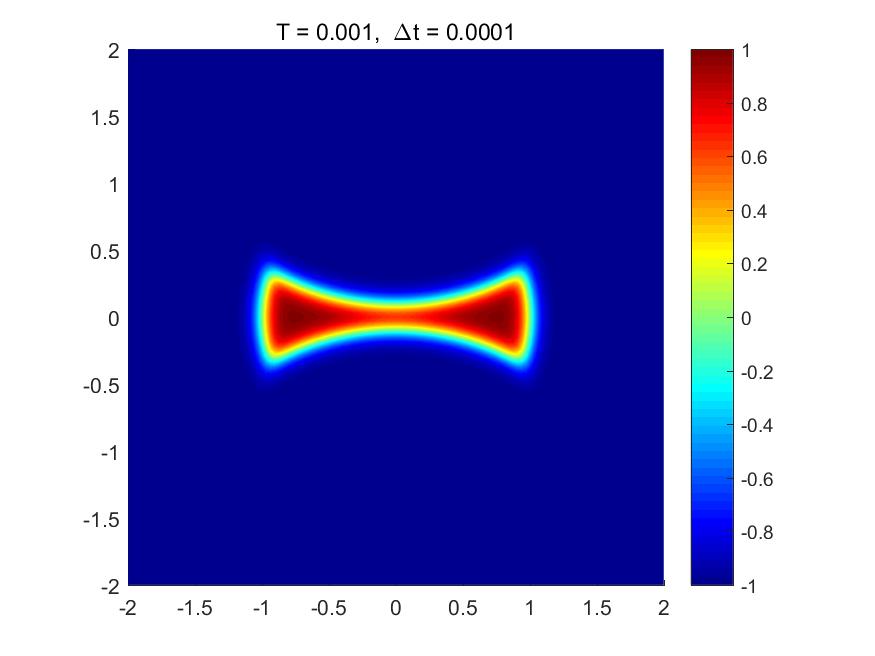}&
\includegraphics[width=5.5cm,height=4.5cm]{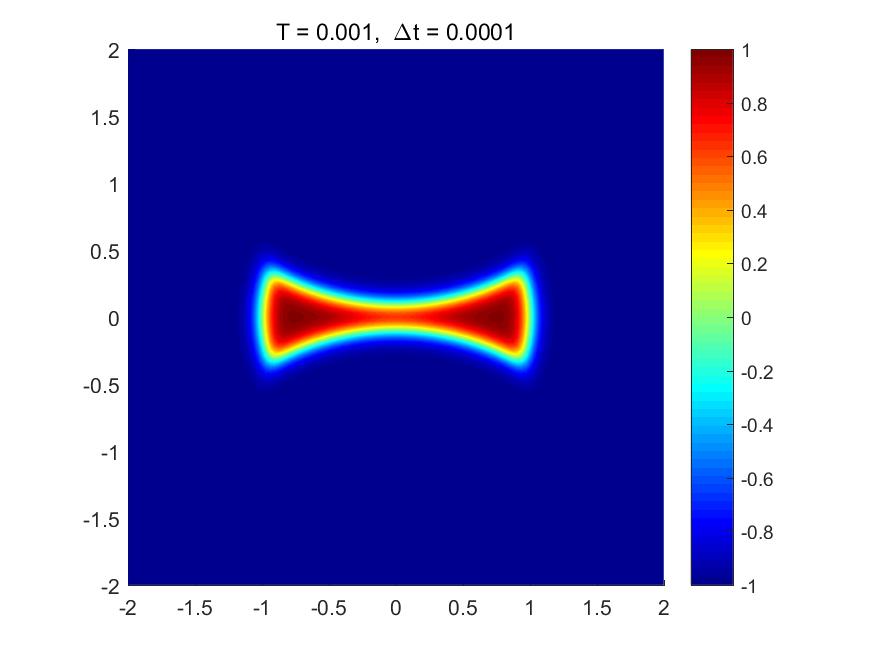}&
\includegraphics[width=5.5cm,height=4.5cm]{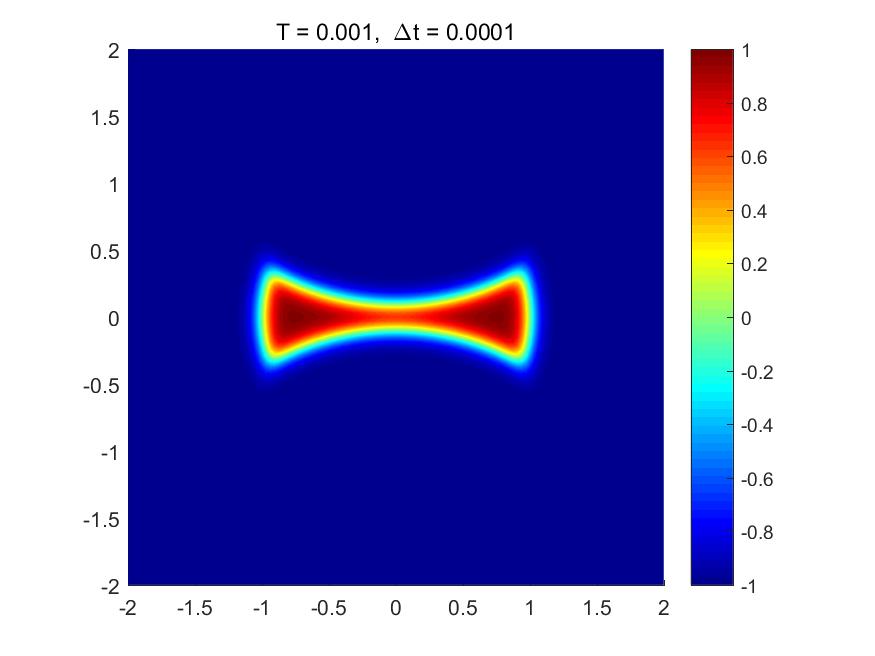}\\\vspace{-0.3cm}
\includegraphics[width=5.5cm,height=4.5cm]{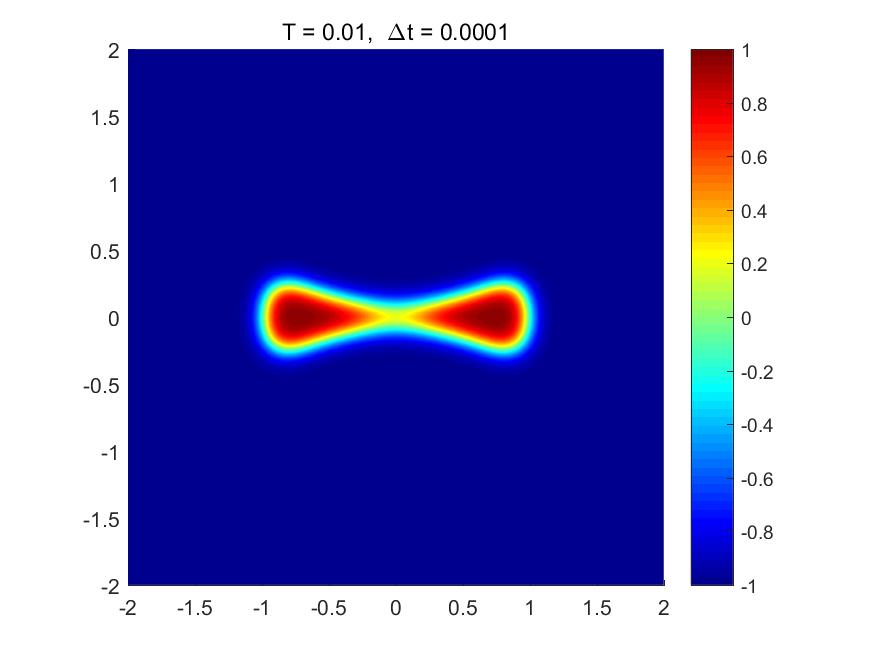}&
\includegraphics[width=5.5cm,height=4.5cm]{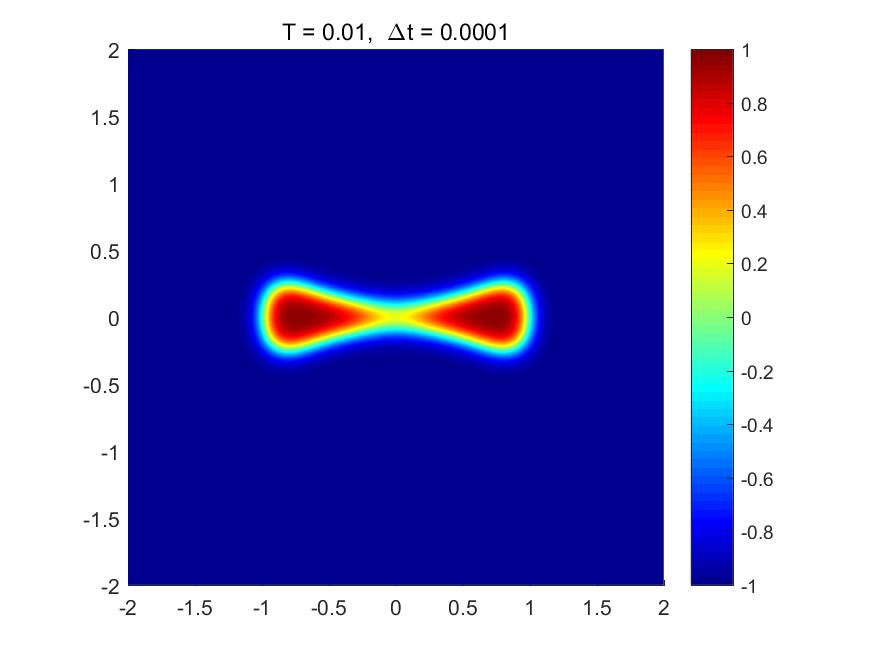}&
\includegraphics[width=5.5cm,height=4.5cm]{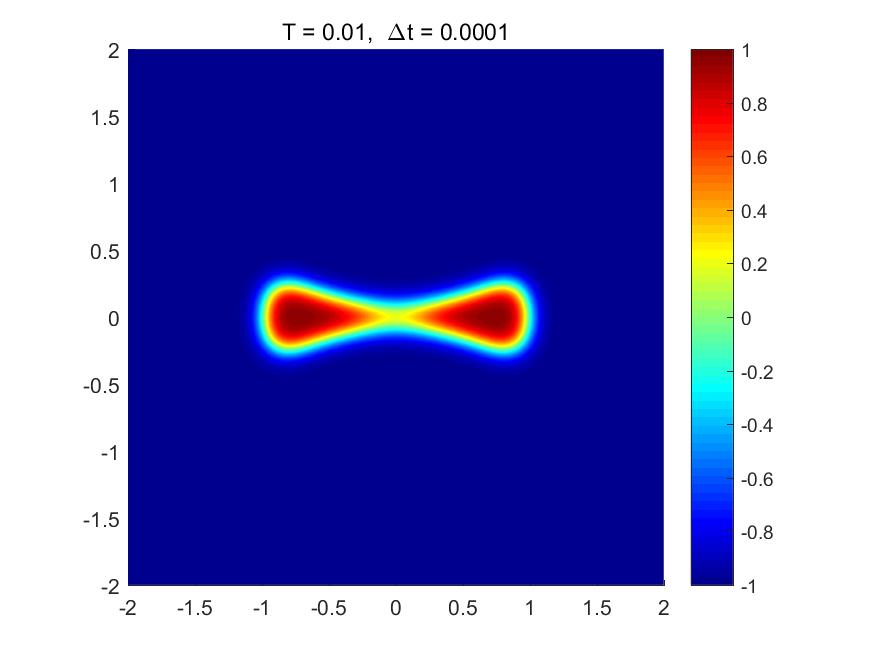}&\\\vspace{-0.3cm}
\includegraphics[width=5.5cm,height=4.5cm]{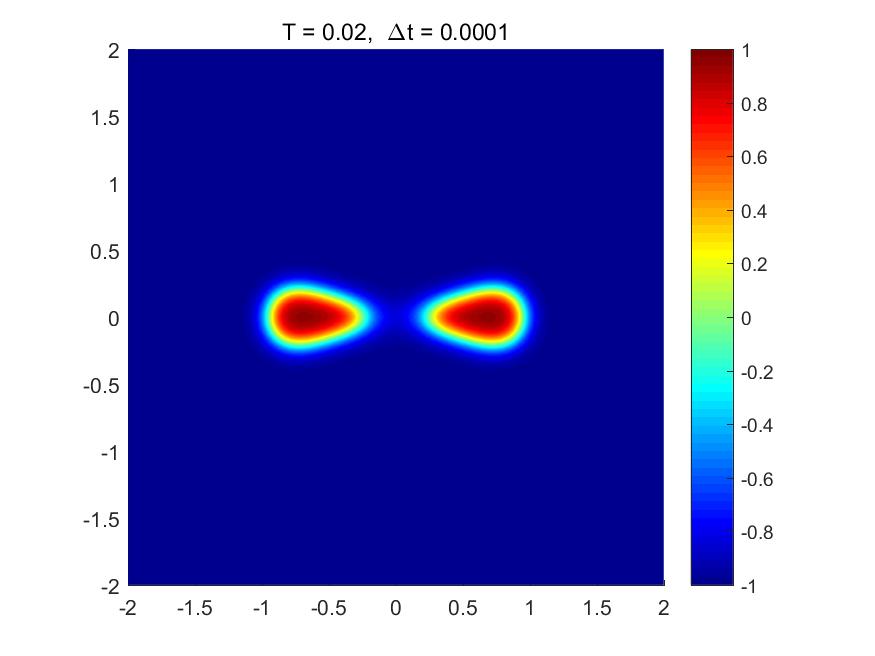}&
\includegraphics[width=5.5cm,height=4.5cm]{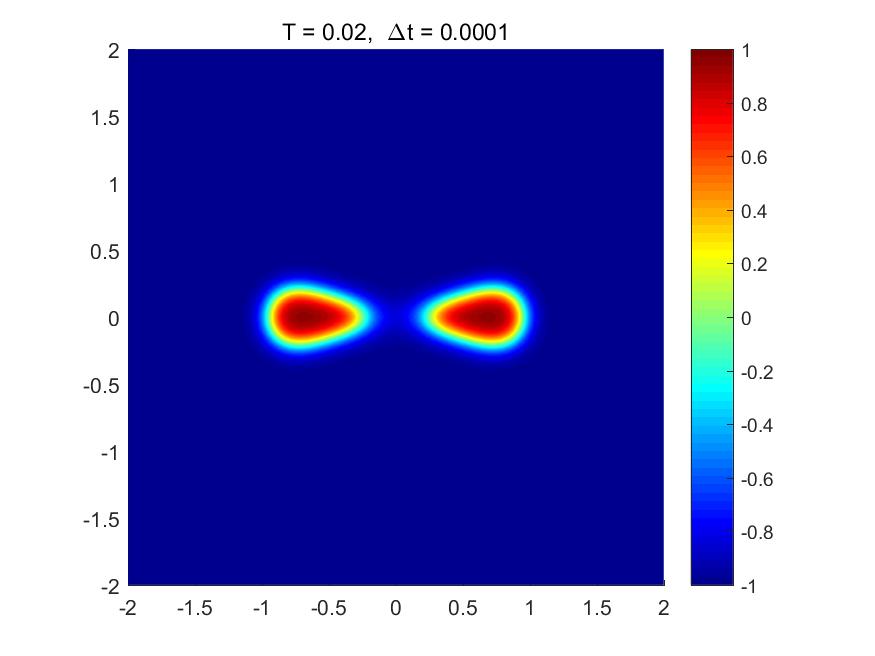}&
\includegraphics[width=5.5cm,height=4.5cm]{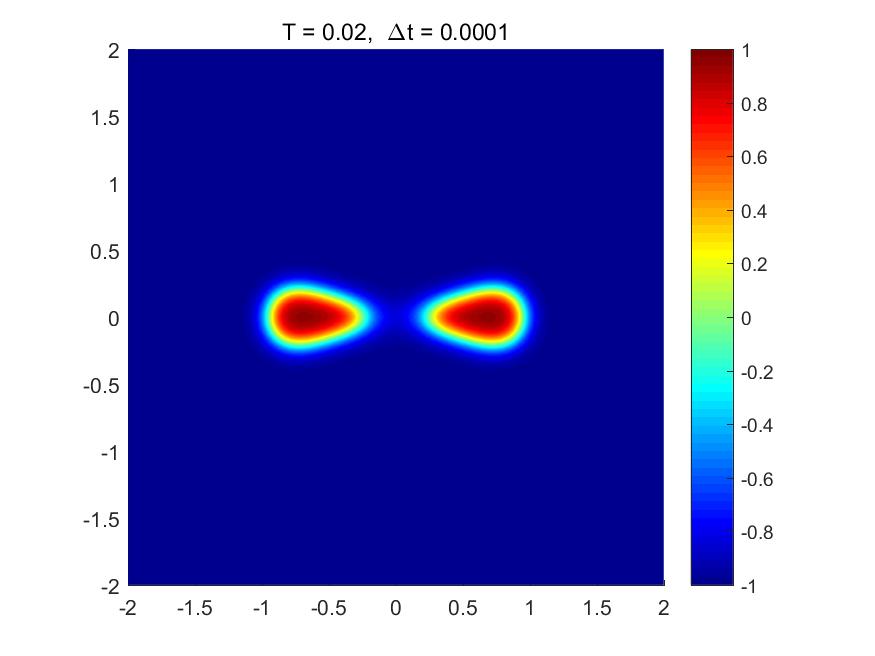}\\\vspace{-0.3cm}
\includegraphics[width=5.5cm,height=4.5cm]{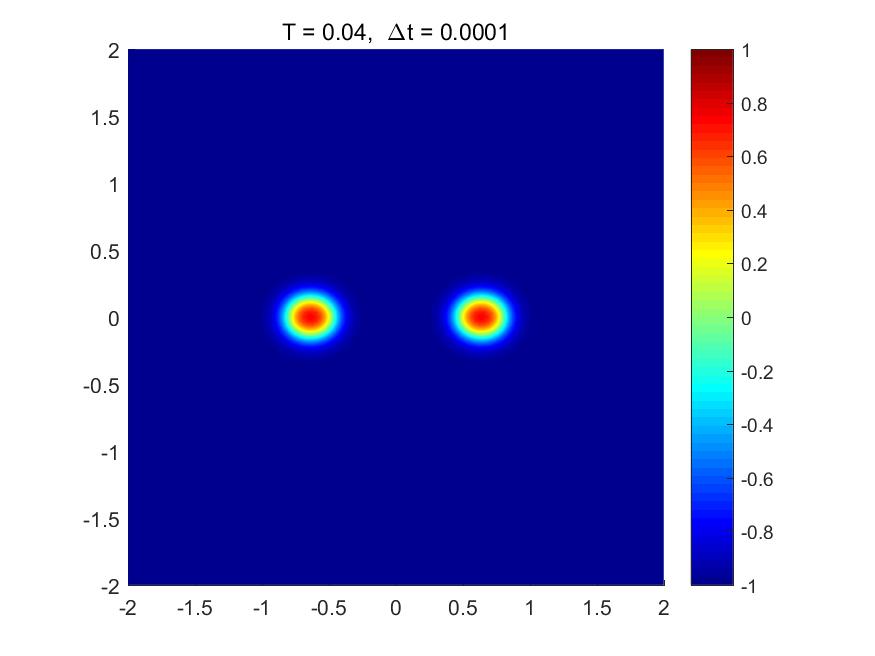}&
\includegraphics[width=5.5cm,height=4.5cm]{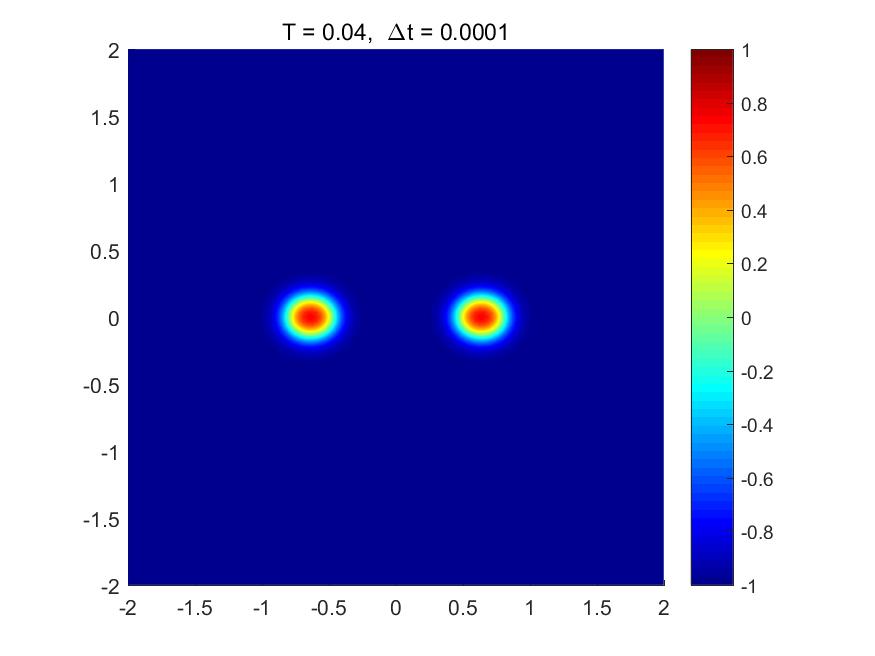}&
\includegraphics[width=5.5cm,height=4.5cm]{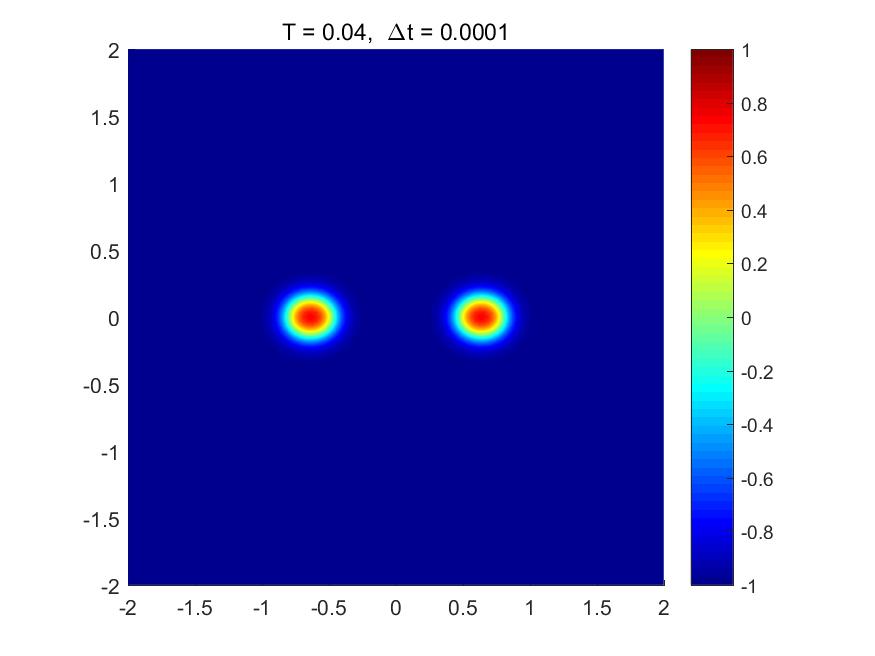}&\\\vspace{-0.3cm}
\includegraphics[width=5.5cm,height=4.5cm]{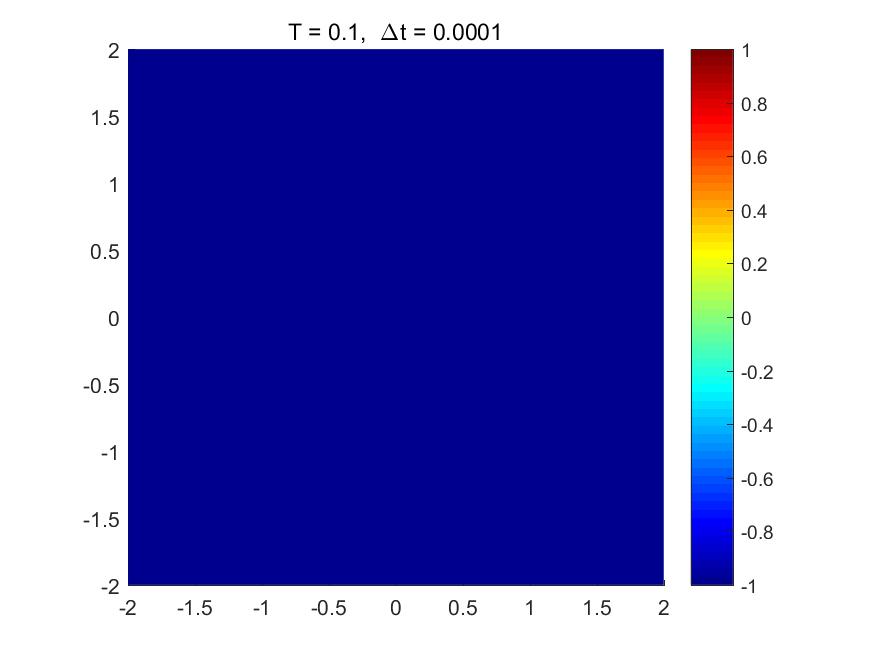}&
\includegraphics[width=5.5cm,height=4.5cm]{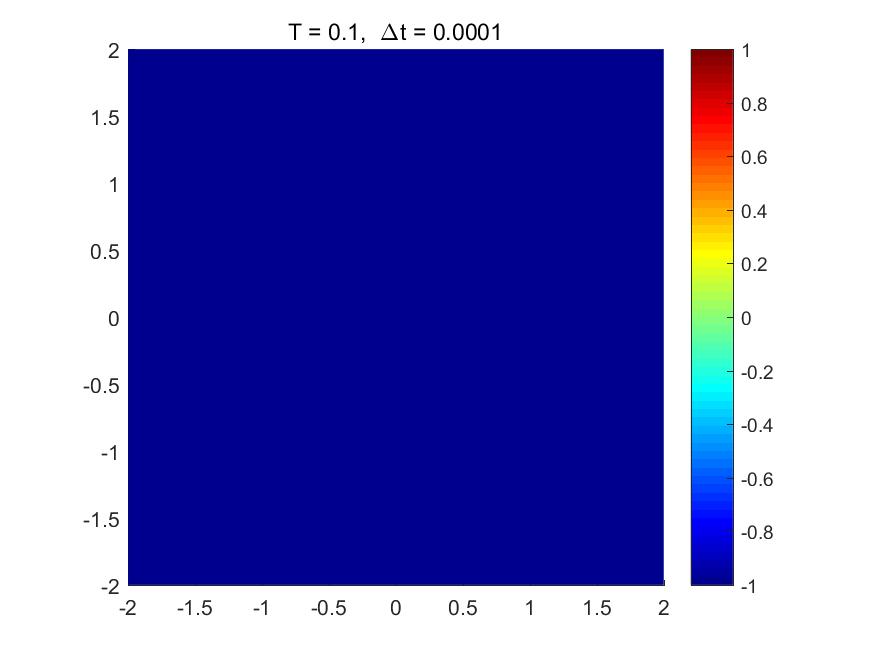}&
\includegraphics[width=5.5cm,height=4.5cm]{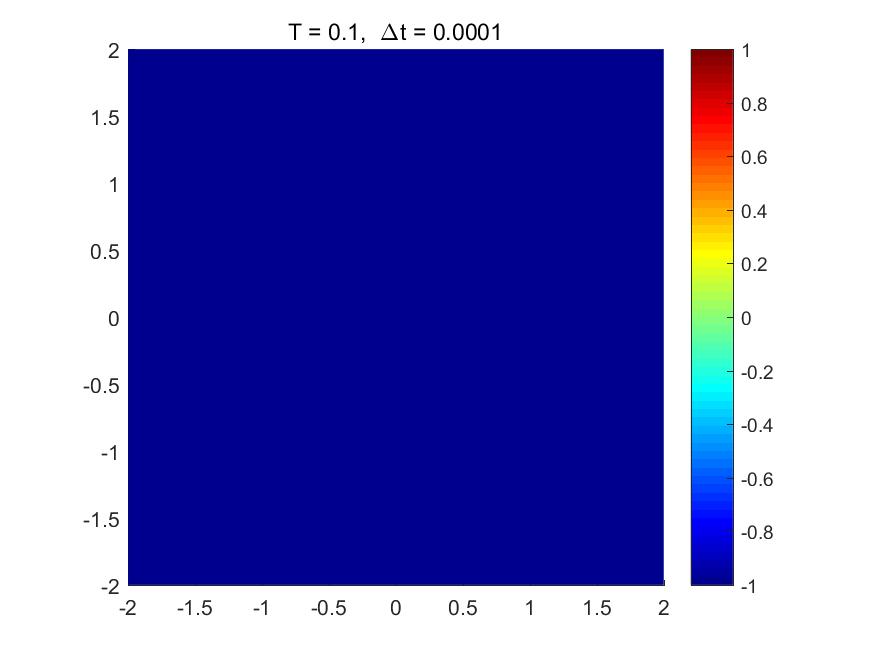}&
\end{array}$\vspace{-0.3cm}
\caption{Contour plots of the IEQ-FEM solutions in $\mathbf{Example\ \ref{exmS1}\ (AC)}$, Left: BDF1-IEQ-FEM1 scheme (\ref{BDF1-FEM1}); Middle: BDF1-IEQ-FEM2 scheme (\ref{BDF1-FEM2}); Right: BDF1-IEQ-FEM3 scheme (\ref{BDF1-FEM3}).}\label{expsu1}
\end{figure}

\begin{figure}[!htbp]
$\begin{array}{cccc}
\includegraphics[width=5.5cm,height=4.5cm]{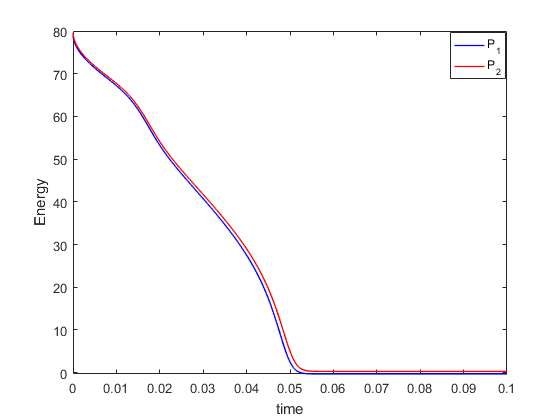}
\includegraphics[width=5.5cm,height=4.5cm]{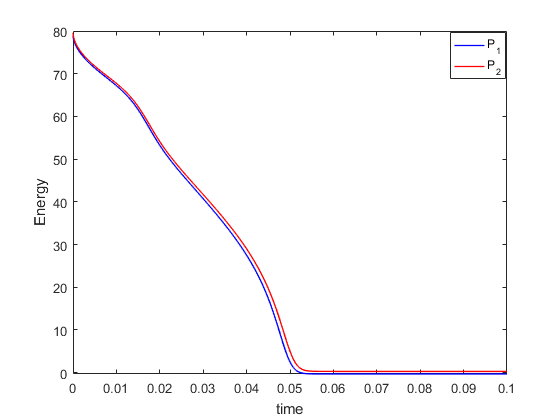}
\includegraphics[width=5.5cm,height=4.5cm]{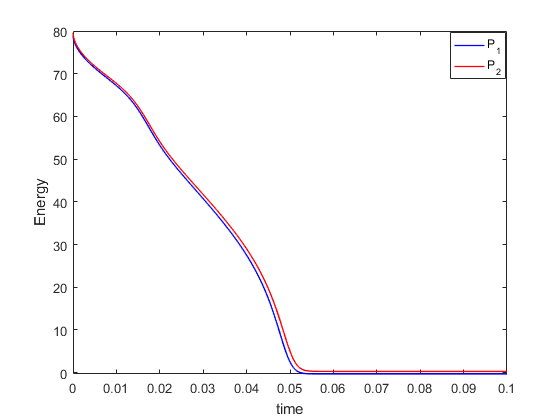}
\end{array}$
\caption{$\mathbf{Example\ \ref{exmS1}\ (AC, T=0.1)}$, Time evolution of the discrete energy, Left: BDF1-IEQ-FEM1 scheme (\ref{BDF1-FEM1}); Middle: BDF1-IEQ-FEM2 scheme (\ref{BDF1-FEM2}); Right: BDF1-IEQ-FEM3 scheme (\ref{BDF1-FEM3}).}\label{energuAC}
\end{figure}


\begin{table}[!ht]
\centering
\caption{$\mathbf{Example\ \ref{exmS1}\ (AC, T=0.1)}$, CPU time (in seconds) of three kinds of methods (11th Gen Intel(R) Core(TM) i5-1135G7 @ 2.40GHz 2.42GHz). (`$--$' represents running out of memory)}\label{Table1:S1}
\begin{tabular}{| c | c | c | c | c |}
\hline
Methods    &  scheme (\ref{BDF1-FEM1})  &scheme (\ref{BDF1-FEM2})  & scheme (\ref{BDF1-FEM3})                   \\  \hline
$P_{1}$   &   277s    &  179s    &    219s                            \\  \hline
$P_{2}$  &  1435s  &   770s     &    828s                             \\  \hline
\end{tabular}
\end{table}

We present contour plots of the corresponding approximate solutions generated by the first-order fully discrete BDF1-IEQ-FEMs (\ref{BDF1-FEM1}),  (\ref{BDF1-FEM2}), and (\ref{BDF1-FEM3}) with $\Delta t=10^{-4}$, $B=1$ in \Cref{expsu1}. 
The initially connected interface undergoes a split into two parts, and subsequently, the two components of the interface evolve into circular shapes. The diameters of the two particles decrease to zero until they collapse. The corresponding pictures of discrete energy are displayed in Figure \ref{energuAC}. These images show a close similarity, indicating that all three methods perform effectively.

The CPU time for these schemes is provided in Table \ref{Table1:S1}. It is evident that the CPU time for the BDF1-IEQ-FEM2 scheme (\ref{BDF1-FEM2}) and the BDF1-IEQ-FEM3 scheme (\ref{BDF1-FEM3}) is shorter than that for the BDF1-IEQ-FEM1 scheme (\ref{BDF1-FEM1}).
Similarly, we also find that solving BDF1-IEQ-FEM1 scheme (\ref{BDF1-FEM1}) by evaluating the inverse of the mass matrix \eqref{AC-BDF1-IEQ-FEM1-mat1} is more expensive than solving the augmented system \eqref{AC-BDF1-IEQ-FEM1-mat}.
Therefore, the focus in the subsequent tests will be primarily on Method 2 and Method 3.
\end{example}
\subsubsection{Numerical solving for the AC equation}
\begin{example}\label{cexmAC2-n} In this example, we consider the AC equation (\ref{AC})
with $\Omega=[-1,1]^{2} $, $\epsilon=0.02$, $\Delta t=10^{-7}$, and the initial condition
\[\begin{aligned}
u_{0}=\left\{\begin{aligned}
& \qquad 1,\qquad\quad \qquad x< x_0,\\
&\quad -1,\qquad\quad \qquad x>x_1,\\
&-\sin\left(\frac{\pi x}{2 x_1} \right),  \qquad x_0\leq x \leq x_1,
\end{aligned}\right.
\end{aligned}\]
where $x_1=-x_0=\frac{\sqrt{2}}{20}$.

 \noindent\textbf{Test Case 1.} 
 We solve this problem by BDF1-IEQ-FEM2 and BDF1-IEQ-FEM3 based on $P_1$ polynomials with mesh size $h=0.1$ but different constants $B$. The evolution of the discrete energy $E(u_h^n, \Pi U^n)$ for BDF1-IEQ-FEM2 and $E(u_h^n, U_h^n)$ for BDF1-IEQ-FEM3 is reported in \Cref{Cexp8u-B}, from which we can find that the solution of BDF1-IEQ-FEM2 only satisfies the energy dissipation law when $B$ is appropriately large, as explained in \Cref{Enerdis}, while that of BDF1-IEQ-FEM3 satisfies the energy dissipation law independent of the choice of the constant $B$.

 \begin{figure}[!htbp]
 $\begin{array}{ccc}
 \includegraphics[width=5.5cm,height=4.5cm]{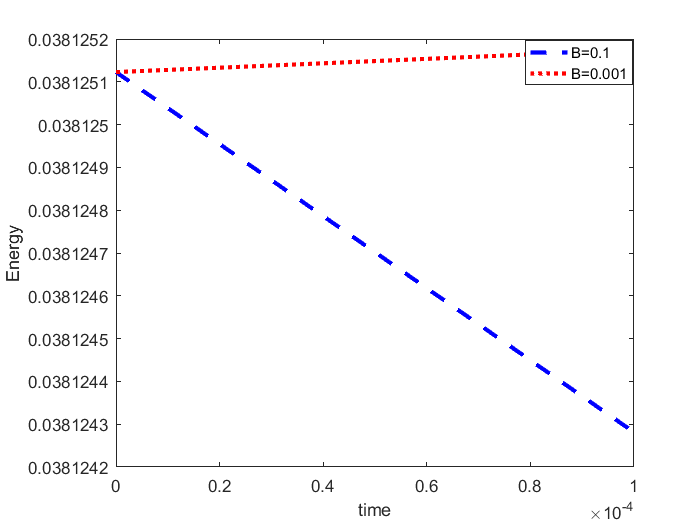}&
 \includegraphics[width=5.5cm,height=4.5cm]{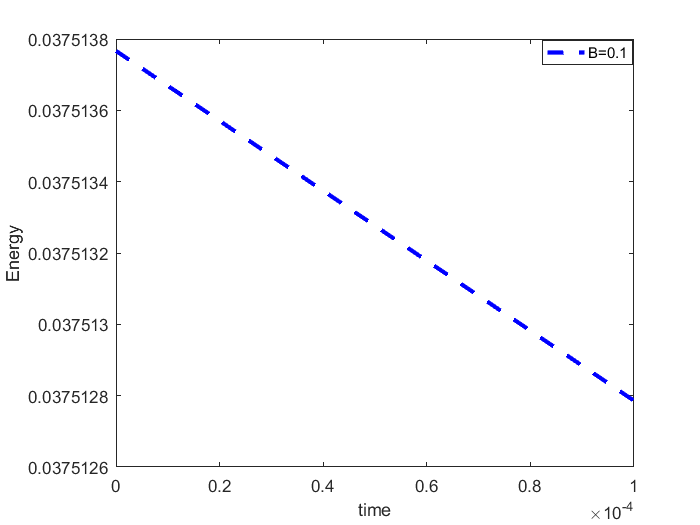}&
 \includegraphics[width=5.5cm,height=4.5cm]{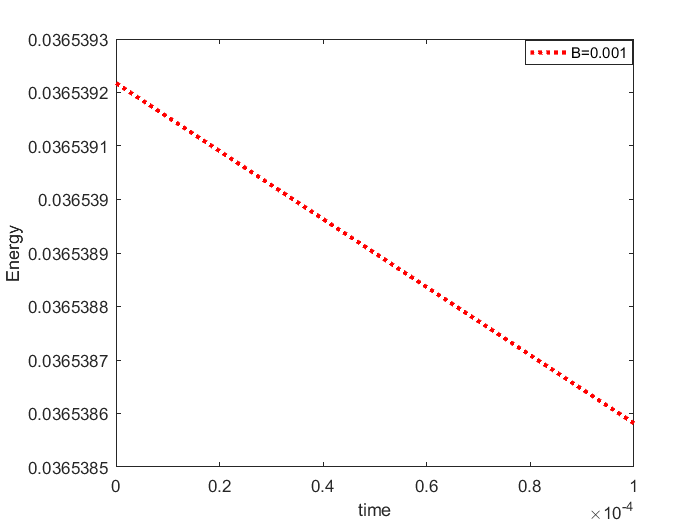}\\
\end{array}$\vspace{-0.2cm}
\caption{$\mathbf{Example\ \ref{cexmAC2-n}\ (AC)}$, Effect of constant $B$ on the discrete energy, Left: BDF1-IEQ-FEM2 scheme \eqref{BDF1-FEM2}, Middle: BDF1-IEQ-FEM3 scheme \eqref{BDF1-FEM3} with $B=0.1$, Right: BDF1-IEQ-FEM3 scheme \eqref{BDF1-FEM3} with $B=0.001$.}\label{Cexp8u-B}
 \end{figure}

 \noindent\textbf{Test Case 2.} 
 We continue solving this problem by BDF1-IEQ-FEM2 based on $P_1$ and $P_2$ elements with the constant $B=0.001$ but different mesh sizes $h$. The evolutions of the discrete energy $E(u_h^n, \Pi U^n)$ are shown in \Cref{Cexp8u-hB}, from which we can find that the solution of BDF1-IEQ-FEM2 satisfies the energy dissipation law by using refined meshes as explained in \Cref{Enerdis}.

 \begin{figure}[!htbp]
 $\begin{array}{cc}
 \includegraphics[width=5.5cm,height=4.5cm]{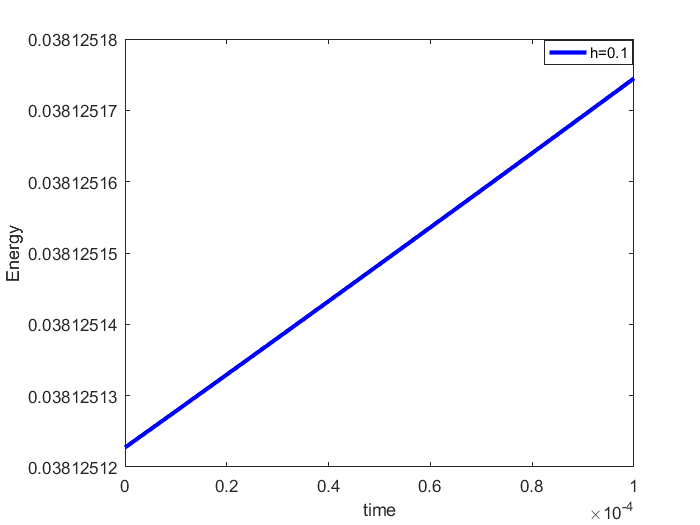}&
 \includegraphics[width=5.5cm,height=4.5cm]{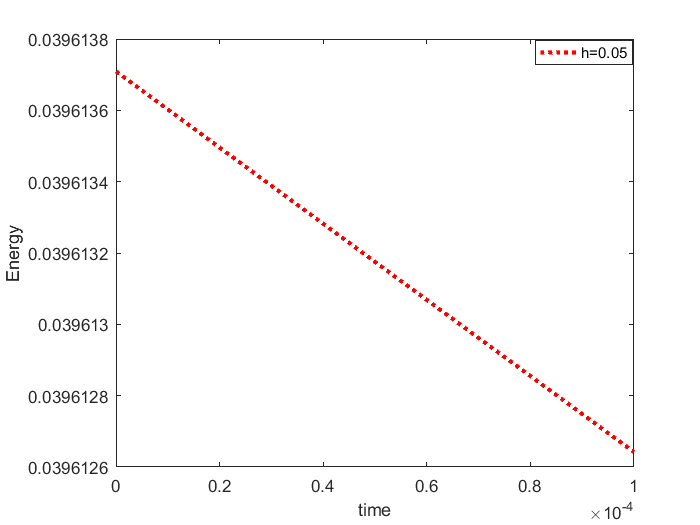}\\
 \includegraphics[width=5.5cm,height=4.5cm]{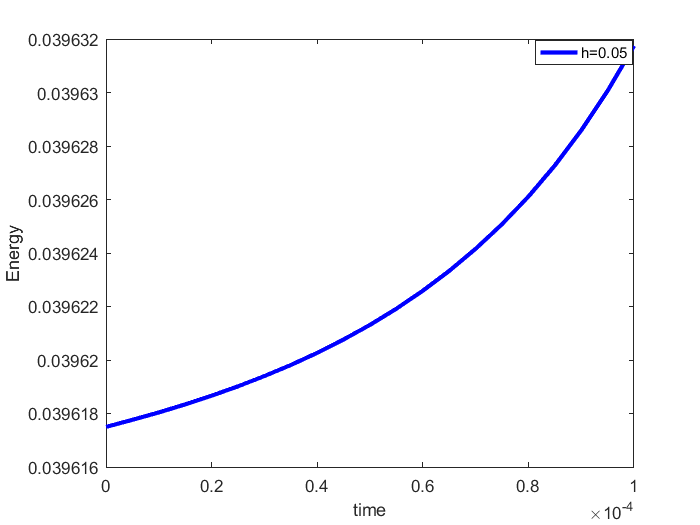}&
  \includegraphics[width=5.5cm,height=4.5cm]{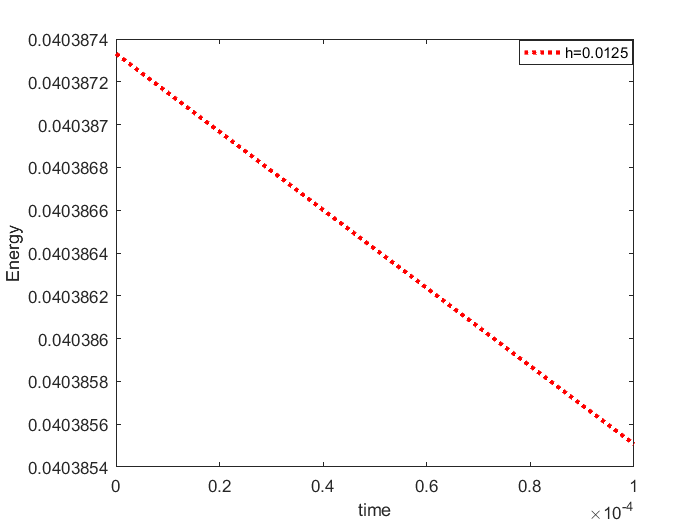}
\end{array}$\vspace{-0.2cm}
\caption{$\mathbf{Example\ \ref{cexmAC2-n}\ (AC)}$, effect of mesh size $h$ on the discrete energy. 
First line: $P_1$ element, Second line: $P_2$ element.
}\label{Cexp8u-hB}
 \end{figure}

 \end{example}

\begin{example}\label{exm3d} \cite{ZFH14} In this example, we apply the BDF2-IEQ-FEM2 scheme \eqref{BDF2-FEM2} with $B=1$, $\Delta t =10^{-4}$ to solve the 3D AC equation (\ref{AC1}) in $\Omega=[-1,1]^{3}$ with parameter $\epsilon=0.05$ and the initial solution
\[u_{0}=\epsilon\cos(1.5\pi x)\cos(1.5\pi y)\big(\sin(\pi z)+\sin(2\pi z)\big).\]

\begin{figure}[!htbp]
$\begin{array}{cccc}
\includegraphics[width=5.5cm,height=4.5cm]{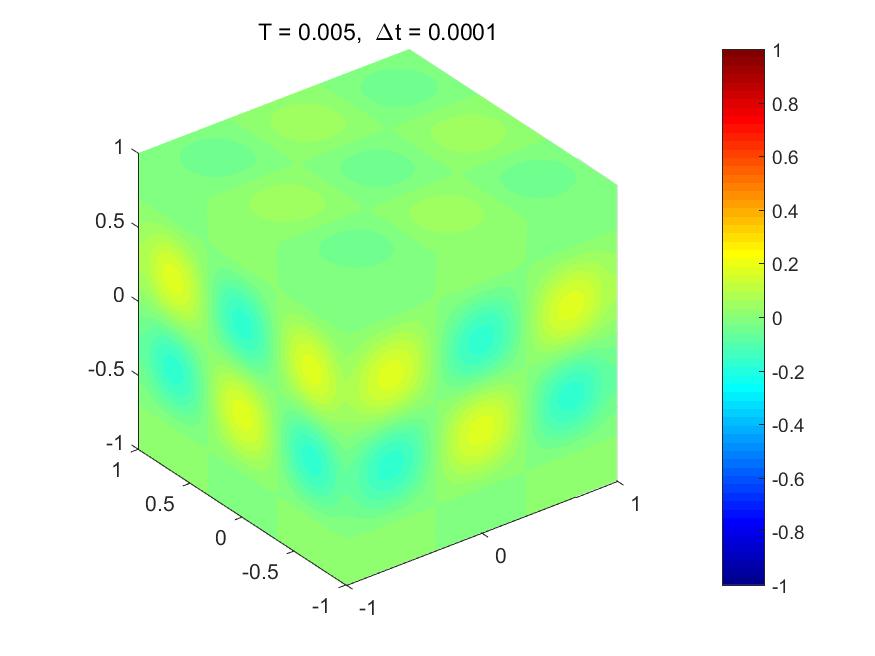}&
\includegraphics[width=5.5cm,height=4.5cm]{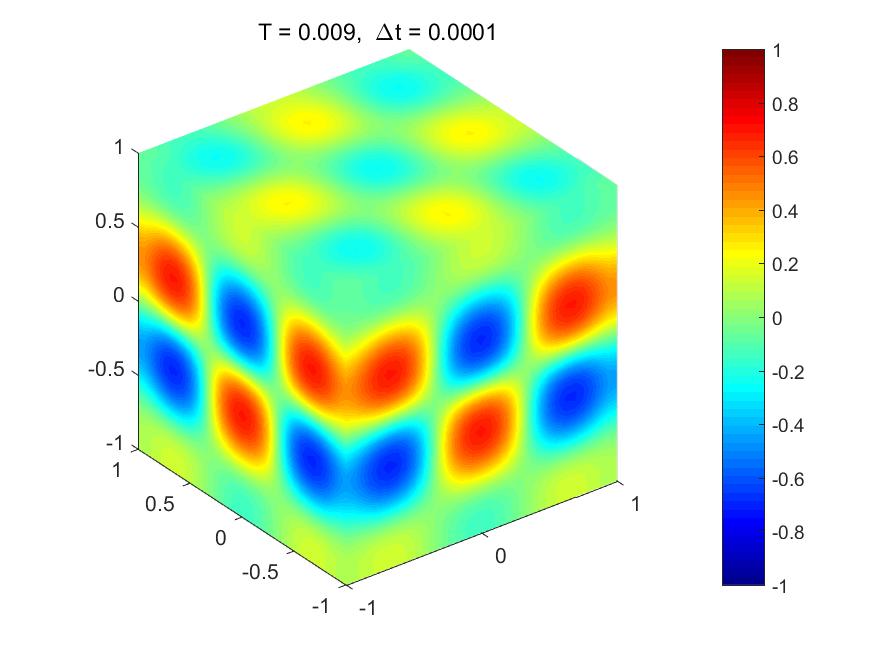}&
\includegraphics[width=5.5cm,height=4.5cm]{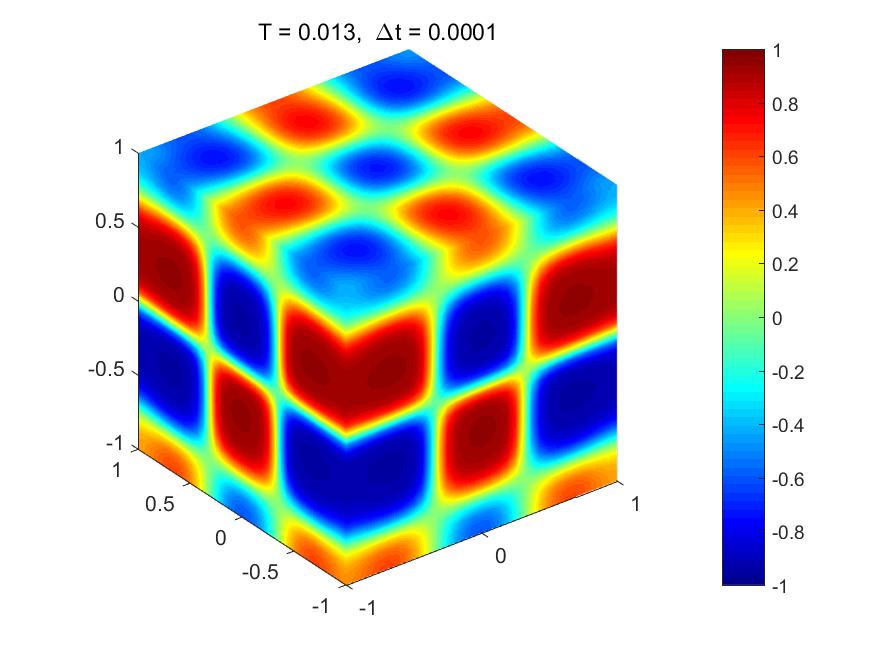}\\
\includegraphics[width=5.5cm,height=4.5cm]{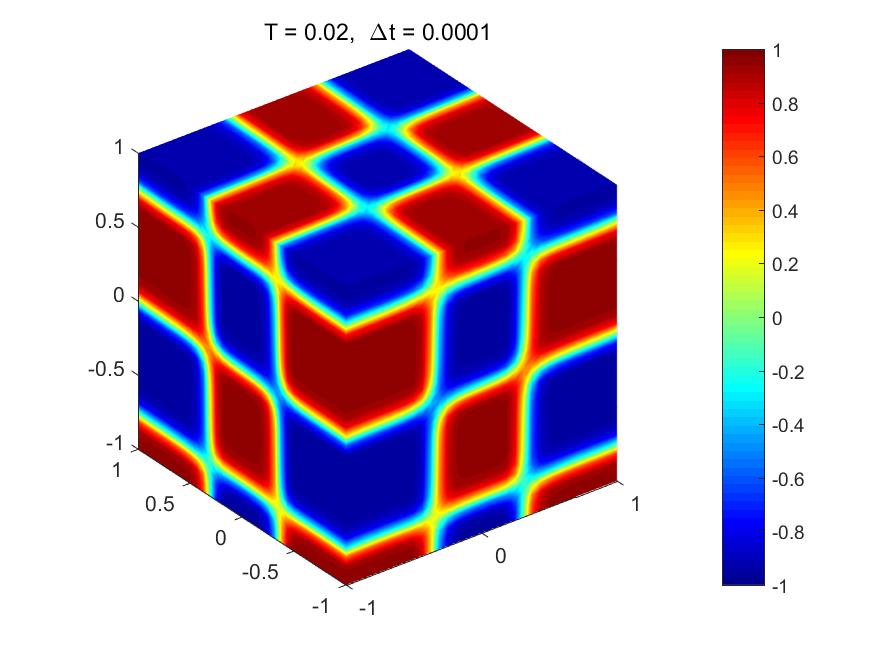}&
\includegraphics[width=5.5cm,height=4.5cm]{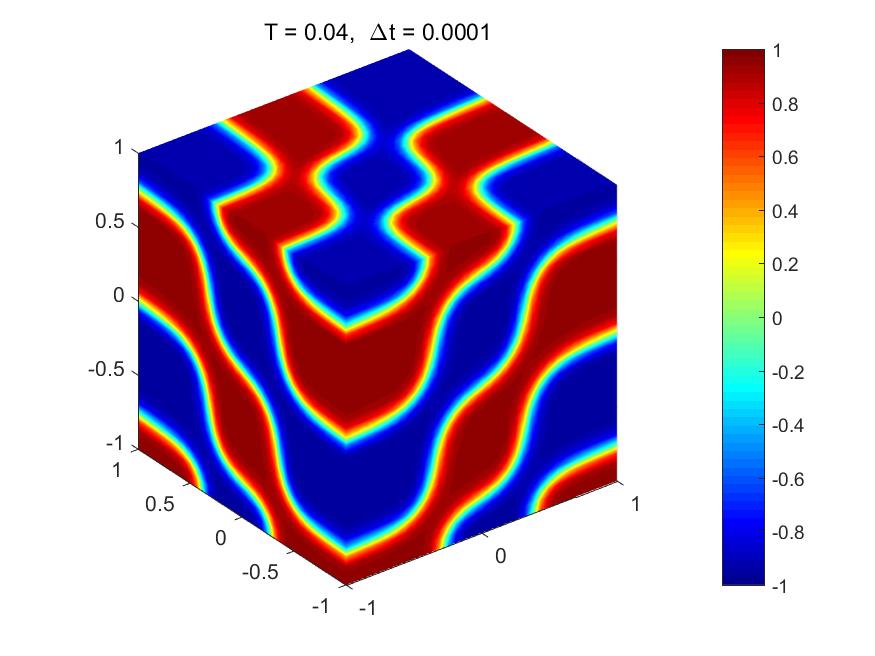}&
\includegraphics[width=5.5cm,height=4.5cm]{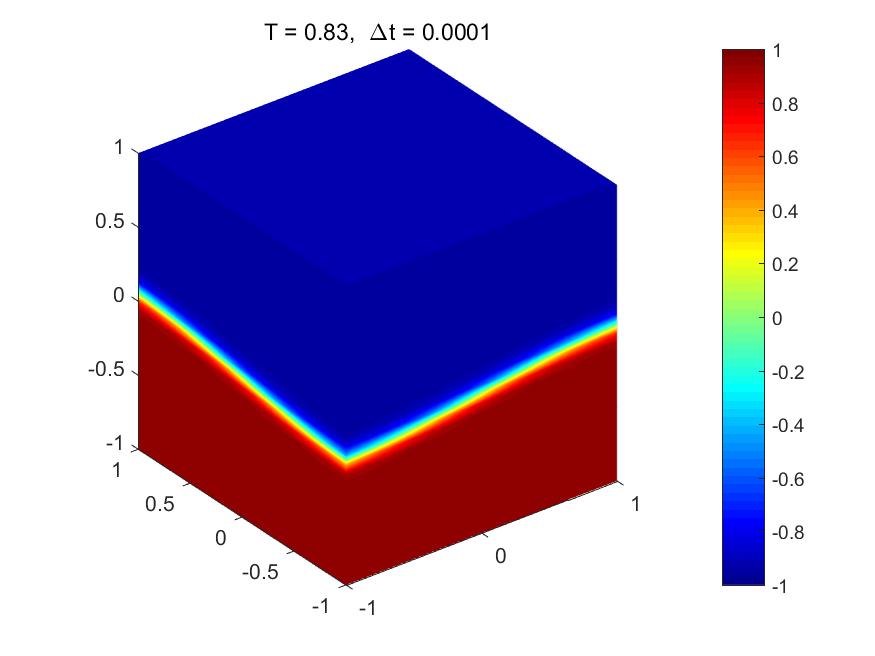}&
\end{array}$\vspace{-0.2cm}
\caption{$\mathbf{Example\ \ref{exm3d}\ (AC3D)}$, BDF2-IEQ-FEM3 scheme (\ref{BDF2-FEM3}), snapshots of numerical solutions.}\label{exp3dAC}
\end{figure}

Figure \ref{exp3dAC} displays the contour plots illustrating the corresponding approximate solutions. The observed coarsening phenomena confirm that the proposed method works effectively for solving the 3D AC equation.
\end{example}
\begin{example}\label{exmS3CN} In the last example, we consider the AC equation $(\ref{AC})$ in $\Omega=[-0.5,0.5]^{2}$ with the logarithmic Flory-Huggins potential
$$F(u)=600(u\ln u+(1-u)\ln(1-u))+1800u(1-u),$$
and the parameter $\epsilon=1$. 
We take the initial data as the $L^2$ projection of the following function
The equation is subject to the following initial data
\bq\label{exmS3e12}
u_{0}=\left\{\begin{aligned}
&0.71, \qquad |x|\leq 0.2, \ |y|\leq 0.2,\\
&0.69, \qquad \text{otherwise}.
\end{aligned}\right.
\eq

\begin{figure}[!htbp]
$\begin{array}{cccc}
 \includegraphics[width=5.5cm,height=4.5cm]{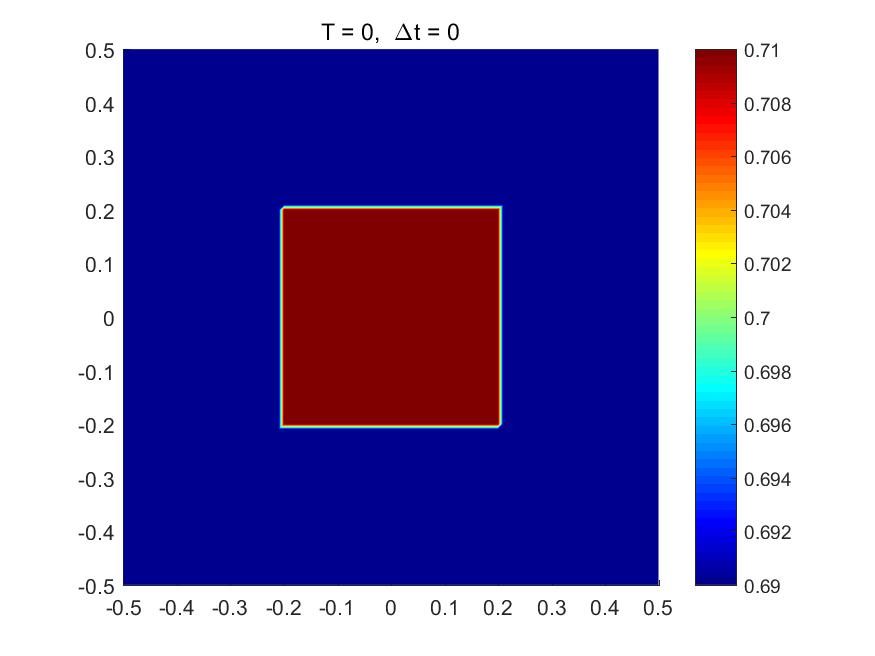}&\vspace{-0.1cm}
 \includegraphics[width=5.5cm,height=4.5cm]{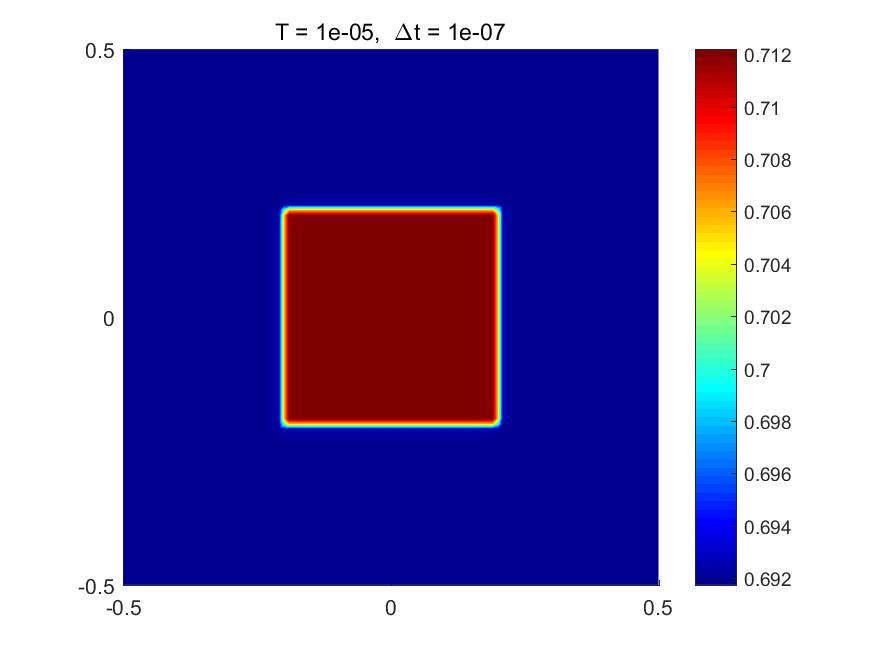}&\vspace{-0.1cm}
 \includegraphics[width=5.5cm,height=4.5cm]{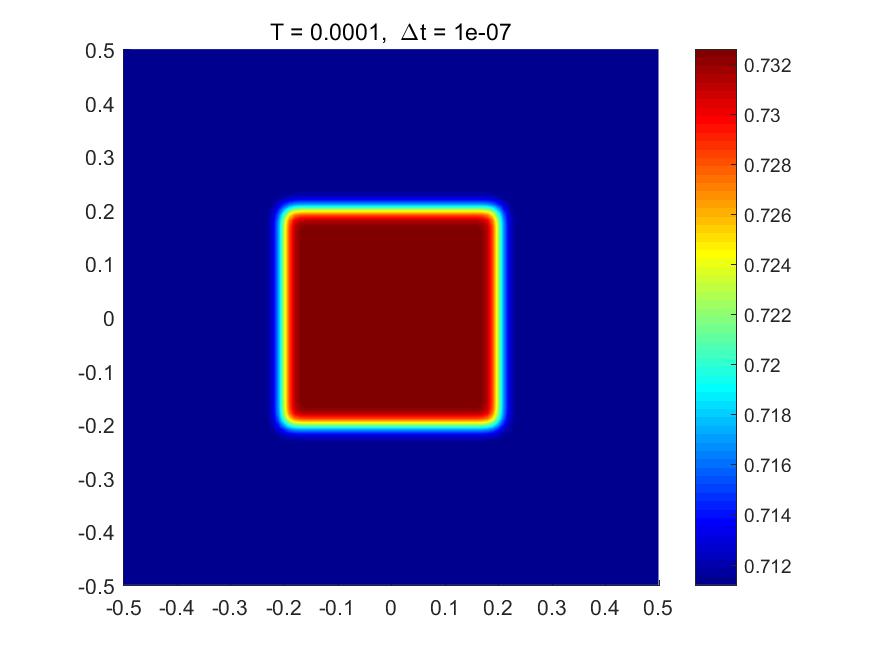}\\
\includegraphics[width=5.5cm,height=4.5cm]{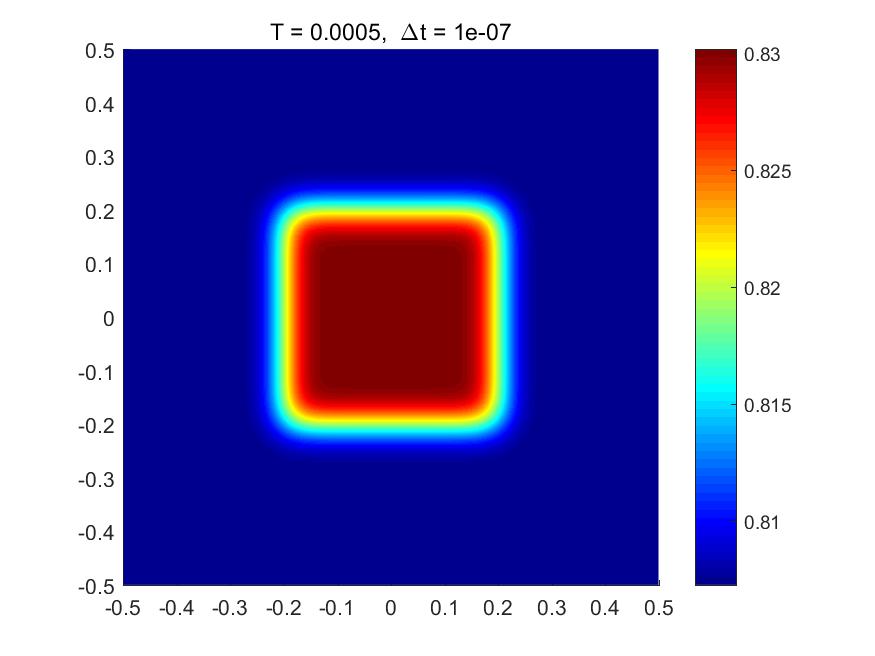}&\vspace{-0.1cm}
\includegraphics[width=5.5cm,height=4.5cm]{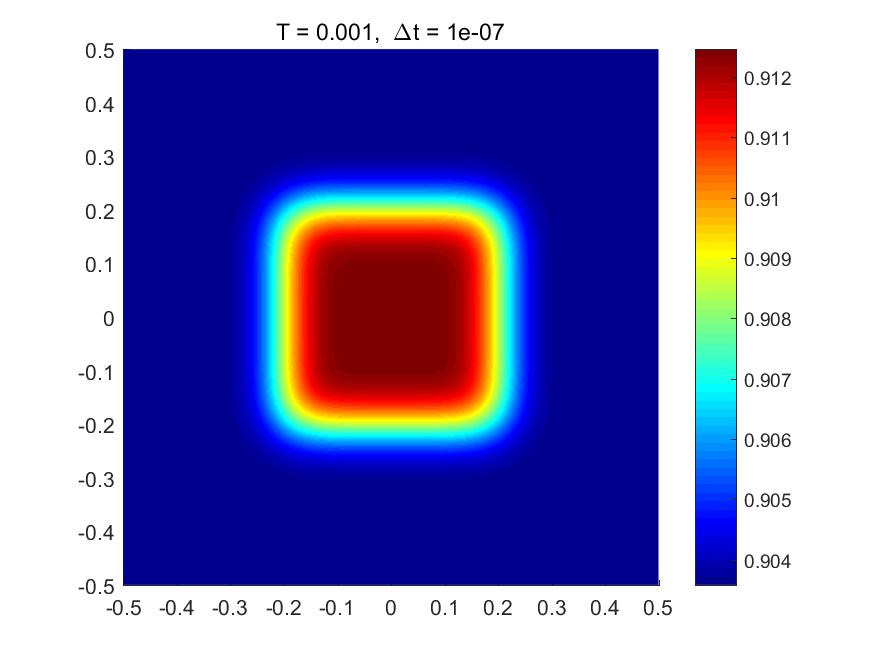}&\vspace{-0.1cm}
\includegraphics[width=5.5cm,height=4.5cm]{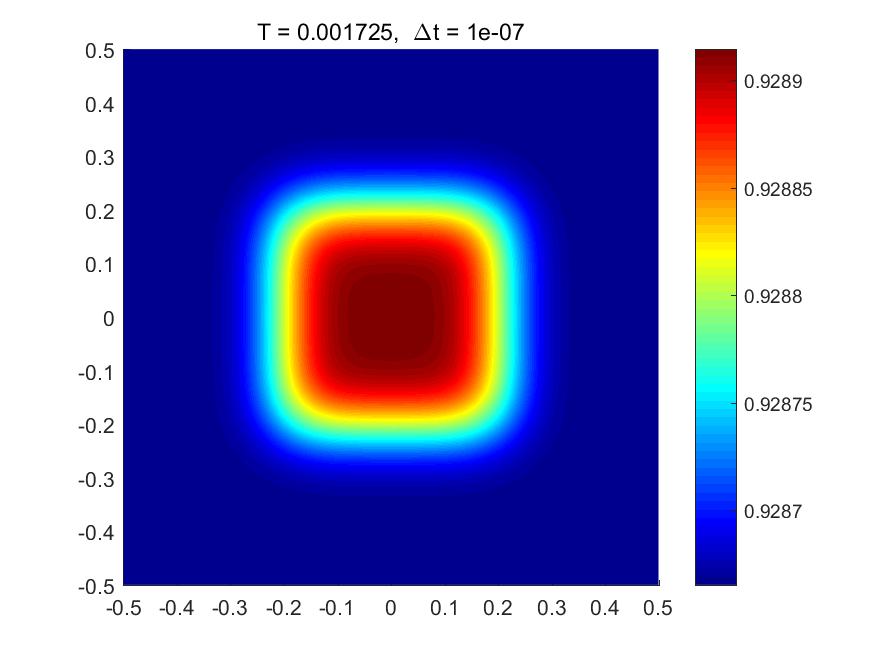}&\vspace{-0.1cm}
\end{array}$\vspace{-0.3cm}
\caption{$\mathbf{Example\ \ref{exmS3CN}\ (AC)}$, CN-IEQ-FEM3 scheme (\ref{CN-FEM3}), snapshots of numerical solutions.}\label{expsuCN-IEQ-FEM3}
\end{figure}

We present contour plots of the corresponding approximate solutions generated by the CN-IEQ-FEM3 scheme (\ref{CN-FEM3}) with $\Delta t=10^{-7}$, $B=100$ in \Cref{expsuCN-IEQ-FEM3}. 
It is shown that the CN-IEQ-FEM3 scheme performs effectively.

\begin{figure}[!htbp]
$\begin{array}{ccc}
\includegraphics[width=5.5cm,height=4.5cm]{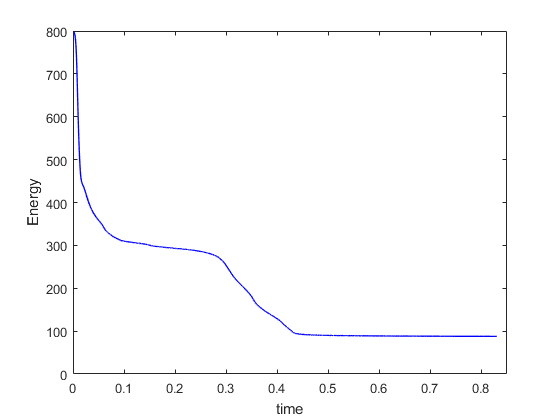}&
\includegraphics[width=5.5cm,height=4.5cm]{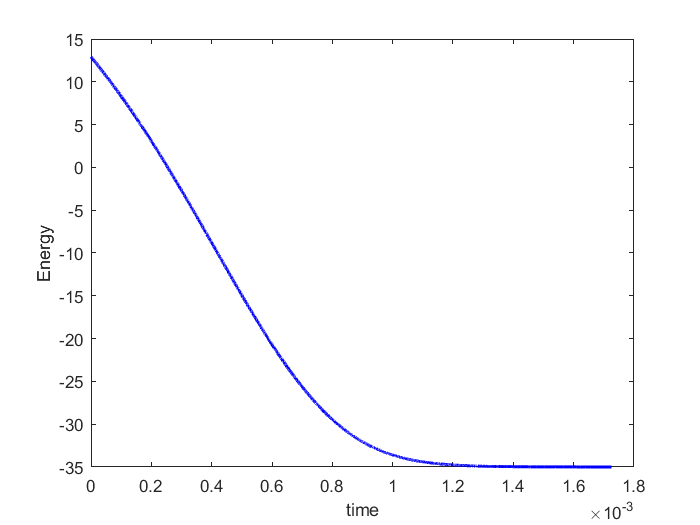}
\end{array}$\vspace{-0.2cm}
\caption{Time evolution of the discrete energy.
Left: $\mathbf{Example\ \ref{exm3d}}$; Right: $\mathbf{Example\ \ref{exmS3CN}}$.}\label{energuAC1}
\end{figure}

In the end, the pictures of discrete energy for $\mathbf{Example\ \ref{exm3d}}$-$\mathbf{Example\ \ref{exmS3CN}}$ are shown in Figure \ref{energuAC1}. From the pictures, it can be observed that the discrete energy decreases with time.
\end{example}

\section{Conclusion}\label{secConsi}

In this paper, we have proposed three types of linearly, first- and second-order, unconditionally energy-stable IEQ-FEMs for solving the CH equation and the AC equation. 
In addition, the methods are mass conservative for the CH equation.
These three types of IEQ-FEMs position the intermediate function introduced by the IEQ approach in different function spaces.
Method 1 places the intermediate function solely in the finite element space, yet it comes with high computational costs. 
Method 2 places the intermediate function in the continuous function space, resulting in the most cost-effective computation. 
Though the derived energy is unconditional energy stable, its approximated energy in the finite element space is only conditionally achievable. 
Method 3 first computes the intermediate function in the continuous function space and then projects it onto the finite element space.  This method strikes a balance, offering a cost-effective computation while ensuring unconditional energy decay in the finite element space. 
Numerical examples have been conducted to assess the accuracy and efficiency, as well as to verify the energy dissipation properties for both equations and mass conservation for the CH equation. 

The error analysis of the proposed methods, particularly Method 3, is highly desirable and deferred in future work.  Also, motivated by \cite{CHYY24}, exploring the adaptive IEQ-FEMs is intriguing, particularly in addressing the challenges posed by the smallness of the parameter $\varepsilon$ in gradient flow models.
Finally, the proposed IEQ-FEMs may be able to combine with much higher-order time discretization, such as the Runge–Kutta methods, but it would be interesting to explore their numerical stability.



\appendix\label{secAppen}

\section{Second order IEQ-FEMs for the CH equation}\label{CH2nd}

In this section we present IEQ-FEMs based on the second-order backward time-differentiation formula (BDF2), demonstrating numerical stability for both the double well and the  logarithmic Flory–Huggins potentials.

In contrast, the Crank–Nicolson (CN) time discretization exhibits high-order accuracy for the CH equation with a double well potential, but instability occurs 
with the logarithmic Flory–Huggins potential. Please note that this phenomenon has also been observed for the IEQ-DG schemes \cite{LY21, Y19}.

\subsection{Second order fully discrete BDF2-IEQ-FEMs}

In this section, we present three stable fully discrete second order backward differentiation formula (BDF2) IEQ-FEMs for the CH equation \eqref{CH} by approximating the intermediate function $U(x,t)$ in different function spaces.
For these schemes, an explicit second order approximation of the numerical solution at $t^{n+1}$ using $u_h^{n-1}$ and $u^n_h$ is given by
\begin{align}\label{CHBDF2v}
u_h^{n, *}=& 2u_h^n-u_h^{n-1}.
\end{align}

\subsubsection{Method 1 (BDF2-IEQ-FEM1)}
First, we use $U_h(x,t) \in V_h$ in the finite element space to approximate $U$ based on \eqref{CHSemiFEM+}.
Given $u_h^{n-1}$, $u_h^n$, $U_h^{n-1}$, $U_h^n \in V_h$ for $n\geq 1$, the BDF2-IEQ-FEM1 scheme is to find $(u^{n+1}_h, w_h^{n+1}, U_h^{n+1}) \in V_h  \times V_h\times V_h$ such that 
\begin{subequations}\label{BDF2-IEQ-FEM1}
	\begin{align}
	\left(  \frac{3u_h^{n+1} - 4u_h^n+u_h^{n-1}}{2\Delta t}, \phi \right)
	= & - A(M(u^{n, *}_h);w_h^{n+1}, \phi),\qquad \forall \phi
\in V_h,\\
	(w_h^{n+1}, \psi) = & A(\epsilon^2;u_h^{n+1},\psi) + \left(H(u^{n,*}_h)U_h^{n+1},\psi \right),\qquad \forall \psi
\in V_h,\\
  \left(\frac{3U_h^{n+1} - 4U_h^n + U_h^{n-1}}{2\Delta t},\tau \right)  =&  \left( \frac{1}{2}H(u^{n,*}_h) \frac{3u_h^{n+1} - 4u_h^n+u_h^{n-1}}{2\Delta t},\tau\right)\qquad \forall \tau
\in V_h.
	\end{align}
\end{subequations}
Here, the initial data $(u_h^0,U_h^0)$ is given by \eqref{u0init} and \eqref{CHinit}, and 
\be
u_h^{-1} = u_h^0, \quad U_h^{-1} = U_h^0.
\ee

For scheme \eqref{BDF2-IEQ-FEM1}, the following results hold. 
\begin{lem}
The second order fully discrete BDF2-IEQ-FEM1 scheme (\ref{BDF2-IEQ-FEM1}) admits a unique solution $(u^{n+1}_h, w_h^{n+1}, U_h^{n+1}) \in V_h  \times V_h \times V_h$ for any $\Delta t>0$. The solution $u_h^n$ satisfies the total mass conservation
\eqref{BDF1MassCons} for any $n>0$. Additionally, it satisfies the energy dissipation law
\bqs
\bar{E}^{n+1}_h = \bar E^n_h  -\Delta t A(M(u^{n, *}_h); w_h^{n+1},  w_h^{n+1})-\frac{1}{4}A(\epsilon^2; u_h^{n+1}-u_h^{n,*},u_h^{n+1}-u_h^{n,*})-\frac{1}{2}\|U_h^{n+1}-U_h^{n,*}\|^2,
\eqs
independent of time step $\Delta t$, where
$$
\bar E^n_h = \frac{1}{2} ( E(u_h^n, U_h^n) + E(u_h^{n,*}, U_h^{n,*})  ),
$$
with $U_h^{n,*}=2U_h^n - U_h^{n-1}$ and $u_h^{n,*}$ is given in (\ref{CHBDF2v}). 
\end{lem}

Similar to the BDF1-IEQ-FEM1 scheme \eqref{BDF1-IEQ-FEM1} 
, the scheme \eqref{BDF2-IEQ-FEM1} can either be solved by the enlarged linear system with the three unknowns, or by a simplified linear system with computing the inverse of the mass matrix.

\subsubsection{Method 2 (BDF2-IEQ-FEM2)}
To save the computational cost, we move to retain the intermediate function $U(x,t)$ in continuous function space $C^0(\Omega)$ with the evolution determined by \eqref{uode}.
Given $u_h^{n-1}$, $u_h^n$, $U^{n-1}$, $U^n$ for $n\geq 1$, the BDF2-IEQ-FEM2 scheme is to find $(u^{n+1}_h, w_h^{n+1}) \in V_h  \times V_h$ and $U^{n+1} \in C^0(\Omega)$ such that  
\begin{subequations}\label{BDF2-IEQ-FEM2}
	\begin{align}
	\left(  \frac{3u_h^{n+1} - 4u_h^n+u_h^{n-1}}{2\Delta t}, \phi \right)
	= & - A(M(u^{n, *}_h);w_h^{n+1}, \phi),\qquad \forall \phi
\in V_h,\\
	(w_h^{n+1}, \psi) = & A(\epsilon^2;u_h^{n+1},\psi) + \left(H(u^{n,*}_h)U^{n+1},\psi \right),\qquad \forall \psi
\in V_h,\\
  \frac{3U^{n+1} - 4U^n + U^{n-1}}{2\Delta t}  =&  \frac{1}{2}H(u^{n,*}_h) \frac{3u_h^{n+1} - 4u_h^n+u_h^{n-1}}{2\Delta t}.
	\end{align}
\end{subequations}
Here, the initial data $(u_h^0,U^0)$ is given by \eqref{u0init} and \eqref{CHinit+}, and 
\be
u_h^{-1} = u_h^0, \quad U^{-1} = U^0.
\ee

By (\ref{BDF2-IEQ-FEM2}c), we have
\bq\label{CHBDF2U}
      U^{n+1}=\frac{1}{2}H(u^{n,*}_h)u_h^{n+1}+\left( \frac{4U^n-U^{n-1}}{3}- \frac{1}{2}H(u^{n,*}_h)\frac{4u_h^n-u_h^{n-1}}{3} \right).
\eq
Plugging $U^{n+1}$ in \eqref{CHBDF2U} into (\ref{BDF2-IEQ-FEM2}b) gives a linear system in terms of $u_h^{n+1}$, $w_h^{n+1}$,
\bq\label{CHBDF2Lin}
\ba
& 	\left(  \frac{3u_h^{n+1} - 4u_h^n+u_h^{n-1}}{2\Delta t}, \phi \right)
	= - A(M(u^{n, *}_h);w_h^{n+1}, \phi),\\
	& (w_h^{n+1}, \psi) =  A(\epsilon^2;u_h^{n+1},\psi) +  \frac{1}{2}\left((H(u^{n,*}_h))^2u_h^{n+1},\psi \right) \\
	& + \frac{1}{3}\left(H(u^{n,*}_h)(4U^{n}-U^{n-1}),\psi \right) - \frac{1}{6}\left((H(u^{n,*}_h))^2(4u_h^{n}-u_h^{n-1}),\psi \right).
\ea
\eq

The scheme (\ref{BDF2-IEQ-FEM2}) is equivalent to the system (\ref{CHBDF2Lin}) with (\ref{CHBDF2U}). The system (\ref{CHBDF2Lin}) is a linear system in terms of $(u_h^{n+1},w_h^{n+1})$, so that the scheme (\ref{BDF2-IEQ-FEM2}) can also avoid solving a linear system coupled with unknown $U^{n+1}$.
Upon solving $u_h^{n+1}$, $w_h^{n+1}$ from (\ref{CHBDF2Lin}), we can obtain $U^{n+1}$ by (\ref{CHBDF2U}), or (\ref{BDF2-IEQ-FEM2}c). 



Then we have the following result.

\begin{thm}\label{CHBDF2thm2}
The second order fully discrete BDF2-IEQ-FEM2 scheme (\ref{BDF2-IEQ-FEM2}) admits a unique solution $(u_h^{n+1},w_h^{n+1})\in V_h$ and $U^{n+1}$ for any $\Delta t>0$, and the solution $u_h^n$ satisfies the total mass conservation \eqref{BDF1MassCons} for any $n>0$, and the energy dissipation law
\bq\label{CHBDF2Stab2}
\bar{E}^{n+1} = \bar E^n  -\Delta t A(M(u^{n, *}_h); w_h^{n+1},  w_h^{n+1})-\frac{1}{4}A(\epsilon^2; u_h^{n+1}-u_h^{n,*},u_h^{n+1}-u_h^{n,*})-\frac{1}{2}\|U^{n+1}-U^{n,*}\|^2,
\eq
independent of time step $\Delta t$, where
$$
\bar E^n = \frac{1}{2} ( E(u_h^n, U^n) + E(u_h^{n,*}, U^{n,*})  ),
$$
with $U^{n,*}=2U^n - U^{n-1}$ and $u_h^{n,*}$ given in (\ref{CHBDF2v}). 
\end{thm}

Note that $\bar E^n$ in \Cref{CHBDF2thm2} can not be computed exactly and it can only approximated by 
\be\label{BDF2Ebar}
\bar E_h^n = \frac{1}{2} ( E(u_h^n, \Pi U^n) + E(u_h^{n,*}, \Pi U^{n,*})  ),
\ee
with $U^{n,*}=2U^n - U^{n-1}$ and $u_h^{n,*}$ given in (\ref{CHBDF2v}). 
Numerically, the computed energy $\bar E_h^n$ may not be dissipating. To guarantee the computed energy dissipation, we further introduce another IEQ-FEM scheme.

\subsubsection{Method 3 (BDF2-IEQ-FEM3)}
To modify the BDF2-IEQ-FEM2 scheme \eqref{BDF2-IEQ-FEM2} such that the computed numerical solution satisfies the energy dissipation law, we provide the BDF2-IEQ-FEM3 scheme, which is to find $(u^{n+1}_h, w_h^{n+1}) \in V_h  \times V_h$ and $U^{n+1} \in C^0(\Omega)$ such that 
\begin{subequations}\label{BDF2-IEQ-FEM3}
	\begin{align}
	\left(  \frac{3u_h^{n+1} - 4u_h^n+u_h^{n-1}}{2\Delta t}, \phi \right)
	= & - A(M(u^{n, *}_h);w_h^{n+1}, \phi),\qquad  \forall \phi
\in V_h,\\
	(w_h^{n+1}, \psi) = & A(\epsilon^2;u_h^{n+1},\psi) + \left(H(u^{n,*}_h)U^{n+1},\psi \right),\qquad \forall \psi
\in V_h,\\
 U_h^n = & \Pi U^n,\\
  \frac{3U^{n+1} - 4U_h^n + U_h^{n-1}}{2\Delta t}  =&  \frac{1}{2}H(u^{n,*}_h) \frac{3u_h^{n+1} - 4u_h^n+u_h^{n-1}}{2\Delta t}.
	\end{align}
\end{subequations}
Here, the initial data $(u_h^0, U^0)$ is given by \eqref{u0init} and \eqref{CHinit+}, and $$u_n^{-1}=u_h^0, \quad U_h^{-1}=U_h^0.$$

The following result holds for the scheme \eqref{BDF2-IEQ-FEM3}.
\begin{thm}\label{CHBDF2thm3}
The second order fully discrete BDF2-IEQ-FEM3 scheme (\ref{BDF2-IEQ-FEM3}) admits a unique solution $(u_h^{n+1},w_h^{n+1})\in V_h$ and $U^{n+1} \in C^0{\Omega}$ for any $\Delta t>0$. The solution $u_h^n$ satisfies the total mass conservation \eqref{BDF1MassCons} for any $n>0$. Additionally, it satisfies the energy dissipation law
\bq\label{CHBDF2Stab3}
\ba
& \bar E(u_h^{n+1}, U_h^{n+1}; u_h^{n+1,*}, U_h^{n+1,*}) \leq 
\bar E(u_h^{n+1}, U^{n+1}; u_h^{n+1,*}, U_h^{n+1,*}) = \bar E(u_h^{n}, U^{n}; u_h^{n,*}, U_h^{n,*}) \\
& -\Delta t A(M(u^{n, *}_h); w_h^{n+1},  w_h^{n+1})-\frac{1}{4}A(\epsilon^2; u_h^{n+1}-u_h^{n,*},u_h^{n+1}-u_h^{n,*})-\frac{1}{2}\|U^{n+1}-U^{n,*}\|^2,
\ea
\eq
independent of time step $\Delta t$, where
$$
\bar E(v^n, V^n; v^{n,*}, V^{n,*}) = \frac{1}{2} ( E(v^n, V^n) + E(v^{n,*}, V^{n,*})  ),
$$
with $V^{n,*}=2V^n - V^{n-1}$ for $V=U, U_h$ and $u_h^{n,*}$ being given in (\ref{CHBDF2v}). 
\end{thm}

Compared with the BDF2-IEQ-FEM2 scheme \eqref{BDF2-IEQ-FEM2}, the BDF2-IEQ-FEM3 scheme \eqref{BDF2-IEQ-FEM3} includes an additional projection step (\ref{BDF2-IEQ-FEM3}c). This extra computational cost is traded for the numerical satisfaction of the energy dissipation law.

\section{Second order IEQ-FEMs for the AC equation}\label{AC2nd}

In this section, we present three stable fully discrete CN-IEQ-FEMs and BDF2-IEQ-FEMs for the AC equation \eqref{AC}.

\subsection{Second order fully discrete CN-IEQ-FEMs}

We discretize the semi-discrete IEQ-FEM schemes (\ref{ACSemiFEM}) and \eqref{ACSemiFEM+} by the second order Crank–Nicolson (CN) time discretization.

We denote
$$
v^{n+1/2}=(v^n+v^{n+1})/2,
$$
where $v$ is a given function.
For these schemes, an explicit second order approximation of the numerical solution at $t^{n+1/2}$ using $u_h^{n-1}$ and $u^n_h$ is given by
\begin{align}\label{CNu8}
u_h^{n, *}=& \frac{3}{2}u_h^n-\frac{1}{2}u_h^{n-1}.
\end{align}
where $u_h^n=u_h(x,t^{n})$.

We first introduce the trivial IEQ-FEM schemes based on the semi-discrete IEQ-FEM schemes \eqref{ACSemiFEM} and \eqref{ACSemiFEM+}.
\subsubsection{Method 1 (CN-IEQ-FEM1)}
Given $u_h^{n}$ and $U_h^n$, the second order fully discrete CN-IEQ-FEM1 scheme based on the semi-discrete IEQ-FEM scheme \eqref{ACSemiFEM} is to find $(u^{n+1}_h, U_h^{n+1})\in V_h \times V_h$ such that  
\begin{subequations}\label{CN-FEM1}
	\begin{align}
	\left(  \frac{u_h^{n+1/2} - u_h^n}{\frac{1}{2}\Delta t}, \phi \right) = & -  A(\epsilon^2; u_h^{n+1/2},\phi) - \left(H(u^{n,*}_h)U_h^{n+1/2},\phi \right),\qquad \forall \phi
\in V_h,\\
 \left( \frac{U_h^{n+1/2} - U_h^n}{\frac{1}{2}\Delta t}, \psi\right)  =&  \left(\frac{1}{2} H(u^{n,*}_h) \frac{u_h^{n+1/2} - u_h^n}{\frac{1}{2}\Delta t}, \psi\right),\qquad \forall \psi
\in V_h.
	\end{align}
\end{subequations}
The initial data is given by (\ref{ACinit}). For the CN-IEQ-FEM1 scheme \eqref{CN-FEM1}, the following result holds.
\begin{lem}
Given $u_h^n \in V_h$ and $U^{n}\in V_h$, the second order fully discrete CN-IEQ-FEM1 scheme (\ref{CN-FEM1}) admits a unique solution $u_h^{n+1}\in V_h$ and $U_h^{n+1}\in V_h$ satisfying the energy dissipation law,
\bq\label{ACCNStab1}
\ba
E(u_h^{n+1}, U_h^{n+1}) = E(u_h^{n}, U_h^{n}) - \frac{1}{\Delta t}\|u_h^{n+1} - u_h^n\|^2,
\ea
\eq
for any $\Delta t>0$.
\end{lem}

\subsubsection{Method 2 (CN-IEQ-FEM2)}

Next, given $u_h^{n} \in V_h$ and $U^n \in C^0(\Omega)$, the second order CN-IEQ-FEM2 scheme based on the semi-discrete IEQ-FEM scheme \eqref{ACSemiFEM+} is to find $u^{n+1}_h\in V_h$ and $U^{n+1} \in C^0(\Omega)$ such that 
\begin{subequations}\label{CN-FEM2}
	\begin{align}
	\left(  \frac{u_h^{n+1/2} - u_h^n}{\frac{1}{2}\Delta t}, \phi \right) = & -  A(\epsilon^2; u_h^{n+1/2},\phi) - \left(H(u^{n,*}_h)U^{n+1/2},\phi \right),\qquad \forall \phi
\in V_h,\\
  \frac{U^{n+1/2} - U^n}{\frac{1}{2}\Delta t}  =&  \frac{1}{2} H(u^{n,*}_h) \frac{u_h^{n+1/2} - u_h^n}{\frac{1}{2}\Delta t}.
	\end{align}
\end{subequations}
The initial data is given by (\ref{ACinit+}). For the CN-IEQ-FEM2 scheme \eqref{CN-FEM2}, it holds the following result.
\begin{thm}
Given $u_h^n \in V_h$ and $U^{n}\in C^0(\Omega)$, the second order fully discrete CN-IEQ-FEM2 scheme (\ref{CN-FEM2}) admits a unique solution $u_h^{n+1}\in V_h$ and $U^{n+1} \in C^0(\Omega)$ satisfying the energy dissipation law,
\bq\label{ACCNStab2}
\ba
E(u_h^{n+1}, U^{n+1}) = E(u_h^{n}, U^{n}) - \frac{1}{\Delta t}\|u_h^{n+1} - u_h^n\|^2,
\ea
\eq
for any $\Delta t>0$.
\end{thm}

\subsubsection{Method 3 (CN-IEQ-FEM3)}  

In this part, we present an IEQ-FEM scheme characterized by low computational cost, ensuring that the resulting numerical solution satisfies the energy dissipation law.

Given $u_h^{n} \in V_h$ and $U^n \in C^0(\Omega)$, the second order CN-IEQ-FEM3 scheme is to find $u^{n+1}_h\in V_h$ and $U^{n+1} \in C^0(\Omega)$ such that  
\begin{subequations}\label{CN-FEM3}
	\begin{align}
	\left(  \frac{u_h^{n+1/2} - u_h^n}{\frac{1}{2}\Delta t}, \phi \right) = & -  A(\epsilon^2; u_h^{n+1/2},\phi) - \left(H(u^{n,*}_h)\bar U^{n+1/2},\phi \right),\qquad \forall \phi
\in V_h,\\
 U_h^n = & \Pi U^n,\\
  \frac{\bar U^{n+1/2} - U_h^n}{\frac{1}{2}\Delta t}  =&  \frac{1}{2} H(u^{n,*}_h) \frac{u_h^{n+1/2} - u_h^n}{\frac{1}{2}\Delta t},
	\end{align}
\end{subequations}
where 
\be
\bar U^{n+1/2} = \frac{1}{2}(U^{n+1}+ U_h^{n}).
\ee
The initial data is given by (\ref{ACinit+}).

By (\ref{CN-FEM3}b), we have 
\begin{align}\label{ACCNU}
\bar U^{n+1/2} =  U_h^n + \frac{1}{2} H(u^{n,*}_h) (u_h^{n+1/2} - u_h^n).
\end{align}
Plugging $\bar U^{n+1/2}$ in (\ref{ACCNU}) into (\ref{CN-FEM3}a), we obtain a linear system in terms of $u_h^{n+1}$,
\bq\label{ACCNLin}
\ba
\left(  \frac{u_h^{n+1/2} - u_h^n}{\frac{1}{2}\Delta t}, \phi \right) = & -  A(\epsilon^2; u_h^{n+1/2},\phi) - \left(H(u^{n,*}_h)U^{n},\phi \right) \\
& - \frac{1}{2}((H(u^{n,*}_h))^2u_h^{n+1/2}, \phi)+\frac{1}{2}((H(u^{n,*}_h))^2u_h^{n}, \phi).
\ea
\eq
Here (\ref{CN-FEM3}) is equivalent to (\ref{ACCNLin}) and (\ref{ACCNU}). Therefore, the scheme (\ref{CN-FEM3}) can also avoid solving the linear system with double unknowns.
Upon solving $u_h^{n+1}$ from (\ref{ACCNLin}), we can obtain $U^{n+1}$ by (\ref{ACCNU}), or (\ref{CN-FEM3}b). 

\begin{rem}
For the CN scheme (\ref{CN-FEM3}), we can plug the expression of $u_h^{n+1/2}$ and $\bar U^{n+1/2}$ into (\ref{CN-FEM3}), then the scheme is in terms of $u_h^{n+1}$ and $U^{n+1}$. In fact, we can also take $\bar u^{n+1/2}$ and $U^{n+1/2}$ as unknowns, and then recover $u_h^{n+1}$ and $U^{n+1}$ by
\bqs
\ba
u_h^{n+1} = & 2u_h^{n+1/2}-u_n^n,\\
U^{n+1} = & 2\bar U^{n+1/2}-U_h^n.
\ea
\eqs
\end{rem}

\begin{thm}\label{ACCNThm}
Given $u_h^n \in V_h$ and $U^{n}\in C^0(\Omega)$, the second order fully discrete CN-IEQ-FEM3 scheme (\ref{CN-FEM3}) admits a unique solution $u_h^{n+1}\in V_h$ and $U^{n+1}\in C^0(\Omega)$ satisfying the energy dissipation law,
\bq\label{ACCNStab}
\ba
E(u_h^{n+1}, U_h^{n+1}) \leq E(u_h^{n+1}, U^{n+1}) = E(u_h^{n}, U_h^{n}) - \frac{1}{\Delta t}\|u_h^{n+1} - u_h^n\|^2,
\ea
\eq
for any $\Delta t>0$.
\end{thm}

\begin{rem}
In contrast to the limitations highlighted for the CN-IEQ-FEMs in the context of the CH equation, it is noteworthy that the CN-IEQ-FEM schemes \eqref{CN-FEM1}, \eqref{CN-FEM2}, and \eqref{CN-FEM3} demonstrate stability for the AC equation. This stability holds for cases involving both the double well potential and the logarithmic Flory–Huggins potential.
\end{rem}

\subsection{Second order fully discrete BDF2-IEQ-FEMs}

We also present the fully discrete second-order backward differentiation formula (BDF2) IEQ-FEMs. 
For these schemes, an explicit second-order approximation of the numerical solution at $t^{n+1}$ using $u_h^{n-1}$ and $u^n_h$ is given by
\begin{align}\label{ACBDF2v}
u_h^{n, *}=& 2u_h^n-u_h^{n-1}.
\end{align}

\subsubsection{Method 1 (BDF2-IEQ-FEM1)}

The second order fully discrete coupled BDF2-IEQ-FEM1 scheme based on the semi-discrete \eqref{ACSemiFEM} is to find $u^{n+1}_h \in V_h$ and $U_h^{n+1} \in V_h$ such that 
\begin{subequations}\label{BDF2-FEM1}
	\begin{align}
	\left(  \frac{3u_h^{n+1} - 4u_h^n+u_h^{n-1}}{2\Delta t}, \phi \right)
	= -&  A(\epsilon^2; u_h^{n+1},\phi) - \left(H(u^{n,*}_h)U_h^{n+1},\phi \right),\qquad \forall \phi
\in V_h,\\
 \left( \frac{3U_h^{n+1} - 4U_h^n + U_h^{n-1}}{2\Delta t} , \psi\right) =&  \left(\frac{1}{2}H(u^{n,*}_h) \frac{3u_h^{n+1} - 4u_h^n+u_h^{n-1}}{2\Delta t}, \psi\right),\qquad \forall \psi
\in V_h,
	\end{align}
\end{subequations}
The initial is given by \eqref{ACinit} and
\be
u_h^{-1} = u_h^0, \quad U_h^{-1} = U_h^0.
\ee
Then the following result holds.
\begin{lem}\label{ACBDF2Thm1}
The second order fully discrete BDF2-IEQ-FEM scheme (\ref{BDF2-FEM1}) admits a unique solution $u_h^{n+1}\in V_h$ and $U^{n+1}\in V_h$ satisfying the energy dissipation law,
\bq\label{ACBDF2Stab1}
\ba
\bar{E}_h^{n+1} = \bar E_h^n - \frac{1}{2\Delta t} \|3u_h^{n+1} - 4u_h^n+u_h^{n-1}\|^2 -\frac{1}{4}A(\epsilon^2; u_h^{n+1}-u_h^{n,*},u_h^{n+1}-u_h^{n,*})-\frac{1}{2}\|U_h^{n+1}-U_h^{n,*}\|^2,
\ea
\eq
for any $\Delta t>0$, where
$$
\bar E^n_h = \frac{1}{2} ( E(u_h^n, U_h^n) + E(u_h^{n,*}, U_h^{n,*})  ),
$$
with $U_h^{n,*}=2U_h^n - U_h^{n-1}$ and $u_h^{n,*}$ being given in (\ref{ACBDF2v}). 
\end{lem}

\subsubsection{Method 2 (BDF2-IEQ-FEM2)}

The second order fully discrete BDF2-IEQ-FEM2 scheme is to find $u^{n+1}_h \in V_h$ and $U^{n+1} \in C^0(\Omega)$ such that 
\begin{subequations}\label{BDF2-FEM2}
	\begin{align}
	\left(  \frac{3u_h^{n+1} - 4u_h^n+u_h^{n-1}}{2\Delta t}, \phi \right)
	= -&  A(\epsilon^2; u_h^{n+1},\phi) - \left(H(u^{n,*}_h)U^{n+1},\phi \right),\qquad \forall \phi
\in V_h,\\
  \frac{3U^{n+1} - 4U^n + U^{n-1}}{2\Delta t}  =&  \frac{1}{2}H(u^{n,*}_h) \frac{3u_h^{n+1} - 4u_h^n+u_h^{n-1}}{2\Delta t}.
	\end{align}
\end{subequations}
The initial solution is given by \eqref{ACinit+} and
\be
u_h^{-1} = u_h^0, \quad U^{-1} = U^0.
\ee
\begin{thm}\label{ACBDF2Thm2}
The second order fully discrete BDF2-IEQ-FEM2 scheme (\ref{BDF2-FEM2}) admits a unique solution $u_h^{n+1}\in V_h$ and $U^{n+1}\in C^0(\Omega)$ satisfying the energy dissipation law,
\bq\label{ACBDF2Stab2}
\ba
\bar{E}^{n+1} = \bar E^n - \frac{1}{2\Delta t} \|3u_h^{n+1} - 4u_h^n+u_h^{n-1}\|^2 -\frac{1}{4}A(\epsilon^2; u_h^{n+1}-u_h^{n,*},u_h^{n+1}-u_h^{n,*})-\frac{1}{2}\|U^{n+1}-U^{n,*}\|^2,
\ea
\eq
for any $\Delta t>0$, where
$$
\bar E^n = \frac{1}{2} ( E(u_h^n, U^n) + E(u_h^{n,*}, U^{n,*})  ),
$$
with $U^{n,*}=2U^n - U^{n-1}$ and $u_h^{n,*}$ being given in (\ref{ACBDF2v}). 
\end{thm}

\subsubsection{BDF2-IEQ-FEM3}

The IEQ-FEM scheme presented below is characterized by low computational cost, ensuring that the resulting numerical solution satisfies the energy dissipation law.
More specifically, the second order fully discrete BDF2-IEQ-FEM3 scheme is to find $u^{n+1}_h \in V_h$ and $U^{n+1}\in C^0(\Omega)$ such that 
\begin{subequations}\label{BDF2-FEM3}
	\begin{align}
	\left(  \frac{3u_h^{n+1} - 4u_h^n+u_h^{n-1}}{2\Delta t}, \phi \right)
	= -&  A(\epsilon^2; u_h^{n+1},\phi) - \left(H(u^{n,*}_h)U^{n+1},\phi \right),\qquad \forall \phi
\in V_h,\\
 U_h^n = & \Pi U^n, \\
  \frac{3U^{n+1} - 4U_h^n + U_h^{n-1}}{2\Delta t}  =&  \frac{1}{2}H(u^{n,*}_h) \frac{3u_h^{n+1} - 4u_h^n+u_h^{n-1}}{2\Delta t}.
\end{align}
\end{subequations}
The initial solution is given by \eqref{ACinit+} and
\be
u_h^{-1} = u_h^0, \quad U_h^{-1} = U_h^0.
\ee

By (\ref{BDF2-FEM3}b), we have
\bq\label{ACBDF2U}
      U^{n+1}=\frac{1}{2}H(u^{n,*}_h)u_h^{n+1}+\left( \frac{4U_h^n-U_h^{n-1}}{3}- \frac{1}{2}H(u^{n,*}_h)\frac{4u_h^n-u_h^{n-1}}{3} \right).
\eq
Plugging $U^{n+1}$ in (\ref{ACBDF2U}), it follows that the scheme (\ref{BDF2-FEM3}) gives a linear system in terms of $u_h^{n+1}$,
\bq\label{ACBDF2Lin}
\ba
	& \left(  \frac{3u_h^{n+1} - 4u_h^n+u_h^{n-1}}{2\Delta t}, \phi \right)
	=  -  A(\epsilon^2; u_h^{n+1},\phi) - \frac{1}{2}\left((H(u^{n,*}_h))^2u_h^{n+1},\phi \right) \\
	& - \frac{1}{3}\left(H(u^{n,*}_h)(4U_h^{n}-U_h^{n-1}),\phi \right) + \frac{1}{6}\left((H(u^{n,*}_h))^2(4u_h^{n}-u_h^{n-1}),\phi \right).
\ea
\eq
Here the scheme (\ref{BDF2-FEM3}) is equivalent to the system formed by (\ref{ACBDF2Lin}) and (\ref{ACBDF2U}).
In the new system, (\ref{ACBDF2Lin}) is an equation with respect to the solution $u_h^{n+1}$ only, which means the scheme (\ref{BDF2-FEM3}) can also avoid solving the coupled system of $u_h^{n+1}$ and $U^{n+1}$. 

\begin{thm}\label{ACBDF2Thm}
The second order fully discrete BDF2-IEQ-FEM3 scheme (\ref{BDF2-FEM3}) admits a unique solution $u_h^{n+1}\in V_h$ and $U^{n+1}\in C^0(\Omega)$ satisfying the energy dissipation law,
\bq\label{ACBDF2Stab}
\ba
& \bar E(u_h^{n+1}, U_h^{n+1}; u_h^{n+1,*}, U_h^{n+1,*}) \leq 
\bar E(u_h^{n+1}, U^{n+1}; u_h^{n+1,*}, U_h^{n+1,*}) = \bar E(u_h^{n}, U^{n}; u_h^{n,*}, U_h^{n,*}) \\
& - \frac{1}{2\Delta t} \|3u_h^{n+1} - 4u_h^n+u_h^{n-1}\|^2 -\frac{1}{4}A(\epsilon^2; u_h^{n+1}-u_h^{n,*},u_h^{n+1}-u_h^{n,*})-\frac{1}{2}\|U^{n+1}-U_h^{n,*}\|^2,
\ea
\eq
for any $\Delta t>0$, where
$$
\bar E(v^n, V^n; v^{n,*}, V^{n,*}) = \frac{1}{2} ( E(v^n, V^n) + E(v^{n,*}, V^{n,*})  ),
$$
with $V^{n,*}=2V^n - V^{n-1}$ for $V=U, U_h$  and $u_h^{n,*}$ being given in (\ref{ACBDF2v}). 
\end{thm}

\section*{Acknowledgments}
Chen's research was supported by NSFC Project (12201010), Natural Science Research Project of Higher Education in Anhui Province (2022AH040027). 
Yi's research was supported partially by NSFC Project (12071400) and  China's National Key R\&D Programs  (2020YFA0713500).
Yin’s research was supported by the University of Texas at El Paso Startup Award.

\end{document}